\documentclass[11pt,a4paper,reqno]{amsart}
\usepackage[T1]{fontenc}
\usepackage[utf8]{inputenc}
\usepackage{lmodern}
\usepackage[english]{babel}
\usepackage[margin=1in]{geometry}
\usepackage{amssymb,mathtools,mathrsfs}
\usepackage{microtype}
\usepackage{xcolor}
\usepackage[numbers,square,semicolon]{natbib}
\usepackage[unicode,colorlinks=true,linktoc=all]{hyperref}
\usepackage{bookmark}
\newtheorem{theorem}{Theorem}[section]
\newtheorem{proposition}[theorem]{Proposition}
\newtheorem{corollary}[theorem]{Corollary}
\newtheorem{lemma}[theorem]{Lemma}
\theoremstyle{definition}
\newtheorem{definition}[theorem]{Definition}
\theoremstyle{remark}
\newtheorem{remark}[theorem]{Remark}
\numberwithin{equation}{section}

\makeatletter
\renewcommand{\section}{\@startsection{section}{1}%
  \z@{.7\linespacing\@plus\linespacing}{.5\linespacing}%
  {\normalfont\bfseries\centering\let\@secnumfont\bfseries}}
\renewcommand{\subsection}{\@startsection{subsection}{2}%
  \z@{.5\linespacing\@plus.7\linespacing}{.25\linespacing}%
  {\normalfont\bfseries}}
\makeatother

\allowdisplaybreaks[2]
\date{}

\DisableLigatures{encoding=*,family=*}
\definecolor{PreviewCitationBlue}{RGB}{0,63,137}
\hypersetup{citecolor=blue,urlcolor=PreviewCitationBlue,
  linkcolor=black,pdflang={en-US}}

\providecommand{\ii}{\mathrm{i}}
\providecommand{\ee}{\mathrm{e}}
\providecommand{\rd}{\mathrm{d}}
\title[ON SEMISTABLE HESSIAN TYPE EQUATIONS]{ON SEMISTABLE HESSIAN TYPE EQUATIONS}
\author{Jixiang Fu}
\address{Shanghai Center for Mathematical Sciences,
Fudan University,
Shanghai 200433, China}
\email{majxfu@fudan.edu.cn}

\author{Xuan Li}
\address{Shanghai Institute for Mathematics and Interdisciplinary Sciences (SIMIS), Shanghai 200433, China}
\email{lixuan@simis.cn}

\author{Dekai Zhang}
\address{School of Mathematical Sciences, Key Laboratory of Mathematics
and Engineering Applications (Ministry of Education), Shanghai Key Laboratory
of PMMP, East China Normal University, Shanghai 200241, China}
\email{dkzhang@math.ecnu.edu.cn}

\author{Ziyi Zhang}
\address{School of Mathematical Sciences, Fudan University,
Shanghai 200433, China}
\email{21210180101@m.fudan.edu.cn}
\hypersetup{pdftitle={ON SEMISTABLE HESSIAN TYPE EQUATIONS}}
\begin{document}
\begin{abstract}
{
Given a right-Noetherian polynomial, we consider the semistable $(1,1)$-class with positive polynomial mass on a compact K\"ahler manifold. We establish a uniform relative $L^\infty$ estimate for generalized Hessian type equations with Orlicz weights on the densities. For a smooth positive density, we establish higher regularity of solutions outside the Hessian null locus, and we also obtain a Collins-Tosatti type characterization for this locus. We prove the existence and uniqueness of weak solutions. The equations considered here include the complex Hessian equation, the Hessian quotient equation, and the LYZ equation.}
\end{abstract}
\maketitle
\markboth{\shorttitle}{\shorttitle}
\tableofcontents

\section{Introduction}\label{sec:introduction}

Yau~\cite{Y78} proved Calabi's conjecture by solving the complex
Monge-Ampère equation on compact Kähler manifolds. For the complex Hessian equation on compact Kähler manifolds,
Hou-Ma-Wu~\cite{HMW2010} established a second-order estimate,
Dinew-Ko{\l}odziej~\cite{DK2017} then proved a gradient estimate,
completing the proof of smooth solvability.
 For the $J$-equation, inverse Hessian and Hessian quotient
equations, and the supercritical LYZ equation,
solvability follows from associated cone condition or $\mathcal C$-subsolution
conditions~\cite{SongWeinkove2008,FangLaiMa2011,szekelyhidi2018,CollinsJacobYau2020}.

\begingroup\emergencystretch=4em
A natural question is whether these conditions can be characterized numerically. Demailly-P{\u a}un~\cite{DP04} gave a numerical characterization of the Kähler cone. For the $J$-equation, Lejmi-Székelyhidi~\cite{LejmiSzekelyhidi2015} conjectured a numerical criterion for solvability in terms of intersections with subvarieties.
Chen~\cite{Chen2021J} proved a uniform version of this criterion and
characterized solvability of the supercritical LYZ equation by uniform
intersection inequalities along test families. 
Song~\cite{Song2020}  later proved the conjecture of ~\cite{LejmiSzekelyhidi2015}. 
Further numerical criteria were established for generalized Monge-Ampère
equations on projective manifolds~\cite{DatarPingali2021}, for equations
with differential-form coefficients~\cite{FangMa2024}, and for the
$\sigma_2$ and $\sigma_2/\sigma_1$ equations on compact Kähler
threefolds~\cite{FuZhangZhang2026}.
\par\endgroup

{
Lin introduced right-Noetherian polynomials in his study of the convexity
of general inverse $\sigma_k$ equations~\cite{Lin2023Convexity}. He
subsequently proved solvability for constant-coefficient strictly
$\Upsilon$-stable equations admitting a $\mathcal C$-subsolution
\cite{Lin2023Solvability}. Fang-Ma developed the theory of ideal
G{\aa}rding polynomials and studied the associated elliptic
equations~\cite{fangma2026,fangmaspectral2026}. On projective manifolds,
Chen-Nie-Xu~\cite{CNX26} proved a numerical criterion for
strongly strictly right-Noetherian polynomials of degree equal to the
dimension, and a uniform criterion for strictly right-Noetherian
polynomials of arbitrary degree.  Fu-Li-Zhang-Zhang~\cite{fu2026} extended this work to compact
Kähler manifolds and proved a Nakai-Moishezon type criterion for
positivity with respect to a strictly right-Noetherian polynomial,
 including
the complex Hessian equation, the Hessian quotient equation, and the
critical LYZ equation, which is the starting point
for the present paper.}

In a semistable class, some of these intersection numbers may vanish.
We study the null locus of such classes and the existence and regularity
of weak solutions with prescribed density. 

{In the Monge-Ampère case, the null locus of a nef and big class is the
union of the positive-dimensional irreducible analytic subvarieties on
which the class has zero top self-\allowbreak intersection. Collins-Tosatti~\cite{CT15} proved that this locus coincides with the non-Kähler locus, extending Nakamaye's theorem~\cite{Nakamaye2000} to transcendental classes (see also~\cite{ELMNP2009,CT22,Dang2026Null,Vu2025Derivative} for related results).

For degenerate complex Hessian equations, Guo-Phong-Tong-Wang~\cite{guo2024}
obtained uniform estimates relative to the envelope in nef classes.
Guedj-Lu~\cite{GuedjLuDegenerateHessian} constructed bounded weak solutions
when the background form is semipositive and big. For the $k$-Hessian
equation with an exponential right-hand side, Pang-Sun-Wang-Zhou~\cite{PangSunWangZhou2025}
obtained bounded weak solutions when the background form is $k$-semipositive
and its class is $k$-big. For $2$-Hessian equation on compact Kähler threefolds,
Fu-Zhang-Zhang~\cite{FuZhangZhang2026} studied the null
divisors and showed that a  $2$-nef class with positive
total $2$-Hessian mass is $2$-big.

For the $J$-flow on compact K\"ahler surfaces, Fang-Lai-Song-Weinkove
proved a uniform $L^{\infty}$ estimate  and the smooth convergence outside the union of
finite  curves under a semipositive assumption~\cite{FLSW2014}.
Murakami constructed weak solutions
of generalized Monge-Ampère equation and supercritical LYZ equation under the semistable assumption~\cite{Murakami2026Boundary}. For the $J$-equation, Liu proved that under the smooth boundary cone condition
and the $J$-big assumption, the  $J$-null locus coincides with
the $J$-non-ample locus~\cite{Liu2026Boundary}.
Under the smooth boundary cone condition, Fu-Zhang constructed a weak solution
of the $J$-equation with bounded potential~\cite{FuZhang2026Boundary}.
In the $J$-semistable case, Liu proved that there are only
finitely many  irreducible $J$-null
subvarieties~\cite{Liu2026Finiteness} (see also~\cite{KhalidSjostromDyrefelt2026,SivaramSjostromDyrefelt2026}).
Fang-Ma proved that Murakami's weak solution is smooth outside
the numerical $J$-null locus and that no point of this locus has a neighborhood
on which the solution is $C^2$~\cite{FangMa2026JNull}. For the supercritical
LYZ equation, Fu-Yau-Zhang-Zhang~\cite{FuYauZhangZhang2026Boundary} proved
that under semistability the destabilizing prime divisors form
a finite exceptional family. Under the semisubsolution condition,
the weak solution has bounded potential and is smooth outside these divisors. For regularity results for degenerate complex Monge-Ampère equations, we refer to~\cite{BEGZ10,DDGPKZ2014,DiNezzaLu2017QuasiProjective,GuyettMcGowan2026}.}

{%
We show that for a semistable class of positive polynomial mass, the
numerical null locus coincides with the $f$-non-Kähler locus.
 We use a uniform estimate relative to the admissible envelope
to construct weak solutions with minimal singularities. For smooth
positive densities, these solutions are smooth outside the null locus.}

\subsection{Main results}

Let $(X,\chi)$ be a connected compact Kähler manifold of dimension $n$, and let $\alpha\in H^{1,1}(X,\mathbb R)$. Fix a smooth closed representative $\theta\in\alpha$ and {set} $\beta=[\chi]$. {%
Let $f$ be a monic real polynomial 
\[
 f(t)=t^d-\sum_{i=0}^{d-1}\binom di c_i t^i,
 \qquad 2\le d\le n.
\]
We assume that $f$ is right-Noetherian. For $0\le j<d$, let
$r(f^{(j)})$ denote the largest real root of $f^{(j)}$.
We write $F=\operatorname{Pol}_n(f)$ for its normalized polarization and
$C_f(\chi)$ for the associated positivity cone. See Section~\ref{sec:polynomial}
for the root condition on $f$, the formula for $F$, and the description
of $C_f(\chi)$ in terms of eigenvalues.}

{We call $\alpha$ $f$-semistable if, for every $\varepsilon>0$,
the class $\alpha+\varepsilon[\chi]$ contains a smooth closed real
$(1,1)$-form in $C_f(\chi)$. For the complex $k$-Hessian equation, $f(t)=t^k$ and
$C_f(\chi)=\Gamma_k(\chi)$, in this case we also call an $f$-semistable
class \emph{$k$-nef} with respect to $\chi$.}

{%
Set $m=n-d$. For real $(1,1)$-forms $\eta$ and $\xi$, define
\begin{equation}\label{eq:intersection-polynomials}
 \begin{aligned}
 \Phi_m^f(\eta,\xi)&=\xi^m,\\
 \Phi_{m+j}^f(\eta,\xi)
 &=\xi^m\left(\eta^j-\sum_{i=0}^{j-1}\binom ji
                 c_{d-j+i}\eta^i\xi^{j-i}\right),\qquad 1\le j\le d.
 \end{aligned}
\end{equation}
We use the same notation for forms and their cohomology classes.
Write $Q=\int_X\Phi_n^f(\alpha,\beta)$.}

{%
For an irreducible analytic subvariety $V$ of dimension $p>m$, the
intersection number $N_p(V)$ is the integral of $\Phi_p^f(\alpha,\beta)$
over $V$. The Hessian null locus $\operatorname{Null}_f(\alpha,\beta)$
is defined as the union of those $V\subsetneq X$ for which $N_p(V)=0$
(see Section~\ref{sec:null-locus}).}

Let $\mathcal Q_f(\alpha,\chi)$ denote the family of $f$-K\"ahler currents
in $\alpha$ with quasi-psh potentials, {as in
Definition~\ref{def:f-positive-current}.}
When this family is nonempty, we define the $f$-non-K\"ahler locus
of $\alpha$ with respect to $\chi$ by
\[
 E_{\mathrm{nK}}^{f}(\alpha,\chi)
 :=\bigcap_{T\in\mathcal Q_f(\alpha,\chi)}
       \{x\in X:\nu(T,x)>0\},
\]
where $\nu(T,x)$ denotes the Lelong number of $T$ at $x$. When $\alpha$ and $\chi$ are fixed, we simply write $E_{\mathrm{nK}}^{f}$.

{
\begin{theorem}\label{thm:null}
Suppose that $\alpha$ is $f$-semistable.
For any $1\le k\le d-1$ such that $r(f^{(k-1)})>r(f^{(k)})$,
the union
\[
 \bigcup_{\substack{V\subset X\ \text{irreducible analytic}\\
                    n-d<\dim V\le n-k,\ N_{\dim V}(V)=0}} V
\]
is an analytic subset of $X$ of dimension at most $n-k$.
If $Q>0$, then
\(\operatorname{Null}_f(\alpha,\beta)=E_{\mathrm{nK}}^f(\alpha,\chi).
\)
\end{theorem}
}
{%
Theorem~\ref{thm:null} generalizes the theorem of
Collins-Tosatti~\cite{CT15} for nef and big classes. For the $J$-equation,
Liu~\cite{Liu2026Boundary} and Fang-Ma~\cite{FangMa2026JNull}
obtained corresponding characterizations of the $J$-null locus.}

{
For the analytic results, assume that $\alpha$ is
$f$-semistable, and $Q>0$.
Let $\mathcal A_f(\theta)$ denote the admissible potentials of
Section~\ref{subsec:generalized-envelopes}. For $t\ge0$, the envelope $V_t$ of all
nonpositive potentials in $\mathcal A_f(\theta+t\chi)$ is defined
in~\eqref{eq:admissible-envelope}. We write $V_\theta$ for this envelope at $t=0$. For $u\in\mathcal A_f(\theta)$,
write $\mu_f(\theta,u)$ for the polynomial Hessian measure defined by the non-$d$-polar product (see Section~\ref{subsec:hessian-products}).}

{%
Write $W=\int_X\chi^n$ and $\rd\nu=\chi^n/W$. For $q>n$, we use the
Luxemburg norm~\eqref{eq:luxemburg-norm} associated with
$\Psi_q(s)=s^{n/d}\log^q(1+s)$, $s\ge0$.}

\begingroup
\begin{theorem}\label{thm:weak}
Fix $q>n$. Let $\rho\ge0$ satisfy $\int_X\rho\chi^n=1$ and $\|\rho\|_{\Psi_q,\nu}<\infty$. Then there exists a unique $u\in\mathcal A_f(\theta)$ with minimal singularities and $\sup_Xu=0$ satisfying
\[
 \mu_f(\theta,u)=Q\rho\chi^n,\qquad V_\theta-C\le u\le V_\theta,
\]
where the equation holds in the sense of non-$d$-polar products and the uniform constant $C>0$ depends only on $(X,\chi)$, $\theta$, $f$, the choice of the barrier $\psi$ from Proposition~\ref{prop:common-locus}, and an upper bound for $\|\rho\|_{\Psi_q,\nu}$. Moreover, the absolutely continuous part of $\theta+\ii\partial\bar\partial u$ satisfies
\[
 F\bigl((\theta+\ii\partial\bar\partial u)_{\rm ac}\bigr)=Q\rho
 \quad\text{a.e. on }X.
\]
\end{theorem}
\endgroup

Theorem~\ref{thm:weak} is new even for $J$-equation and $k$-Hessian equation. Fang-Ma~\cite{FangMa2026JNull} ask when the weak solution of the
semistable $J$-equation is globally bounded. The estimate in Theorem~\ref{thm:weak}
shows that the solution is indeed uniformly bounded relative to the admissible envelope. 
For the $k$-Hessian equation, as Lu-Nguyen~\cite{LN19Pure} point out, the lack of integrability
of $k$-subharmonic functions when $k<n$ (see \cite{Blocki2005}) makes relative $L^\infty$
estimates delicate. 
The main ingredient in overcoming this difficulty is the quasi-psh barrier provided by the Demailly-P{\u{a}}un type criterion established in \cite{fu2026} (see also Lemma~\ref{lem:tangential-current}). Our proof for the relative $L^\infty$ estimate is inspired by the auxiliary Monge-Ampère equation method of Guedj-Lu~\cite{guedj2025}. We compare the approximating solutions $u_t$ with combinations of the envelope $V_t$, the barrier $\psi$, and solutions of auxiliary Monge-Ampère equations in the fixed Kähler class $[\chi]$. The exponential integrability of $\psi$ then allows us to complete the sublevel set iteration.

{For $t>0$, write $\theta_t=\theta+t\chi$ and
$Q_t=\int_X\Phi_n^f(\theta_t,\chi)$.}

\begin{theorem}\label{thm:regularity}
Let $g$ be a smooth positive function on $X$ satisfying $\int_Xg\,\chi^n=1$, and let $u_g$ be the corresponding normalized weak solution given by Theorem~\ref{thm:weak}. For sufficiently small $t>0$, let $u_{t,g}$ be the smooth solution of
\[
 F(\theta_t+\ii\partial\bar\partial u_{t,g})=Q_tg,
 \qquad\theta_t+\ii\partial\bar\partial u_{t,g}\in C_f(\chi),
 \qquad\sup_Xu_{t,g}=0.
\]
Then $u_{t,g}\to u_g$ smoothly on compact subsets of $X\setminus\operatorname{Null}_f(\alpha,\beta)$ as $t\searrow0$.
\end{theorem}

We first prove a weighted gradient estimate for the exponential equation.
Using the stability estimate along a path of densities, we obtain a weighted
gradient estimate for $u_{t,g}$. 
{ Let $\psi$ denotes the function in
Proposition~\ref{prop:common-locus},
the gradient estimate and Proposition~\ref{prop:weighted-hessian} give}
\[
 \sup_{0<t\le t_0}\sup_{X\setminus E_{\mathrm{nK}}^{f}}
 \ee^{a\psi}
 \left(|\partial u_{t,g}|_\chi^2+
       |\ii\partial\bar\partial u_{t,g}|_\chi\right)\le C.
\]

In Section~\ref{sec:polynomial}, we recall the polynomial inequalities and
show that Demailly's regularization preserves the cone condition.
We use this regularization in Section~\ref{sec:null-locus} to identify
the numerical null locus with the $f$-non-K\"ahler locus.
In Sections~\ref{sec:intrinsic}--\ref{sec:limits}, we define the polynomial
Hessian measure, establish the relative $L^\infty$ estimate, and prove
existence and uniqueness of normalized weak solutions with minimal singularities.
In Section~\ref{sec:regularity}, we prove that the approximating solutions converge smoothly  outside the null locus.
\label{intro:end}

\medskip
\noindent\textbf{Declaration on the use of AI.}
ChatGPT 6 Astra was used to assist with the preparation of this manuscript, especially with the eigenvalue calculation of Proposition~\ref{prop:weighted-hessian}. The authors take full responsibility for the mathematical content and the final text.

\section{Preliminaries}
\label{sec:polynomial}

\subsection{Polynomial positivity}

We first recall the theory of right-Noetherian polynomials, introduced by Lin \cite{Lin2023Convexity} for $d=n$ and developed to $d\le n$ by Fang-Ma \cite{fangma2026}. For earlier work, see \cite{FangLaiMa2011,CollinsSzekelyhidi2017,Yuan2006,Chen2021J}.

Let $f(t)=t^d-\sum_{i=0}^{d-1}\binom di c_i t^i$ be a monic real polynomial of degree $d\le n$. For $0\le j\le d-1$, let $r(f^{(j)})$ denote the largest real root of $f^{(j)}$, when it exists. We say that $f$ is right-Noetherian if these roots exist and
\[
 r(f)\ge r(f')\ge\cdots\ge r(f^{(d-1)}).
\]

Let $\sigma_j$ denote the $j$-th elementary symmetric polynomial, with $\sigma_0=1$. The polarization of $f$ is defined as the unique symmetric multiaffine polynomial $F$ in $n$ variables such that $F(t,\ldots,t)=f(t)$, i.e.,
\[
 F(\lambda)=\operatorname{Pol}_n(f)(\lambda)
 =\frac{\sigma_d(\lambda)}{\binom nd}
 -\sum_{i=0}^{d-1}\binom di c_i\frac{\sigma_i(\lambda)}{\binom ni}.
\]

We denote by $C_f$ the connected component of $\{F>0\}$ containing $t\mathbf1$ for all sufficiently large $t$, and by $\overline C_f$ its closure. A real $(1,1)$-form belongs to $C_f(\chi)$ if, at every point, its eigenvalues with respect to $\chi$ lie in $C_f$. Correspondingly, we define $\overline C_f(\chi)$ by requiring the eigenvalues to lie in $\overline C_f$.

We also write \(G\coloneqq F^{1/d},\; E_j\coloneqq \sigma_j/\binom{n}{j}\)
and let \(\Gamma_d\) denote the Gårding cone of \(\sigma_d\), given by
\[
 \Gamma_d=\{\lambda\in\mathbb R^n:\sigma_j(\lambda)>0,\ 1\le j\le d\}.
\]

The set $C_f$ is convex \cite{Lin2023Convexity,fangma2026}, and $G$ is concave on $C_f$ \cite{fangma2026}. Moreover, the following inclusions hold:
\begin{equation}\label{eq:centered-branch}
 \overline C_f+\mathbb R_{>0}^n\subset C_f,
 \qquad C_f-c_{d-1}{\bf1}\subset\Gamma_d.
\end{equation}

For $I\subset\{1,\ldots,n\}$, write $\partial_I\coloneqq\prod_{i\in I}\frac{\partial}{\partial\lambda_i}$, with $\partial_\varnothing F\coloneqq F$. Any point $\lambda\in\mathbb R^n$ belongs to $C_f$ if and only if (see \cite[Lemma~2.4]{Lin2023Convexity} and \cite[Theorem~8.2]{fangma2026})
\begin{equation}\label{eq:partial-characterization}
 \partial_I F(\lambda)>0
 \quad\text{for all }I\subset\{1,\ldots,n\}
       \text{ such that }\partial_I F\not\equiv0.
\end{equation}

For a real $(1,1)$-form $\eta$, write $F(\eta)=F\bigl(\lambda_\chi(\eta)\bigr)$, where $\lambda_\chi(\eta)$ denotes the eigenvalues of $\eta$ with respect to $\chi$. It is easy to see that
\[
 F(\eta)\chi^n=\eta^d\wedge\chi^{n-d}-\sum_{i=0}^{d-1}\binom di c_i\,\eta^i\wedge\chi^{n-i}.
\]

We set $G(\eta)=F(\eta)^{1/d}$. Then, $G$ and $\log F$, viewed as functions of real $(1,1)$-forms, are concave on $C_f(\chi)$.
For a smooth function $H$ and a real $(1,1)$-form $\xi$, let $D_\xi^jH(\eta)\coloneqq\left.\frac{\mathrm d^j}{\mathrm dt^j}\right|_{t=0}H(\eta+t\xi)$ denote the $j$-th derivative in the direction $\xi$. When $H=F$, we expand $F(\eta+t\xi)\chi^n$ in powers of $t$ to obtain
\[
 \frac{(d-j)!}{d!}D_\xi^jF(\eta)\chi^n
 =\Phi_{n-j}^f(\eta,\chi)\wedge\xi^j,\qquad 0\le j\le d.
\]
In particular,
\[
 D_\chi^jF(\eta)=\operatorname{Pol}_n\bigl(f^{(j)}\bigr)\bigl(\lambda_\chi(\eta)\bigr).
\]

The following lemma will be used in the weighted gradient estimates.

\begin{lemma}\label{lem:root-trace}
Let $\eta\in C_f(\chi)$ and set $\xi=\eta-c_{d-1}\chi$. Then the following hold.
\begin{enumerate}
\item $C_f(\chi)+\Gamma_d(\chi)\subset C_f(\chi)$.
\item If $\gamma\in\Gamma_d(\chi)$, then $D_\gamma G(\eta)\ge E_d(\gamma)^{1/d}>0$.
\item If $d\ge2$, then
\begin{equation}\label{eq:root-trace}
 |\xi|_\chi\le\operatorname{tr}_\chi\xi
 \le nG(\eta)\bigl(D_\chi G(\eta)\bigr)^{d-1}.
\end{equation}
\end{enumerate}
\end{lemma}

\begin{proof}
We first prove that $C_f+\Gamma_d\subset C_f$. Let $\lambda\in C_f$ and $\mu\in\Gamma_d$. For every $I\subset\{1,\ldots,n\}$ with $|I|\le d$, the highest-degree part of $\partial_I F$ is $\partial_I E_d$. Hence
\[
 \lim_{t\to\infty}t^{-(d-|I|)}\partial_I F(\lambda+t\mu)
 =\partial_I E_d(\mu)>0.
\]
By \eqref{eq:partial-characterization}, it follows that $\lambda+t\mu\in C_f$ for all sufficiently large $t$. Since $C_f$ is convex and $\lambda\in C_f$, we obtain $\lambda+t\mu\in C_f$ for all $t\ge0$.

Now let $\eta\in C_f(\chi)$ and $\gamma\in\Gamma_d(\chi)$. Choose a $\chi$-unitary frame in which $\eta+\gamma$ is diagonal, and let $a$ and $b$ be the diagonal vectors of $\eta$ and $\gamma$, respectively. By the Schur-Horn theorem \cite{MarshallOlkinArnold2011} and the convexity of $C_f$ and $\Gamma_d$, we have $a\in C_f$ and $b\in\Gamma_d$. This implies $\lambda_\chi(\eta+\gamma)=a+b\in C_f$, and therefore $\eta+\gamma\in C_f(\chi)$.

By (1), we have $\eta+t\gamma\in C_f(\chi)$ for every $\gamma\in\Gamma_d(\chi)$ and $t\ge0$. By concavity of $G$, we have
\[
 G(\eta+t\gamma)\le G(\eta)+tD_\gamma G(\eta).
\]
Dividing by $t$ and letting $t\to\infty$, we obtain $D_\gamma G(\eta)\ge\lim_{t\to\infty}G(\eta+t\gamma)/t=E_d(\gamma)^{1/d}>0$. This proves (2).

For (3), take any semipositive real $(1,1)$-form $\kappa$. By a pointwise Khovanskii-Teissier type inequality (see \cite[Theorem~13.1]{fangma2026}), we have
\[
 \left(\frac{\Phi_n^f(\eta,\chi)}{\chi^n}\right)^{d-2}
 \frac{\Phi_{n-d+1}^f(\eta,\chi)\wedge\kappa^{d-1}}{\chi^n}
 \le
 \left(\frac{\Phi_{n-1}^f(\eta,\chi)\wedge\kappa}{\chi^n}\right)^{d-1}.
\]
By taking $\kappa=\chi$ and substituting
\[
 \frac{\Phi_n^f(\eta,\chi)}{\chi^n}=G(\eta)^d,
 \qquad
 \frac{\Phi_{n-1}^f(\eta,\chi)\wedge\chi}{\chi^n}
 =G(\eta)^{d-1}D_\chi G(\eta)
\]
into the preceding inequality and using $\Phi_{n-d+1}^f(\eta,\chi)\wedge\chi^{d-1}/\chi^n=\operatorname{tr}_\chi\xi/n$, we obtain $\operatorname{tr}_\chi\xi\le nG(\eta)\bigl(D_\chi G(\eta)\bigr)^{d-1}$. By \eqref{eq:centered-branch}, we also have $\xi\in\Gamma_d\subset\Gamma_2$, and hence
\[
 |\xi|_\chi^2=(\operatorname{tr}_\chi\xi)^2-2\sigma_2(\xi)<(\operatorname{tr}_\chi\xi)^2
\]
This completes the proof.

\end{proof}

We now turn to developing notions of positivity in the polynomial setting.

\begin{definition}\label{def:f-positive-nef}
We say that a class $\alpha\in H^{1,1}(X,\mathbb R)$ is $f$-positive (with respect to $\chi$) if it contains a smooth closed form in $C_f(\chi)$. It is called $f$-semistable (or $f$-nef) if $\alpha+\varepsilon\beta$ is $f$-positive for every $\varepsilon>0$.
\end{definition}

For a real $(1,1)$-current $T$ with measure coefficients, we write $T=T_{\rm ac}+T_{\rm s}$ for its Lebesgue decomposition with respect to $\chi^n$.

We recall the Demailly-P{\u{a}}un type criterion for Hessian type operators established in \cite{fu2026} (see also Lemma~\ref{lem:tangential-current}).

\begin{proposition}[{\cite[Proposition~4.3]{fu2026}}]\label{prop:semistable-current}
Suppose that $\alpha$ is $f$-semistable and $\int_X\Phi_n^f(\alpha,\beta)>0$. For a smooth representative $\theta\in\alpha$, there are a quasi-psh function $u\not\equiv-\infty$ and constants $b_T,\delta_T>0$ such that $T=\theta+\ii\partial\bar\partial u$ satisfies $T\ge-b_T\chi$ and $(T-\delta_T\chi)_{\rm ac}\in\overline C_f(\chi)$ a.e.
\end{proposition}


\begin{definition}\label{def:f-positive-current}
Let $T$ be a real $(1,1)$-current with measure coefficients. We call $T$ a \emph{$C_f$-positive current} (or simply \emph{$f$-positive}) with respect to $\chi$ if $T_{\rm ac}\in\overline C_f(\chi)$ a.e. and $T_{\rm s}\ge0$. We call $T$ an $f$-K\"ahler current if $T-\delta\chi$ is $f$-positive for some $\delta>0$. We denote by $\mathcal Q_f(\alpha,\chi)$ the set of $f$-K\"ahler currents in $\alpha$ with quasi-psh potentials.
\end{definition}

By Proposition~\ref{prop:semistable-current}, $\mathcal Q_f(\alpha,\chi)$ is nonempty whenever $\alpha$ is $f$-semistable and $\int_X\Phi_n^f(\alpha,\beta)>0$.

\subsection{Regularization of quasi-psh functions under the cone condition}
\label{sec:regularization}

In this section, we show that Demailly's Bergman regularization preserves the cone condition. The forms considered here need not be closed. First, we need the following spectral truncation lemma.

\begin{lemma}\label{lem:truncation}
Let $(X,\omega)$ be a Hermitian manifold, and let $\mathcal C(\omega)$ be a set of real $(1,1)$-forms determined by their $\omega$-eigenvalues and stable under addition of semipositive forms. For every $\varepsilon>0$, there exists $b>0$, depending only on $\mathcal C$ and $\varepsilon$, such that every positive $(1,1)$-current $T$ with $T_{\rm ac}\in\mathcal C(\omega)$ a.e. satisfies
\[
 (T_{\rm ac})^{(b)}+\varepsilon\omega\in\mathcal C(\omega)
 \quad\textnormal{a.e.}
\]
Here $(T_{\rm ac})^{(b)}$ is obtained by replacing each $\omega$-eigenvalue $\lambda_i$ of $T_{\rm ac}$ with $\min\{\lambda_i,b\}$.
\end{lemma}

\begin{proof}
Suppose the conclusion fails for some $\varepsilon>0$. There exist $b_j\to\infty$, positive currents $T_j$, and points $x_j\in X$ such that $T_{j,\rm ac}(x_j)\ge0$ and
\[
 T_{j,\rm ac}(x_j)\in\mathcal C(\omega_{x_j}),
 \qquad
 \bigl(T_{j,\rm ac}(x_j)\bigr)^{(b_j)}
 +\varepsilon\omega_{x_j}\notin\mathcal C(\omega_{x_j}).
\]

Let $\lambda_{j,1}(x)\le\cdots\le\lambda_{j,n}(x)$ be the eigenvalues of $T_{j,\rm ac}(x)$ with respect to $\omega_x$. Up to extracting a subsequence, we may assume that $\lambda_{j,i}(x_j)\to l_i\in[0,\infty]$ for every $i$. Fix $k$ sufficiently large. Then for all sufficiently large $j$, we have
\[
 \min\{\lambda_{j,i}(x_j),b_j\}+\varepsilon
 \ge\lambda_{k,i}(x_k),
 \qquad 1\le i\le n.
\]
The stability under addition of semipositive forms then leads to
\[
 \bigl(T_{j,\rm ac}(x_j)\bigr)^{(b_j)}
 +\varepsilon\omega_{x_j}\in\mathcal C(\omega_{x_j}),
\]
a contradiction.
\end{proof}

\begin{remark}
For the $J$-equation, the error term $\varepsilon$ can be estimated explicitly under spectral truncation, which Chen uses to control the fiber integral near the diagonal (see \cite[pp.~559--560]{Chen2021J}). For $\mathcal C(\omega)=\Gamma_n(\omega)$, any $b\ge0$ works.
\end{remark}

Let $U\subset\mathbb C^n$ be a bounded pseudoconvex domain, and let $\varphi\in\mathrm{PSH}(U)$. For an integer $k\ge1$, consider the weighted Bergman space
\[
 \mathcal H_{k\varphi}(U)\coloneqq
 \left\{h\in\mathcal O(U):\|h\|_{k\varphi}^2\coloneqq
 \int_U|h|^2\ee^{-k\varphi}\rd V<\infty\right\}.
\]
Let $(f_j)$ be an orthonormal basis of $\mathcal H_{k\varphi}(U)$. The Bergman kernel is defined by
\[
 B_{k\varphi}(z,\zeta)\coloneqq\sum_jf_j(z)\overline{f_j(\zeta)},\qquad
 B_{k\varphi}(z)\coloneqq B_{k\varphi}(z,z).
\]

We prove the local Bergman kernel curvature estimate for a given cone condition, motivated by \cite[Lemma~3.3]{Liu2026MinimalSlope}.

\begin{lemma}\label{lem:bergman-form-bound}
Let $U\subset\mathbb C^n$ be a bounded pseudoconvex domain, and let $\omega_0$ be the standard Kähler form on $\mathbb C^n$. Let $w$ be psh on a neighborhood $V$ of $\overline U$. Suppose that
\[
 \ii\partial\bar\partial w\ge\gamma\quad\text{on }V,\qquad \gamma(x)\in\mathcal S\subseteq\{0\le\eta\le b\omega_0\}\quad\text{for a.e. }x\in V,
\]
where $\gamma$ is a measurable $(1,1)$-form on $V$, $\mathcal S$ is a compact convex set, and $b>0$ is a constant. For any $a>0$ and integer $k\ge1$, let $\varphi=w+a|z|^2$ and $B_k=B_{k\varphi}$. Then, at every $x\in U$ with $B_k(x)>0$, there exists $\bar\gamma_x\in\mathcal S$ such that
\begin{equation}\label{eq:bergman-form-bound}
 \ii\partial\bar\partial\bigl(k^{-1}\log B_k-a|z|^2\bigr)
 \ge\bar\gamma_x-\frac{b^2}{a}\omega_0.
\end{equation}
\end{lemma}

\begin{proof}
Define $\varphi_k\coloneqq k^{-1}\log B_k$ and $H_{\varphi_k}\coloneqq((\varphi_k)_{i\bar j})$. Fix $x\in U$ with $B_k(x)>0$, and set
\[
 e_x(z)\coloneqq\frac{1}{B_k(x)^{1/2}}\sum_jf_j(z)\overline{f_j(x)},
\]
which is the unique element of $\mathcal H_{k\varphi}(U)$ satisfying
\[
 \|e_x\|_{k\varphi}=1,\qquad
 \langle h,e_x\rangle_{k\varphi}=\frac{h(x)}{B_k(x)^{1/2}}
 \quad\text{for every }h\in\mathcal H_{k\varphi}(U).
\]
In particular, $\rd\mu_x\coloneqq|e_x|^2\ee^{-k\varphi}\rd V$ is a probability measure on $U$. We claim that (this also follows from Berndtsson's convexity for psh variations \cite[Theorem~1.1]{Berndtsson2006Subharmonicity})
\begin{equation}\label{eq:bergman-harmonic}
 H_{\varphi_k}(x)\ge\left(\int_U(aI+\gamma)^{-1}\rd\mu_x\right)^{-1}.
\end{equation}

Let $\ell$ be a constant complex linear form, and set
\[
 F_\ell(z)=\ell(z-x)e_x(z),\qquad
 s_\ell(z)=\overline{\ell(z-x)}e_x(z).
\]
By the Bergman-Fuks formula, for $\xi\ne0$,
\[
 \min_{\substack{h\in\mathcal H_{k\varphi}(U)\\h(x)=0,\ \partial_\xi h(x)=1}}
 \|h\|_{k\varphi}^{2}
 =\frac1{kB_k(x)H_{\varphi_k}(x)(\xi,\bar\xi)}.
\]
Since $F_\ell(x)=0$ and $\partial_\xi F_\ell(x)=\ell(\xi)B_k(x)^{1/2}$, we obtain
\[
 |\ell(\xi)|^2\le kH_{\varphi_k}(x)(\xi,\bar\xi)\,\|F_\ell\|_{k\varphi}^2.
\]
Choosing $\ell$ with $\ell(\xi)\ne0$ shows that $H_{\varphi_k}(x)$ is positive definite. Taking the supremum over $\xi$ yields
\[
 \frac1k\langle H_{\varphi_k}(x)^{-1}\ell,\ell\rangle
 \le\|F_\ell\|_{k\varphi}^2.
\]
By the reproducing property, $\langle s_\ell,h\rangle_{k\varphi}=\langle e_x,\ell(z-x)h\rangle_{k\varphi}=0$ for every $h\in\mathcal H_{k\varphi}(U)$. Thus $s_\ell$ is the minimal solution of $\bar\partial s=e_x\overline{\ell(\rd z)}$. Since $\ii\partial\bar\partial\varphi\ge a\omega_0+\gamma$, H\"ormander's $L^2$-estimates for the $\bar\partial$-equation \cite[Theorem~5.1]{demailly2009analytic} give
\[
 \|F_\ell\|_{k\varphi}^2=\|s_\ell\|_{k\varphi}^2
 \le\frac1k\int_U\langle(aI+\gamma)^{-1}\ell,\ell\rangle\rd\mu_x.
\]

We have
\[
 H_{\varphi_k}(x)\ge\rho_x^{-1},\qquad
 \rho_x\coloneqq\int_U(aI+\gamma)^{-1}\rd\mu_x.
\]

Set $\bar\gamma_x\coloneqq\int_U\gamma\,\rd\mu_x$. Since $\mathcal S$ is compact and convex, $\bar\gamma_x\in\mathcal S$. The identity $\rho_x^{-1}-2aI+a^2\rho_x=\rho_x^{-1}(I-a\rho_x)^2\ge0$ yields
\[
 \begin{aligned}
 H_{\varphi_k}(x)-aI
 &\ge\rho_x^{-1}-aI\ge aI-a^2\rho_x\\
 &=\bar\gamma_x-\int_U\gamma^2(aI+\gamma)^{-1}\rd\mu_x
 \ge\bar\gamma_x-\frac{b^2}{a}I.
 \end{aligned}
\]
This finishes the proof.
\end{proof}

We recall that a quasi-psh function $u$ has analytic singularities if there exist a modification $\pi:\widetilde X\to X$, an effective divisor $D$ with simple normal crossings, and a constant $c>0$ such that $u\circ\pi=c\log|s|^2+h$ locally on $\widetilde X$, where $s$ is a local equation for $D$ and $h$ is smooth.

For a quasi-psh function $v$, the multiplier ideal sheaf $\mathcal I(v)$ consists of the holomorphic germs $g$ for which $|g|^2\ee^{-v}$ is locally integrable. Any orthonormal basis of $\mathcal H_{k\varphi}(U)$ generates the multiplier ideal sheaf $\mathcal I(k\varphi)$ on $U$ (see \cite[Proposition~5.7]{demailly2009analytic}).

Although it remains a problem to obtain the global regularization of $f$-positive currents (see Definition~\ref{def:f-positive-current}), Proposition~\ref{prop:semistable-current} allows us to work with quasi-psh potentials and adapt Demailly's Bergman kernel approximation \cite{Demailly1992Reg,Dem15}.

\begin{proposition}\label{prop:strict-regularization}
Assume that $\mathcal C(\omega)$ is a closed convex set as above and is stable under addition of semipositive forms. Let $u$ be a quasi-psh function and set $T=\theta+\ii\partial\bar\partial u$. Assume that $T_{\rm ac}\in\mathcal C(\omega)$ a.e. For every sequence $\varepsilon_j\searrow0$, there exist currents $T_j\coloneqq\theta+\ii\partial\bar\partial u_j$ with analytic singularities such that $u_j\searrow u$ and
\[
 (T_j+\varepsilon_j\omega)_{\rm ac}\in\mathcal C(\omega)
 \quad\textnormal{a.e.}
\]
\end{proposition}

\begin{proof}
Since $u$ is quasi-psh, we may replace $\theta$ by $\theta+b\omega$ and $\mathcal C(\omega)$ by $(\mathcal C(\omega)+b\omega)\cap\{\eta\ge0\}$, with $b$ large enough that $T\ge0$. Fix $\varepsilon>0$. By Lemma~\ref{lem:truncation}, we can choose $M=M(\varepsilon)$ such that $\gamma_\varepsilon=(T_{\rm ac})^{(M)}$ satisfies
\begin{equation}\label{eq:regularization-lower-form}
 0\le\gamma_\varepsilon\le M\omega,\qquad
 T\ge\gamma_\varepsilon,\qquad
 \gamma_\varepsilon+\varepsilon\omega\in\mathcal C(\omega).
\end{equation}

Choose coordinate balls $U'_{\varepsilon,i}\subset\!\subset U_{\varepsilon,i}\subset\!\subset V_{\varepsilon,i}$, with the $U'_{\varepsilon,i}$ covering $X$ and $\omega=\omega_0$ at the center of each $V_{\varepsilon,i}$. Taking the balls sufficiently small, we may assume that $\gamma_\varepsilon$ lies a.e. in the compact convex set
\[
 \mathcal S_\varepsilon=\left\{\eta:0\le\eta\le2M\omega_0,\ 
 \eta+2\varepsilon\omega_0\in\mathcal C(\omega_0)\right\}.
\]
We can also assume that there is a real quadratic function $h_{\varepsilon,i}$ on a neighborhood of $\overline V_{\varepsilon,i}$ such that
\[
 0\le\ii\partial\bar\partial h_{\varepsilon,i}-\theta\le\varepsilon\omega,\qquad
 \ii\partial\bar\partial(u+h_{\varepsilon,i})\ge T\ge\gamma_\varepsilon.
\]

Choose $a_\varepsilon>0$ with $4M^2/a_\varepsilon\le\varepsilon/2$ and set $s_{\varepsilon,i}=h_{\varepsilon,i}+a_\varepsilon|z|^2$. Let $B_{\varepsilon,i,m}$ be the Bergman kernel on $V_{\varepsilon,i}$ for the weight $\ee^{-m(u+s_{\varepsilon,i})}$, and set $v_{\varepsilon,i,m}=m^{-1}\log B_{\varepsilon,i,m}-s_{\varepsilon,i}$. By Lemma~\ref{lem:bergman-form-bound}, for each $x$ with $B_{\varepsilon,i,m}(x)>0$, there exists $\bar\gamma_{\varepsilon,x}\in\mathcal S_\varepsilon$ such that
\[
\begin{aligned}
 \theta+\ii\partial\bar\partial v_{\varepsilon,i,m}
 &=\ii\partial\bar\partial\bigl(m^{-1}\log B_{\varepsilon,i,m}-a_\varepsilon|z|^2\bigr)
   +\theta-\ii\partial\bar\partial h_{\varepsilon,i}\\
 &\ge\bar\gamma_{\varepsilon,x}-\tfrac{\varepsilon}{2}\omega_0-\varepsilon\omega
 \ge\bar\gamma_{\varepsilon,x}-2\varepsilon\omega,
\end{aligned}
\]
so $\theta+\ii\partial\bar\partial v_{\varepsilon,i,m}+5\varepsilon\omega\in\mathcal C(\omega)$. For each fixed $\varepsilon$, \cite[Theorem~2.2.1, Step~3]{DPS01} provides a constant $A_\varepsilon>0$ such that, for every $i$ and all positive integers $m_1,m_2$,
\[
 v_{\varepsilon,i,m_1+m_2}\le\frac{A_\varepsilon}{m_1+m_2}
 +\frac{m_1}{m_1+m_2}v_{\varepsilon,i,m_1}+\frac{m_2}{m_1+m_2}v_{\varepsilon,i,m_2}.
\]
By taking $m_1=m_2=2^k$, we see that $v_{\varepsilon,i,2^k}+A_\varepsilon/2^k$ decreases with $k$. Choose smooth functions $-1\le\rho_{\varepsilon,i}\le1$ equal to $1$ on $U'_{\varepsilon,i}$ and to $-1$ near $\partial U_{\varepsilon,i}$. Set $\tau_m=m^{-1/2}$. For large $m$, we have $\ii\partial\bar\partial(\tau_m\rho_{\varepsilon,i})\ge-\varepsilon\omega$. Then choose $K>0$ to be fixed below, and define
\[
 v_{\varepsilon,m}(x)\coloneqq\frac{A_\varepsilon}{m}+K\tau_m+
 \max_{\tau_m/2}\{v_{\varepsilon,i,m}(x)+\tau_m\rho_{\varepsilon,i}(x)\mid x\in U_{\varepsilon,i}\}.
\]
where $\max_{\tau_m/2}$ denotes the regularized maximum. The overlap estimate in \cite[Lemma~3.6]{Demailly1992Reg} yields
\begin{equation}\label{eq:regularization-overlap}
 |v_{\varepsilon,i,m}-v_{\varepsilon',l,m}|\le C_{\varepsilon,\varepsilon'}\frac{\log m}{m}
 \quad\text{on }U_{\varepsilon,i}\cap U_{\varepsilon',l}.
\end{equation}
Take $m$ large enough so that $C_{\varepsilon,\varepsilon}(\log m)/m<\tau_m$. Then $v_{\varepsilon,i,m}+\tau_m\rho_{\varepsilon,i}$ does not contribute to $v_{\varepsilon,m}$ near $\partial U_{\varepsilon,i}$. A simple calculation shows that $\ii\partial\bar\partial v_{\varepsilon,m}\ge\sum_i\lambda_i\,\ii\partial\bar\partial(v_{\varepsilon,i,m}+\tau_m\rho_{\varepsilon,i})$, where $\lambda_i\ge0$ and $\sum_i\lambda_i=1$. Hence $\theta+\ii\partial\bar\partial v_{\varepsilon,m}+6\varepsilon\omega\in\mathcal C(\omega)$ outside $Z_m=V(\mathcal I(mu))$.

By \cite[proof of Theorem~2.2.1, Step~2]{DPS01}, the local Bergman kernel for $\varphi=u+s_{\varepsilon,i}$ can be written locally as
\begin{equation}\label{eq:bergman-kernel-expansion}
 B_{\varepsilon,i,m}(z,w)=\sum_{j=1}^{N}a_j(z,w)f_j(z)\overline{f_j(w)},
\end{equation}
where the $a_j$ are holomorphic in $z$ and antiholomorphic in $w$, and $f_1,\ldots,f_N$ generate $\mathcal I(m\varphi)$ locally. Hence $v_{\varepsilon,i,m}$ has analytic singularities, and so does $v_{\varepsilon,m}$. Also, by \cite[Proposition~3.1]{Demailly1992Reg}, $v_{\varepsilon,m}\to u$ in $L^1(X)$ as $m\to\infty$.

Applying \eqref{eq:regularization-overlap} with $\varepsilon'=\varepsilon/2$, we obtain, after enlarging $C_{\varepsilon,\varepsilon/2}$,
\[
\begin{aligned}
 v_{\varepsilon/2,l,2m}+\frac{A_{\varepsilon/2}}{2m}
 &\le v_{\varepsilon,i,2m}+\frac{A_\varepsilon}{2m}
   +C_{\varepsilon,\varepsilon/2}\frac{\log(2m)}{2m}\\
 &\le v_{\varepsilon,i,m}+\frac{A_\varepsilon}{m}
   +C_{\varepsilon,\varepsilon/2}\frac{\log(2m)}{2m}.
\end{aligned}
\]
Since $\log(2m)/(2m)=o(\tau_m)$, we fix $K$ sufficiently large so that, for each $\varepsilon>0$, there exists $m(\varepsilon)$ such that for $m\ge m(\varepsilon)$,
\[
 v_{\varepsilon/2,2m}\le v_{\varepsilon,m},\qquad
 v_{\varepsilon,2m}\le v_{\varepsilon,m},\qquad
 \|v_{\varepsilon,m}-u\|_{L^1(X)}<\varepsilon.
\]

Let $m_k=2^{k+1}$. Choose $\varepsilon_k\searrow0$ (with $\varepsilon_{k+1}\in\{\varepsilon_k,\varepsilon_k/2\}$) so that $m_k\ge m(\varepsilon_k)$. Set $u_k\coloneqq v_{\varepsilon_k,m_k}$. Then $u_k\searrow u$.
\end{proof}

\begin{remark}
For every $0<\lambda<\lambda'$, the functions $u_j$ constructed above satisfy $\mathcal I(\lambda' u_j)\subset\mathcal I(\lambda u)$ when $j$ is sufficiently large (see \cite[Corollary~1.12]{Dem15}).
\end{remark}

\section{The Hessian null locus}\label{sec:null-locus}

\subsection{The Hessian non-K\"ahler locus and existence of an analytic barrier}

Assume that $\mathcal Q_f(\alpha,\chi)$ is nonempty. For $T\in\mathcal Q_f(\alpha,\chi)$, let $E_+(T)\coloneqq\{x\in X:\nu(T,x)>0\}$ denote the locus where the Lelong number of $T$ is positive. The $f$-non-K\"ahler locus of $\alpha$ (with respect to the form $\chi$) is then defined as
\[
 E_{\mathrm{nK}}^{f}(\alpha,\chi)\coloneqq\bigcap_{T\in\mathcal Q_f(\alpha,\chi)}E_+(T).
\]

\begin{proposition}\label{prop:common-locus}
There exist constants $C,\kappa>0$ and a quasi-psh function $\psi$ with analytic singularities such that
\[
 \theta+\ii\partial\bar\partial\psi\in\mathcal Q_f(\alpha,\chi),
 \qquad \operatorname{Sing}(\psi)=E_{\mathrm{nK}}^{f}.
\]
Moreover,
\refstepcounter{equation}\label{eq:barrier-gradient}%
$|\partial\psi|_\chi^2\le C\ee^{-\kappa\psi}$ on $X\setminus E_{\mathrm{nK}}^{f}$.\hfill\textup{(\theequation)}
\end{proposition}

\begin{proof}
We write $\mathcal Q_f^{\mathrm{an}}(\alpha,\chi)$ for the currents in $\mathcal Q_f(\alpha,\chi)$ with analytic singularities. By Proposition~\ref{prop:strict-regularization}, each $T\in\mathcal Q_f(\alpha,\chi)$ admits a regularization $T'\in\mathcal Q_f^{\mathrm{an}}(\alpha,\chi)$ with $E_+(T')\subset E_+(T)$. Hence
\[
 E_{\mathrm{nK}}^{f}=\bigcap_{T\in\mathcal Q_f^{\mathrm{an}}(\alpha,\chi)}E_+(T).
\]
By the Noetherian property of analytic sets, the intersection can be taken over finitely many currents $T_i=\theta+\ii\partial\bar\partial u_i$. For $w=\max_\varepsilon(u_1,\ldots,u_N)$, we have $\theta+\ii\partial\bar\partial w\in\mathcal Q_f(\alpha,\chi)$ and $\{w=-\infty\}=\bigcap_i E_+(T_i)=E_{\mathrm{nK}}^{f}$.

Let $\psi$ be a regularization of $w$ from Proposition~\ref{prop:strict-regularization}. We still have $\{\psi=-\infty\}=E_{\mathrm{nK}}^{f}$. Let $v_i$ denote the local potential constructed from the $k$-th Bergman kernel so that $\psi=\max_\varepsilon(v_1,\ldots,v_r)$. By \eqref{eq:bergman-kernel-expansion}, on each smaller chart we have
\[
 |\partial v_i|_\chi^2=k^{-2}\ee^{-2kv_i}|\partial\ee^{kv_i}|_\chi^2\le C\ee^{-2kv_i}.
\]
Writing $\lambda_i=(\partial_{t_i}\max_\varepsilon)(v_1,\ldots,v_r)$, we have $\partial\psi=\sum_i\lambda_i\partial v_i$, with $\lambda_i\ge0$ and $\sum_i\lambda_i=1$. Since $v_i\ge\psi-C_0$ whenever $\lambda_i>0$, the preceding estimate implies \eqref{eq:barrier-gradient} with $\kappa=2k$.
Moreover, $\ee^\psi$ is H\"older continuous with exponent $1/k$.
\end{proof}

Let $V\subset X$ be an irreducible analytic subvariety of dimension $p$, and let $\kappa$ be a Kähler form on $X$. We consider the intersection numbers
\[
 I_{pq}(V;\kappa)\coloneqq\int_V\kappa^{p-q}\Phi_q^f(\alpha,\beta),
 \qquad \max\{1,n-d\}\le q\le\min\{p,n-1\},
\]
where $d=\deg f$. For $n-d<p\le n$, we write $N_p(V)\coloneqq\int_V\Phi_p^f(\alpha,\beta)$, so that $N_p(V)=I_{pp}(V;\kappa)$.

In \cite{fu2026}, a Nakai-Moishezon type criterion was established for $f$-positivity in terms of the mixed intersection numbers $I_{pq}(V;\kappa)$. In particular, if $\alpha$ is numerically semistable in the sense that $\int_X\Phi_n^f(\alpha,\beta)$ and all $I_{pq}(V;\kappa)$ are nonnegative, then $\alpha+s[\kappa]$ is numerically strictly positive, and hence $f$-positive for every $s>0$. Therefore $\alpha$ is $f$-semistable.

Now, given an $f$-semistable class $\alpha$, we show that it suffices to test the numbers $N_p(V)$ to conclude $f$-positivity. In other words, the convex set of $f$-positive classes is a connected component of $\{\alpha : N_p(V)>0 \text{ for all } V\}$.

\begin{proposition}\label{prop:intersection-positivity}
Suppose that $\alpha$ is $f$-semistable. Then $I_{pq}(V;\kappa)\ge0$, and its vanishing is independent of the choice of $\kappa$. If $N_p(V)>0$, then $I_{pq}(V;\kappa)>0$ for every $q$.
\end{proposition}

\begin{proof}
By semistability, we can choose $s_j\searrow0$ and $\omega_j\in C_f(\chi)$ with $[\omega_j]=\alpha+s_j\beta$. By \eqref{eq:partial-characterization}, the forms $\Phi_q^f(\omega_j,\chi)$ are strongly positive. Therefore, by integrating against $\kappa^{p-q}$ on $V$ and letting $j\to\infty$, we obtain $I_{pq}(V;\kappa)\ge0$. If $a\kappa_2\le\kappa_1\le b\kappa_2$, the same argument shows that
\[
 a^{p-q}I_{pq}(V;\kappa_2)\le I_{pq}(V;\kappa_1)
 \le b^{p-q}I_{pq}(V;\kappa_2).
\]
Thus it suffices to check the positivity for $\kappa=\chi$.

Let $\pi:W\to V$ be a resolution and $\rho$ a Kähler form on $W$. Set $\gamma=\pi^*\chi$ and $g=\frac{(p-n+d)!}{d!}\,f^{(n-p)}$. To apply \cite[(4.10)]{fu2026} on $W$, we choose $\varepsilon_j\searrow0$, let $M_j=\max\{1,\sup_X\lambda_{\max,\chi}(\omega_j)\}$, and define
\[
 \eta_j=\pi^*\omega_j+\varepsilon_j\rho,\qquad
 h_j=\gamma+\zeta_j\rho,\qquad
 0<\zeta_j<\min\{j^{-1},\varepsilon_j/M_j\}.
\]

At $x\in W$, let $\lambda_1\le\cdots\le\lambda_n$ be the eigenvalues of $\omega_j$ with respect to $\chi$ at $\pi(x)$, and let $\mu_1\le\cdots\le\mu_p$ be those of $\eta_j$ with respect to $h_j$ at $x$. Since $\omega_j\in C_f(\chi)$, it follows from \eqref{eq:partial-characterization} that $(\lambda_1,\ldots,\lambda_p)\in C_g$. Moreover, for $t=\lambda_l$, the form $\eta_j-th_j=\pi^*(\omega_j-t\chi)+(\varepsilon_j-t\zeta_j)\rho$ has at most $l-1$ negative eigenvalues. We then infer that $\mu_l\ge\lambda_l$ for $1\le l\le p$, which implies $\eta_j\in C_g(h_j)$.

Now set $U_q=\int_W\gamma^{p-q}\Phi_q^g(\pi^*\alpha,[\gamma])$ for $n-d\le q\le p$. By hypothesis and the projection formula, $U_p=N_p(V)>0$, and we shall show inductively that $U_{p-1},\ldots,U_{n-d+1}>0$. For $0\le r\le\deg g-2$, let $\tilde g$ be the monic normalization of $g^{(r)}$ and let $\tilde d$ denote its degree. As the coefficient of $t^{\tilde d-1}$ in $\tilde g(t)$ is $-\tilde d\,c_{d-1}$, we obtain from \cite[(4.10)]{fu2026} that
\[
 \left(\int_W(\eta_j-c_{d-1}h_j)h_j^{p-1}\right)
 \left(\tilde d\int_Wh_j\Phi_{p-1}^{\tilde g}(\eta_j,h_j)\right)
 \ge \left(\int_W\Phi_p^{\tilde g}(\eta_j,h_j)\right)\int_Wh_j^p.
\]
Letting $j\to\infty$, we see that the first factor on the left has a nonnegative limit since $\eta_j-c_{d-1}h_j\in\Gamma_{\tilde d}(h_j)$. Also, the second factor tends to $\tilde d U_{p-r-1}$ and the right-hand side tends to $U_{p-r}\int_W\gamma^p$. The induction hypothesis $U_{p-r}>0$ and the identity $U_{n-d}=\int_V\chi^p>0$ then imply $U_{p-r-1}>0$, which completes the proof.
\end{proof}

Suppose that $\alpha$ is $f$-semistable and $\int_X\Phi_n^f(\alpha,\beta)>0$. We define the Hessian numerical null locus $\operatorname{Null}_f(\alpha,\beta)$ as the union of all irreducible analytic subvarieties $V\subseteq X$ with $\dim V=p>n-d$ and $N_p(V)=0$.

\begin{lemma}\label{lem:null-forcing}
If $I_{pq}(V;\kappa)=0$, then $V\subset E_{\mathrm{nK}}^{f}$. In particular, $\operatorname{Null}_f(\alpha,\beta)\subseteq E_{\mathrm{nK}}^{f}(\alpha,\chi)$.
\end{lemma}

\begin{proof}
Suppose that $V$ intersects $X\setminus E_{\mathrm{nK}}^{f}$. We can then choose a coordinate ball $B\Subset X\setminus E_{\mathrm{nK}}^{f}$ intersecting $V_{\rm reg}$, and a neighborhood $U$ of $E_{\mathrm{nK}}^{f}$ disjoint from $B$. Write $\omega_j=\theta+s_j\chi+\ii\partial\bar\partial v_j$.

The potential $\psi$ given by Proposition~\ref{prop:common-locus} satisfies $\theta+\ii\partial\bar\partial\psi-\eta\chi\in\overline C_f$ outside $E_{\mathrm{nK}}^{f}$ for some $\eta>0$. Choose $c_j$ so that $\psi+c_j>v_j+2\varepsilon$ on $X\setminus U$, and set $w_j=\max_\varepsilon\{\psi+c_j,v_j\}$. Then $w_j$ is smooth and the form $\widetilde\omega_j\coloneqq\theta+s_j\chi+\ii\partial\bar\partial w_j$ belongs to $\overline C_f$ on $X$ and satisfies $\widetilde\omega_j-\eta\chi\in\overline C_f$ on $B$. By the binomial expansion, $\Phi_q^f(\widetilde\omega_j,\chi)\ge\eta^{q-m}\chi^q$ there. Hence
\[
 \int_V\kappa^{p-q}\Phi_q^f(\widetilde\omega_j,\chi)
 \ge\eta^{q-m}\int_{V_{\rm reg}\cap B}\kappa^{p-q}\chi^q>0.
\]
Letting $j\to\infty$, we obtain a contradiction since $I_{pq}(V;\kappa)=0$.
\end{proof}

\subsection{Extension of potentials with analytic singularities under cone condition}

Before proceeding with the proof of Theorem~\ref{thm:null}, we need the following singular Demailly-P{\u{a}}un type criterion.

\begin{lemma}\label{lem:tangential-current}
Let $V\subseteq X$ be an irreducible analytic subvariety of dimension $p>n-d$ such that $N_p(V)>0$. Let $\pi:W\to V$ be a resolution and $\rho$ a Kähler form on $W$ with $\rho\ge\gamma\coloneqq\pi^*\chi$. Set $g=\frac{(p-n+d)!}{d!}f^{(n-p)}$. Then there exist $c\in\mathbb R$, $t,\delta>0$, and a quasi-psh function $v$ on $W$ such that, for $\widetilde\gamma=\gamma+t\rho$, the current $T\coloneqq\pi^*\theta-c\gamma+\ii\partial\bar\partial v$ satisfies
\[
 T\ge\delta\widetilde\gamma,\qquad (T-\delta\widetilde\gamma)_{\rm ac}\in\overline{C}_{g(\,\cdot+c)}(\widetilde\gamma).
\]
\end{lemma}

\begin{proof}
Let $\widetilde d=\deg g=p-n+d$ and $\sigma=\pi^*\theta$. By the proof of Proposition~\ref{prop:intersection-positivity} and the projection formula,
\[
 J_j\coloneqq\frac{\widetilde d!}{(\widetilde d-j)!}\int_W\rho^j\Phi_{p-j}^g([\sigma],[\gamma])
 \ge\frac{\widetilde d!}{(\widetilde d-j)!}\int_V\chi^j\Phi_{p-j}^f(\alpha,\beta)>0,
 \qquad 0\le j\le\widetilde d.
\]

Also, as in the proof of Proposition~\ref{prop:intersection-positivity}, choose $s_j\searrow0$ and $\omega_j\in C_f(\chi)$ with $[\omega_j]=\alpha+s_j\beta$, and set $\eta_j=\pi^*\omega_j+\varepsilon_j\rho$ and $\gamma_j=\gamma+\zeta_j\rho$. For each $j$, the tuples $(\lambda_1,\ldots,\lambda_p)$ there range over a compact subset of $C_g$. By \cite[Lemma~2.2]{fu2026}, we may choose $L_j\to\infty$ so that these tuples lie in $C_{P_j}$, where $P_j(t)=(t+L_j)^{n-d}g(t)$. Since $\mu_l\ge\lambda_l$, we have $\eta_j\in C_{P_j}(\gamma_j)$ for all sufficiently small $\zeta_j>0$.

Write $-pc_j$ for the coefficient of $t^{p-1}$ in $P_j$, and choose $\zeta_j$ so that $|c_j|\zeta_j\to0$. Since $L_j^{-(n-d)}P_j\to g$ coefficientwise and $([\eta_j],[\gamma_j])\to([\sigma],[\gamma])$, a simple calculation shows (cf. \cite[(6.4)]{fu2026})
\[
 \lambda_j\coloneqq
 \frac{\int_W\Phi_p^{P_j}(\eta_j,\gamma_j)}
      {p\int_W\rho\,\Phi_{p-1}^{P_j}(\eta_j,\gamma_j)}
 \longrightarrow \frac{J_0}{J_1}>0.
\]

By \cite[Proposition~4.2]{fu2026}, there is a current $S_j$ satisfying
\[
 [S_j]=[\eta_j],\qquad S_j-\lambda_j\rho\ge c_j\gamma_j,\qquad
 (S_j-\lambda_j\rho)_{\rm ac}\in\overline{C}_{P_j}(\gamma_j).
\]

Since $\lambda_j-\varepsilon_j-s_j+c_j\zeta_j\to J_0/J_1>0$, we may now fix $j$ sufficiently large and choose a sufficiently small $\delta>0$ such that $\lambda_j-\varepsilon_j-s_j+c_j\zeta_j\ge(1+\zeta_j)\delta$.

Now we have $[S_j]=[\sigma]+s_j[\gamma]+\varepsilon_j[\rho]$ and then define $T=S_j-(s_j\gamma+\varepsilon_j\rho)-c_j\gamma$. Since $\gamma_j-\gamma=\zeta_j\rho$ and $\gamma\le\rho$, our choice of $\delta$ yields
\[
\begin{aligned}
 0\le S_j-\lambda_j\rho-c_j\gamma_j
 &=T-\bigl((\lambda_j-\varepsilon_j)\rho-s_j\gamma\bigr)-c_j(\gamma_j-\gamma)\\
 &\le T-(1+\zeta_j)\delta\rho\le T-\delta\gamma_j.
\end{aligned}
\]
By the cone inclusion $C_{P_j}\subset C_g$ \cite[Lemma~3.1]{CNX26}, we have $(S_j-\lambda_j\rho-c_j\gamma_j)_{\rm ac}\in\overline{C}_{g(\,\cdot+c_j)}(\gamma_j)$. Therefore $(T-\delta\gamma_j)_{\rm ac}\in\overline{C}_{g(\,\cdot+c_j)}(\gamma_j)$ and $T\ge\delta\gamma_j$, with $[T]=[\sigma]-c_j[\gamma]$, as desired.
\end{proof}

By regularizing the potential obtained in Lemma~\ref{lem:tangential-current} and taking a further resolution, we have the following.

\begin{proposition}\label{prop:divisorial-preparation}
Let $V\subseteq X$ be an irreducible analytic subvariety of dimension $n>p>n-d$ such that $N_p(V)>0$. There is a modification $\pi:\widetilde X\to X$, biholomorphic over a general point of $V$, such that the strict transform $W$ of $V$ is smooth. Denote the pullback $\pi^*\chi$ by $\gamma$. There are a Kähler form $\widetilde\gamma\ge\gamma$ on $\widetilde X$, constants $c\in\mathbb R$, $\delta>0$, and a quasi-psh function $v$ on $W$ such that the current $T\coloneqq(\pi^*\theta-c\gamma)|_W+\ii\partial\bar\partial v$ satisfies
\[
 T\ge\delta\widetilde\gamma|_W,\qquad
 (T-\delta\widetilde\gamma|_W)_{\rm ac}\in
 \overline{C}_{g(\,\cdot+c)}(\widetilde\gamma|_W)\quad\text{a.e.}
\]
Moreover, $v=\sum_l a_l\log|s_l|^2|_W+r$ with $a_l>0$ and $r\in C^\infty(W)$, where the $s_l$ are the defining sections of divisors on $\widetilde X$ with simple normal crossings transverse to $W$.
\end{proposition}

\begin{proof}
Take an embedded desingularization $\pi_1$ of $V$ with smooth strict transform $W_1$ (see \cite{Hir77,Wlodarczyk2009}) and choose a Kähler form $\rho_1\ge\pi_1^*\chi$ on the ambient manifold containing $W_1$. Apply Lemma~\ref{lem:tangential-current} with $\rho_1|_{W_1}$ and set $\omega=\pi_1^*\chi+t\rho_1$. By Proposition~\ref{prop:strict-regularization}, we choose a regularization $v_k$ such that, for some $\eta>0$,
\[
 T_k\coloneqq\bigl(\pi_1^*\theta-c\pi_1^*\chi\bigr)|_{W_1}
       +\ii\partial\bar\partial v_k\ge\eta \omega|_{W_1},\qquad
 (T_k-\eta \omega|_{W_1})_{\rm ac}\in\overline{C}_{g(\,\cdot+c)}(\omega|_{W_1})\quad\text{a.e.}
\]

Again, we take an embedded log resolution $\pi_2$ of the singularity of $v_k$, and set $\pi=\pi_1\circ\pi_2$ and $\nu=\pi_2|_W$, so that $\nu^*v_k$ can be expressed as desired. Then there exists an effective exceptional divisor $E$ such that $\rho\coloneqq C\pi_2^*\omega-\Theta_E$ is a Kähler form for a suitable smooth metric on $\mathcal O(E)$ and sufficiently large $C$.

Applying Lemma~\ref{lem:truncation} to $T_k-\eta\omega|_{W_1}$, we obtain $b>0$ and a form $\xi\coloneqq((T_k-\eta\omega|_{W_1})_{\rm ac})^{(b)}\ge0$ such that
\[
 \xi+\frac{\eta}{4}\omega|_{W_1}\in\overline{C}_{g(\,\cdot+c)}(\omega|_{W_1}).
\]

Fix $a,\varepsilon>0$ to be determined and denote
\[
 v=\nu^*v_k+a\log|s_E|^2|_W,\qquad
 \widetilde\gamma=\pi_2^*\omega+\varepsilon\rho .
\]
Now we consider the current $T=(\pi^*\theta-c\gamma)|_W+\ii\partial\bar\partial v$. Outside $E$, we have $T_{\rm ac}=\nu^*((T_k)_{\rm ac})-a\Theta_E|_W$. Let $\lambda$ be an eigenvalue of $\xi+\frac{\eta}{4}\omega|_{W_1}$ with respect to $\omega|_{W_1}$. Since $(T_k)_{\rm ac}\ge\xi+\eta\omega|_{W_1}$, we have
\[
\begin{aligned}
 T_{\rm ac}-(\lambda+\eta/4)\widetilde\gamma|_W
 &\ge \nu^*(\xi+\eta\omega|_{W_1})-a\Theta_E|_W
       -(\lambda+\eta/4)\widetilde\gamma|_W\\
 &\ge \nu^*\!\left(\xi+\frac{\eta}{4}\omega|_{W_1}-\lambda\omega|_{W_1}\right)\\
 &\quad +(\eta/2-aC)\nu^*(\omega|_{W_1})
       +\bigl(a-\varepsilon(\lambda+\eta/4)\bigr)\rho|_W .
\end{aligned}
\]
Choose $a$ so that $aC\le\eta/2$ and then $\varepsilon$ so that $\varepsilon(b+\eta/2)<a$. Since $\lambda\le b+\eta/4$, the last two terms are semipositive. By an argument similar to that in the proof of Proposition~\ref{prop:intersection-positivity}, we have $T_{\rm ac}\ge\frac{\eta}{4}\widetilde\gamma|_W$ and $T_{\rm ac}-\frac{\eta}{4}\widetilde\gamma|_W\in\overline{C}_{g(\,\cdot+c)}(\widetilde\gamma|_W)$, which completes the proof.
\end{proof}

Let us revisit a fundamental result on the local extension under cone condition (see \cite[Lemma~9.1]{CNX26} and \cite[Lemma~6.2]{fu2026}).

\begin{lemma}\label{lem:smooth-local-extension}
Let $W\subset X$ be a smooth compact submanifold of dimension $p$, and let $u\in C^\infty(W,\mathbb R)$. If $p\le n-d$, then $u$ admits a smooth extension $\widetilde u$ near $W$ such that $\theta+\ii\partial\bar\partial\widetilde u\in C_f(\chi)$. If $p>n-d$, the same holds provided that
\[
 \theta|_W+\ii\partial\bar\partial u\in C_g(\chi|_W),
 \qquad g=\frac{(d-n+p)!}{d!}f^{(n-p)}.
\]
\end{lemma}

\begin{proof}
Choose a smooth extension $u_0$ of $u$ near $W$ and set $\widetilde u=u_0+L\,\operatorname{dist}_\chi(\cdot,W)^2$. Write $h=\chi|_W$ and $\eta=\theta|_W+\ii\partial\bar\partial u$. Fix $\varepsilon>0$ small enough that $\eta-\varepsilon h\in C_g(h)$ when $p>n-d$. Along $W$, for $\tau\in TW$ and $\nu\perp TW$, we have
\[
 (\theta+\ii\partial\bar\partial\widetilde u)(\tau+\nu,\overline{\tau+\nu})
 \ge(\eta-\varepsilon h)(\tau,\bar\tau)+(cL-C_\varepsilon)|\nu|_\chi^2,
\]
where $c>0$ and $C_\varepsilon>0$ are independent of $L$. Let $\lambda$ be the $h$-eigenvalues of $\eta-\varepsilon h$, and write $F=\operatorname{Pol}_n(f)$. For $I\subseteq\{1,\ldots,p\}$ and $J\subseteq\{p+1,\ldots,n\}$, write $a=|I|$ and $b=|J|$. When $a+b\le d$, the leading term of $\partial_I\partial_J F(\lambda,t\mathbf 1_{n-p})$ is
\[
 \begin{cases}
 \displaystyle\frac{p!\,d!}{n!\,(d-n+p)!}\,
 \partial_I\operatorname{Pol}_p(g)(\lambda)\,t^{n-p-b},
 & p>n-d,\quad a\le d-n+p,\\[7pt]
 \displaystyle\frac{\binom{n-p-b}{d-a-b}}{\binom nd}\,t^{d-a-b},
 & \text{otherwise}.
 \end{cases}
\]
Therefore, $\partial_I\partial_J F(\lambda,t\mathbf 1_{n-p})>0$ for sufficiently large $t$, and by \eqref{eq:partial-characterization} we have $(\lambda,t\mathbf 1_{n-p})\in C_f$. This implies $\theta+\ii\partial\bar\partial\widetilde u\in C_f(\chi)$.
\end{proof}

\begin{theorem}\label{thm:null-locus-equality}
Suppose that $\alpha$ is $f$-semistable and that $\int_X\Phi_n^f(\alpha,\beta)>0$. Then $\operatorname{Null}_f(\alpha,\beta)=E_{\mathrm{nK}}^{f}(\alpha,\chi)$.
\end{theorem}

\begin{proof}
We use an extension argument as in~\cite{CT15}.
From Lemma~\ref{lem:null-forcing}, we already have $\operatorname{Null}_f(\alpha,\beta)\subset E_{\mathrm{nK}}^f$. Now let $Y$ be an irreducible component of $E_{\mathrm{nK}}^f$ of dimension $p$. To prove the reverse inclusion, it suffices to show that $p>n-d$ and $N_p(Y)=0$.

Suppose that $p>n-d$ and $N_p(Y)>0$. By Proposition~\ref{prop:divisorial-preparation}, there is a modification $\pi:\widetilde X\to X$ with smooth strict transform $W$ of $Y$, biholomorphic over a general point of $Y$. Write $\gamma=\pi^*\chi$ and $g=\frac{(d-n+p)!}{d!}f^{(n-p)}$. For some K\"ahler form $\widetilde\gamma\ge\gamma$, $c\in\mathbb R$ and $\eta>0$, the current $T\coloneqq(\pi^*\theta-c\gamma)|_W+\ii\partial\bar\partial v$ satisfies
\[
 T\ge\eta\widetilde\gamma|_W,\qquad
 (T-\eta\widetilde\gamma|_W)_{\rm ac}\in\overline{C}_{g(\,\cdot+c)}(\widetilde\gamma|_W)
 \quad\text{a.e. on }W.
\]

If $p\le n-d$, also take an embedded desingularization $\pi:\widetilde X\to X$ of $Y$ with smooth strict transform $W$ and let $\gamma=\pi^*\chi$. By \cite[Lemma~3.5]{DP04}, there exist an effective exceptional divisor $E$ and a K\"ahler form $\widetilde\gamma\ge\gamma$ such that, for $v=(\log|s_E|^2)|_W$ and $c$ sufficiently negative,
\[
 \bigl((\pi^*\theta-c\gamma)|_W+\ii\partial\bar\partial v\bigr)_{\rm ac}\ge\widetilde\gamma|_W.
\]

In either case, write $v=\sum_l a_l(\log|s_l|^2)|_W+r$. Apply the construction of Lemma~\ref{lem:smooth-local-extension} to the smooth tangential form
\[
 \widetilde T\coloneqq\Bigl(\pi^*\theta-c\gamma-\sum_l a_l\Theta_l\Bigr)|_W+\ii\partial\bar\partial r,
\]
where $\Theta_l$ is the curvature of the metric defining $|s_l|^2$, using the polynomial $f(\,\cdot+c)$. Since $\widetilde T$ is positive, we obtain an extension $\varphi=\sum_l a_l\log|s_l|^2+\widetilde r$ of $v$ on a neighborhood $U$ of $W$ such that, for $S=\pi^*\theta-c\gamma+\ii\partial\bar\partial\varphi$ and some $\delta>0$,
\[
 S\ge\delta\widetilde\gamma,\qquad
 (S-\delta\widetilde\gamma)_{\rm ac}\in\overline{C}_{f(\,\cdot+c)}(\widetilde\gamma).
\]

Choose $W\subset U_0\Subset U$. Let $D$ be the reduced exceptional divisor of $\pi$ and $A'$ the union of the components of $E_{\mathrm{nK}}^f$ other than $Y$. It follows from \cite[Lemma~2.1]{DP04} that there exists a quasi-psh function $\zeta$ with logarithmic poles precisely on $A'$ and smooth elsewhere. By adding a small positive multiple of $\log|s_D|^2+\zeta\circ\pi$ to $\varphi$ near $\overline U_0$, we may assume that $\xi=S_{\rm ac}-(\delta/2)\widetilde\gamma$ is a positive current and belongs to $\overline{C}_{f(\,\cdot+c)}(\widetilde\gamma)$. Now we have $\lambda(\xi)\ge 0$ a.e. and $\widetilde\gamma\ge\gamma$. It follows that
\[
 \pi^*\theta+\ii\partial\bar\partial\varphi\ge c\gamma,\qquad
 (\pi^*\theta+\ii\partial\bar\partial\varphi)_{\rm ac}-(\delta/2)\gamma\in\overline{C}_f(\gamma)
 \quad\text{a.e.}
\]

Let $\psi$ be the potential given in Proposition~\ref{prop:common-locus}. Near $\partial U_0$, one sees immediately that the poles of $\psi\circ\pi$ are contained in those of $\varphi$, since $\pi^{-1}(Y)\subset W\cup D$. Then, by taking a common log resolution of the singularities of $\varphi$ and $\psi\circ\pi$, we may choose $0<\lambda<1$ sufficiently small that $\varphi\le\lambda\psi\circ\pi+O(1)$ near $\partial U_0$.

Choose $\eta>0$ such that $(\theta+\ii\partial\bar\partial\psi)_{\rm ac}-\eta\chi\in\overline C_f(\chi)$. Since $\alpha$ is $f$-semistable, for every $s>0$ there is a smooth $v_s$ with $\theta+s\chi+\ii\partial\bar\partial v_s\in C_f(\chi)$. Take $s$ sufficiently small that $t=\lambda\eta-(1-\lambda)s>0$. For $u=\lambda\psi+(1-\lambda)v_s$, we have
\[
 \theta+\ii\partial\bar\partial u-t\chi
 =\lambda(\theta+\ii\partial\bar\partial\psi-\eta\chi)
 +(1-\lambda)(\theta+s\chi+\ii\partial\bar\partial v_s),
\]
with its absolutely continuous part belonging to $\overline C_f(\chi)$. Therefore, after subtracting a constant from $\varphi$, we may assume that $\varphi\le u\circ\pi-1$ near $\partial U_0$. For $0<\varepsilon<1/2$, the function $\max_\varepsilon\{\varphi,u\circ\pi\}$ on $U_0$ therefore extends by $u\circ\pi$ outside. One sees that the current $\pi^*\theta+\ii\partial\bar\partial\max_\varepsilon\{\varphi,u\circ\pi\}$ descends to a current in $\mathcal Q_f(\alpha,\chi)$ on $X$ which is smooth near a general point of $Y$. This leads to a contradiction to $Y\subset E_{\mathrm{nK}}^f$. Therefore $E_{\mathrm{nK}}^f\subset\operatorname{Null}_f(\alpha,\beta)$, completing the proof.
\end{proof}

\begin{theorem}\label{thm:finite-null-subvarieties}
Suppose that $\alpha$ is $f$-semistable. Let $1\le k\le d-1$ and set $p=n-k$. If $r(f^{(k-1)})>r(f^{(k)})$, then the union
\[
 \bigcup_{\substack{V\subset X\ \text{irreducible analytic}\\
                     n-d<\dim V\le p,\ N_{\dim V}(V)=0}} V
\]
is an analytic subset of $X$ of dimension at most $p$.
\end{theorem}

\begin{proof}
Denote $r_i=r(f^{(i)})$. We let $a=r_k$, and consider
\[
 h(t)=f(t)-\sum_{i=0}^{k-1}\frac{f^{(i)}(a)}{i!}(t-a)^i.
\]

Noting that we have $r(h^{(i)})=a$ for $i<k$, while $r(h^{(i)})=r_i$ for $k\le i<d$. In particular $h$ is right-Noetherian. Since $r(h^{(i)})\le r_i$ for every $0\le i<d$, we obtain that $C_f\subset C_h$ by \cite[Theorem~10.10]{FangMaGarding2026} (see also \cite[Theorem~1.3]{Lin2023Convexity}). Hence $\alpha$ is also $h$-semistable. Now we claim that $\int_X\Phi_n^h(\alpha,\beta)>0$ and that the union in the statement is precisely $\operatorname{Null}_h(\alpha,\beta)$, which suffices for our purposes by Theorem~\ref{thm:null-locus-equality}.

Indeed, since $h^{(j)}=f^{(j)}$ for every $j\ge k$, we have $\Phi_q^h=\Phi_q^f$ whenever $n-d<q\le p$ by \eqref{eq:intersection-polynomials}. This implies that the union of all $f$-null subvarieties of dimension at most $p$ is contained in $\operatorname{Null}_h(\alpha,\beta)$. On the other hand, fix $p<q\le n$, set $j=n-q$, and write $f_q=(d-j)!f^{(j)}/d!$ and $h_q=(d-j)!h^{(j)}/d!$. Since $j<k$, we have $r_j\ge r_{k-1}>a$, so $\delta_q=h_q(r_j)>0$. Also, as $h_q$ is strictly increasing on $(a,\infty)$, we infer that the largest root of the polynomial $h_q-\delta_q$ is $r_j$. Denote by $C_{f_q}$ and $C_{h_q-\delta_q}$ the cones defined by the polarizations of $f_q$ and $h_q-\delta_q$ in $q$ variables, respectively. Comparing the roots of $h_q-\delta_q$ and $f_q$, and using \cite[Theorem~10.10]{FangMaGarding2026} again, we have the inclusion $C_{f_q}\subset C_{h_q-\delta_q}$, and hence
\begin{equation}\label{eq:auxiliary-polynomial-lower}
 \operatorname{Pol}_q(h_q)(\lambda)>\delta_q
 \qquad\text{for }\lambda\in C_{f_q}.
\end{equation}

Let $\eta\in C_f(\chi)$. By \eqref{eq:partial-characterization}, we see that each $q$-tuple of its eigenvalues with respect to $\chi$ belongs to $C_{f_q}$. Applying \eqref{eq:auxiliary-polynomial-lower} to each such tuple, we obtain
\[
 \Phi_q^h(\eta,\chi)\ge\delta_q\chi^q.
\]

Now, for any $f$-semistable class $\alpha$ and every irreducible $q$-dimensional subvariety $Y\subseteq X$ with $q>p$, we have
\[
 \int_Y\Phi_q^h(\alpha,\beta)\ge\delta_q\int_Y\beta^q>0.
\]

This implies that every $h$-null subvariety of $\alpha$ has dimension at most $p$ and that $\int_X\Phi_n^h(\alpha,\beta)>0$, finishing the proof.
\end{proof}

\section{Admissible potentials}\label{sec:intrinsic}

Building on the work of Bedford-Taylor on the Monge-Ampère operator \cite{BT76,BT82}, Błocki \cite{Blocki2005} and Dinew-Kołodziej \cite{DK2014} developed pluripotential methods to study bounded weak solutions of Hessian equations. The relative pluripotential theory \cite{BEGZ10,DDNL18} was also extended to Hessian equations on compact Kähler manifolds by Lu-Nguyen \cite{LN14Pure,LN19Pure}.

In this section, we involve some basic notions needed to study weak solutions of generalized Hessian type equations.

\subsection{Admissible potentials and envelopes}\label{subsec:generalized-envelopes}

We say that a locally integrable function $u$ is ($f$-)admissible, and write $u\in\mathcal A_f(\theta)$, if $\theta+\ii\partial\bar\partial u$ is an $f$-positive current as in Definition~\ref{def:f-positive-current}.

Let $\mathrm{SH}_d(X,\chi)$ be the class of $(\chi,d)$-subharmonic functions, defined distributionally by
\[
 (\chi+\ii\partial\bar\partial v)\wedge\eta_1\wedge\cdots
       \wedge\eta_{d-1}\wedge\chi^{n-d}\ge0
\]
for all smooth $\Gamma_d$-positive $(1,1)$-forms $\eta_1,\ldots,\eta_{d-1}$.

\begin{lemma}\label{lem:measure-hessian}
Assume $d\ge2$. If $u/A\in\mathrm{SH}_d(X,\chi)$ for some $A>0$, then $T:=A\chi+\ii\partial\bar\partial u$ has finite Radon measure coefficients and $T_{\rm s}\ge0$.
\end{lemma}

\begin{proof}
Take smooth $u_j\searrow u$ with $u_j/A\in\mathrm{SH}_d(X,\chi)$ as in \cite{LN14Pure}, and let $\omega_j=A\chi+\ii\partial\bar\partial u_j$. Since $\Gamma_d\subset\Gamma_2$,
\[
 |\omega_j|_\chi\le\operatorname{tr}_\chi\omega_j,\qquad
 \int_X\operatorname{tr}_\chi\omega_j\,\chi^n
       =nA\int_X\chi^n.
\]
By pairing with smooth test functions, it follows that $\ii\partial\bar\partial u_j$ has bounded total variation, and so does $\ii\partial\bar\partial u$.

It remains to prove that $T_{\rm s}\ge0$. Set $S=\ii\partial\bar\partial u$. Locally, we write $S_{\rm s}=\eta\,|S_{\rm s}|$, where $|S_{\rm s}|$ is the total variation measure and $\eta$ is the density. Since $\rd S=0$, \cite[Theorem~1.1]{dephilippisrindler2016} implies that, for $|S_{\rm s}|$-almost every $x$, there exists a nonzero real covector $\xi$ such that
\[
 \xi\wedge\eta(x)=0.
\]
As $\eta(x)$ is a nonzero real $(1,1)$-form, one has $\eta(x)=a\,\ii\zeta\wedge\bar\zeta$, where $a\in\mathbb R$ and $\zeta=\xi^{1,0}\ne0$. Moreover, since $S_{\rm s}\wedge\chi^{n-1}$ is positive, it follows that $a|\zeta|_\chi^2=\operatorname{tr}_\chi\eta(x)\ge0$. Thus $a\ge0$, and therefore $T_{\rm s}=S_{\rm s}\ge0$.
\end{proof}

\begin{remark}
When $d=1$, $\ii\partial\bar\partial u$ may not have measure coefficients. In $\mathbb C^n$ with $n\ge2$, the function $u(z)=-|z|^{2-2n}$ is subharmonic, but $|\ii\partial\bar\partial u|\asymp|z|^{-2n}$ on the punctured ball.
\end{remark}

Recall that an upper semicontinuous function $u$ satisfies $\theta+\ii\partial\bar\partial u\in\overline C_f(\chi)$ in the viscosity sense if, for every $x\in X$ and every $C^2$ function $q$ with $q\ge u$ near $x$ and $q(x)=u(x)$, one has $(\theta+\ii\partial\bar\partial q)(x)\in\overline C_f(\chi_x)$. Combining the above lemma with \cite{HL13Linear}, we have the following.

\begin{proposition}\label{prop:admissibility-equivalence}
Assume $d\ge2$, and let $u\in L^1(X)$ satisfy $u=u^*$, where $u^*(x):=\lim_{r\searrow0}\operatorname*{ess\,sup}_{B(x,r)}u$. Set $T=\theta+\ii\partial\bar\partial u$. The following are equivalent:
\begin{enumerate}
\item $u\in\mathcal A_f(\theta)$, i.e., $T$ has measure coefficients, $T_{\rm ac}\in\overline C_f(\chi)$ a.e., and $T_{\rm s}\ge0$.
\item In any local $\chi$-unitary coframe, $T\wedge\Omega\ge s(\Omega)\chi^n$ in the distributional sense for every semipositive real $(n-1,n-1)$-form $\Omega$ with constant coefficients in that coframe for which $s(\Omega):=\inf_{\eta\in\overline C_f(\chi)}(\eta\wedge\Omega)/\chi^n$ is finite.
\item $\theta+\ii\partial\bar\partial u\in\overline C_f(\chi)$ in the viscosity sense.
\end{enumerate}
Moreover, if (i) holds, then for every $A>0$ with $(A+c_{d-1})\chi-\theta\ge0$, we have $u/A\in\mathrm{SH}_d(X,\chi)$.
\end{proposition}

\begin{proof}
Condition (i) implies (ii). By definition, $T_{\rm ac}\in\overline C_f(\chi)$ a.e. and $T_{\rm s}\ge0$. Thus $T_{\rm ac}\wedge\Omega\ge s(\Omega)\chi^n$ a.e. and $T_{\rm s}\wedge\Omega\ge0$, as required in (ii).

Conversely, assume (ii) and fix $A>0$ with $(A+c_{d-1})\chi-\theta\ge0$. Convolving the coefficients of $T$ in a local $\chi$-unitary coframe preserves the inequalities in (ii), so the resulting forms $T_\varepsilon$ belong to $\overline C_f(\chi)$. Since $\overline C_f(\chi)-c_{d-1}\chi\subset\overline\Gamma_d(\chi)$, we have
\[
 T_\varepsilon+A\chi-\theta
 =(T_\varepsilon-c_{d-1}\chi)+\bigl((A+c_{d-1})\chi-\theta\bigr)
 \in\overline\Gamma_d(\chi).
\]
Testing against $\eta_1\wedge\cdots\wedge\eta_{d-1}\wedge\chi^{n-d}$ and letting $\varepsilon\searrow0$, we obtain $u/A\in\mathrm{SH}_d(X,\chi)$. Lemma~\ref{lem:measure-hessian} now shows that $T$ has measure coefficients and $T_{\rm s}\ge0$. Moreover, $T_\varepsilon\to T_{\rm ac}$ a.e. by the Lebesgue differentiation theorem, and hence $T_{\rm ac}\in\overline C_f(\chi)$ a.e. Hence $u\in\mathcal A_f(\theta)$.

Finally, we prove the equivalence of (ii) and (iii). By convexity, $\overline C_f(\chi)$ is the intersection of the half-spaces defined by the inequalities in (ii). Fix a strictly positive test form $\Omega$ with $s(\Omega)>-\infty$. Since $q\mapsto(\ii\partial\bar\partial q\wedge\Omega)/\chi^n$ is locally uniformly elliptic and $u=u^*$, it follows immediately from \cite[Corollary~9.4]{HL13Linear} that $T\wedge\Omega\ge s(\Omega)\chi^n$ holds in the distributional sense if and only if it holds in the viscosity sense. For semipositive $\Omega$, apply this to $\Omega+\varepsilon\chi^{n-1}$ and let $\varepsilon\searrow0$, using $s(\Omega+\varepsilon\chi^{n-1})\to s(\Omega)$.
\end{proof}

\begin{corollary}\label{cor:admissible-closure}
Assume $d\ge2$. The following properties hold.
\begin{enumerate}
\item The class $\mathcal A_f(\theta)$ is convex. If $\theta_j\to\theta$ smoothly and $u_j\in\mathcal A_f(\theta_j)$ converge to $u$ in $L^1$, then $u^*\in\mathcal A_f(\theta)$. In particular, a decreasing limit belongs to $\mathcal A_f(\theta)$ whenever it is not identically $-\infty$.
\item If $u,v\in\mathcal A_f(\theta)$, then $\max\{u,v\}\in\mathcal A_f(\theta)$.
\item If $\{u_\lambda\}_{\lambda\in\Lambda}$ is a nonempty family in $\mathcal A_f(\theta)$ that is locally uniformly bounded above, then $(\sup_{\lambda\in\Lambda}u_\lambda)^*\in\mathcal A_f(\theta)$.
\end{enumerate}
\end{corollary}

\begin{proof}
For (i), convexity follows directly from the definition of $\mathcal A_f(\theta)$. Also, the $L^1$ closedness follows from the distributional characterization in Proposition~\ref{prop:admissibility-equivalence}.

For (ii), let $w=\max\{u,v\}$. Suppose that a $C^2$ function $q$ satisfies $q\ge w$ near $x$ and $q(x)=w(x)$. If $w(x)=u(x)$, then $q\ge u$ near $x$ and $q(x)=u(x)$. Proposition~\ref{prop:admissibility-equivalence} therefore shows that $(\theta+\ii\partial\bar\partial q)(x)\in\overline C_f(\chi_x)$. The same argument applies when $w(x)=v(x)$. Applying Proposition~\ref{prop:admissibility-equivalence} again, we obtain $w\in\mathcal A_f(\theta)$.

For (iii), let $w=(\sup_{\lambda\in\Lambda}u_\lambda)^*$. By Choquet's lemma, there is a sequence $\{u_{\lambda_j}\}$ such that $w=(\sup_j u_{\lambda_j})^*$. Set $v_j=\max\{u_{\lambda_1},\ldots,u_{\lambda_j}\}$. By (ii), $v_j\in\mathcal A_f(\theta)$. There exists $A>0$ such that $v_j/A\in\mathrm{SH}_d(X,\chi)$ for all $j$, which implies that $v_j\to w$ in $L^1$. Hence $w\in\mathcal A_f(\theta)$ by (i).
\end{proof}

For $t\ge0$, we set $\theta_t=\theta+t\chi$ and define the admissible envelope
\begin{equation}\label{eq:admissible-envelope}
 V_t\coloneqq\left(\sup\{v\in\mathcal A_f(\theta_t):v\le0\}\right)^*,
\end{equation}
where the star denotes upper semicontinuous regularization. We say that $u\in\mathcal A_f(\theta_t)$ has minimal singularities if $u=V_t+O(1)$ on $X$.

It then follows from Corollary~\ref{cor:admissible-closure}(iii) that $V_t\in\mathcal A_f(\theta_t)$. Moreover, $V_t$ is bounded for $t>0$, and by definition we have $V_t\searrow V_\theta$ as $t\searrow0$.

\subsection{Hessian products and capacity}\label{subsec:hessian-products}

We recall the Hessian products that will be used to define the polynomial measure (see \cite{Lu2013Degenerate,LN14Pure,DinewLu2015}). For bounded $u_1,\ldots,u_l\in\mathrm{SH}_d(X,\chi)$, the products
\[
 (\chi+\ii\partial\bar\partial u_1)\wedge\cdots
 \wedge(\chi+\ii\partial\bar\partial u_l)\wedge\chi^{n-d},
 \qquad 0\le l\le d,
\]
are defined recursively as closed positive currents and are weakly continuous along bounded decreasing sequences. For a Borel set $E\subset X$, the Hessian capacity of $E$ is defined as
\[
 \operatorname{Cap}_d(E)\coloneqq
 \sup_{\substack{v\in\mathrm{SH}_d(X,\chi)\\-1\le v\le0}}
 \int_E(\chi+\ii\partial\bar\partial v)^d\wedge\chi^{n-d}.
\]

Recall that a set $E\subset X$ is $d$-polar if $E\subset\{v=-\infty\}$ for some $v\in\mathrm{SH}_d(X,\chi)$.
A function $u\in\mathrm{SH}_d(X,\chi)$ is quasi-continuous with respect to this capacity if for every $\varepsilon>0$, there is an open set $G$ with $\operatorname{Cap}_d(G)<\varepsilon$ such that $u|_{X\setminus G}$ is continuous. We say that $u_j$ converges to $u$ in capacity if
\[
 \operatorname{Cap}_d\bigl(\{|u_j-u|>\delta\}\bigr)\longrightarrow0
 \quad\text{for every }\delta>0.
\]

For an unbounded $\varphi$ with $\varphi/A\in\mathrm{SH}_d(X,\chi)$, set $\varphi_k=\max(\varphi,-k)$. For $0\le p\le d$, the sequence
\[
 \mathbf1_{\{\varphi>-k\}}(A\chi+\ii\partial\bar\partial\varphi_k)^p\wedge\chi^{n-d}
\]
is non-decreasing in $k$ and converges weakly to a closed positive current which puts no mass on $d$-polar sets.
For $0\le j\le d$, we define $\Phi_{n-j}^f(\theta_\varphi,\chi)$ by expanding the right-hand side of \eqref{eq:intersection-polynomials} using $\theta_\varphi=(A\chi+\ii\partial\bar\partial\varphi)+(\theta-A\chi)$ and taking the non-$d$-polar products defined above. The resulting currents are independent of $A$. For $j=0$, we write $\mu_f(\theta,\varphi)\coloneqq\Phi_n^f(\theta_\varphi,\chi)$.

We shall show that they are positive when $\varphi$ is admissible (see also \cite{EGZ2011Viscosity,Murakami2026Boundary,Lu2013Viscosity}).

\begin{proposition}\label{prop:nonpolar-positivity}
If $\varphi\in\mathcal A_f(\theta)$, then $\Phi_{n-j}^f(\theta_\varphi,\chi)$ is a positive current for $0\le j\le d$.
\end{proposition}

\begin{proof}
We work on a coordinate ball and write $\theta_\varphi=\ii\partial\bar\partial v$ and $\chi=\ii\partial\bar\partial h$, where we may assume that $h\ge0$. Suppose first that $v$ is continuous. We choose concentric coordinate balls $B'\Subset B$ and define
\[
 v^\varepsilon(x)=\sup_{y\in\overline B}\left(v(y)-\frac{|x-y|^2}{2\varepsilon}\right),\qquad x\in B'.
\]
If $y$ attains the supremum, then $|x-y|\le C\sqrt\varepsilon$. By the uniform continuity of $v$, we obtain $0\le |x-y|^2/(2\varepsilon)\le v(y)-v(x)\le\delta_\varepsilon$, with $\delta_\varepsilon\searrow0$.

We now take a smooth function $q$ touching $v^\varepsilon$ from above at $x$. Choose a complex linear map $L$ with $L^*\chi_x=\chi_y$ and $L=I+O(|x-y|)$. Set $\tau(z)=x+L(z-y)$ and $p(z)=|\tau(z)-z|^2/(2\varepsilon)$. For $z$ near $y$, we have
\[
 \ii\partial\bar\partial(q\circ\tau+p)(y)
 =L^*(\ii\partial\bar\partial q(x))+\ii\partial\bar\partial p(y)
 \in\overline C_f(\chi_y),
 \qquad
 v(z)-p(z)\le v^\varepsilon(\tau(z))\le q(\tau(z)).
\]
Thus $q\circ\tau+p$ touches $v$ from above at $y$. Since $\ii\partial\bar\partial p(y)\le C|x-y|^2\chi_y/\varepsilon\le C\delta_\varepsilon\chi_y$, it follows that $\ii\partial\bar\partial q(x)+C\delta_\varepsilon\chi_x\in\overline C_f(\chi_x)$. Thus $w=v^\varepsilon+C\delta_\varepsilon h$ is admissible and decreases to $v$ as $\varepsilon\searrow0$.

Fix $\varepsilon>0$. Since $w$ is semiconvex, we may assume that $w$ is psh by considering the translated cone. Let $w_\delta=w*\rho_\delta$ denote the local convolution. For $0<\eta<1$, set $\omega_0=\eta\chi_x$, so that $\omega_0\le\chi$ near $x$. Since the eigenvalues of $(\ii\partial\bar\partial w)_{\rm ac}$ are nonnegative, we have $\ii\partial\bar\partial w_\delta(x)\in\overline C_f(\omega_0)$, hence $\eta^{-1}\ii\partial\bar\partial w_\delta(x)\in\overline C_f(\chi_x)$. Thus, on each $U\Subset B'$ and for all sufficiently small $\delta$,
\[
 \Phi_{n-j}^f(\eta^{-1}\ii\partial\bar\partial w_\delta,\chi)\ge0,
 \qquad 0\le j\le d.
\]
Letting first $\delta\searrow0$, then $\eta\nearrow1$, and finally $\varepsilon\searrow0$, we obtain $\Phi_{n-j}^f(\ii\partial\bar\partial v,\chi)\ge0$ by Bedford-Taylor continuity.

For bounded admissible $v$, we can find continuous admissible $v_l\searrow v$ locally by \cite[Theorem~2.5]{HLP2016}, and the result also follows by letting $l\to\infty$. Finally, choose $b>0$ with $b\chi\in C_f(\chi)$ and set $v_k=\max\{v,bh-k\}$. Since $v_k$ is locally bounded and admissible, it follows from locality \cite[Lemma~2.7]{LN19Pure} that, for $0\le j\le d$,
\[
 \mathbf 1_{\{v>bh-k\}}\Phi_{n-j}^f(\ii\partial\bar\partial v,\chi)
 =\mathbf 1_{\{v>bh-k\}}\Phi_{n-j}^f(\ii\partial\bar\partial v_k,\chi)\ge0.
\]
Letting $k\to\infty$ proves the claim, since the non-$d$-polar products put no mass on $\{v=-\infty\}$.
\end{proof}

We now establish the maximum principle.

\begin{lemma}\label{lem:maximum-principle}
If $u,v\in\mathcal A_f(\theta)$, then
\[
 \mu_f(\theta,\max\{u,v\})
 \ge \mathbf 1_{\{u\ge v\}}\mu_f(\theta,u)
    +\mathbf 1_{\{u<v\}}\mu_f(\theta,v).
\]
\end{lemma}

\begin{proof}
Suppose first that $u,v$ are bounded. Set $w_\varepsilon=\max\{u+\varepsilon,v\}$.
By locality on $\{u+\varepsilon>v\}$ and $\{u+\varepsilon<v\}$, and since
$\mu_f(\theta,w_\varepsilon)\ge0$, we obtain
\[
 \mu_f(\theta,w_\varepsilon)
 \ge \mathbf 1_{\{u+\varepsilon>v\}}\mu_f(\theta,u)
    +\mathbf 1_{\{u+\varepsilon<v\}}\mu_f(\theta,v).
\]
Letting $\varepsilon\searrow0$ and using Bedford-Taylor continuity, we obtain
the desired inequality. In general, choose $b>0$ with $b\chi\in C_f(\chi)$. On a coordinate ball,
take a smooth function $h$ satisfying $\theta+\ii\partial\bar\partial h=b\chi$, and set
\[
 u_k=\max\{u,h-k\},\qquad v_k=\max\{v,h-k\}.
\]

Applying the locally bounded case to $u_k,v_k$ and using locality, we then obtain
the desired inequality on $\{\min(u,v)>h-k\}$. Finally, we let $k\to\infty$ and extend
across $\{u=-\infty\}\cup\{v=-\infty\}$, which proves the inequality.
\end{proof}

\begin{lemma}[{\cite[Proposition~3.12]{Lu2013Degenerate}}]\label{lem:weighted-coefficients}
Let $u_j,u\in\mathrm{SH}_d(X,\chi)$ be uniformly bounded and suppose that $u_j$ converges to $u$ in capacity. Let $h_j$ be uniformly bounded quasi-continuous functions converging to $h$ in capacity. Then, for $0\le p\le d$,
\[
 h_j(\chi+\ii\partial\bar\partial u_j)^p\wedge\chi^{n-d}
 \longrightarrow h(\chi+\ii\partial\bar\partial u)^p\wedge\chi^{n-d}
\]
weakly as currents.
\end{lemma}

The following result is classical and also holds for signed measures.

\begin{lemma}\label{prop:tv-closure}
Let $\theta_j$ be smooth closed $(1,1)$-forms converging smoothly to $\theta$. Let $u_j\to u$ in capacity with $u_j/A,u/A\in\mathrm{SH}_d(X,\chi)$ for some $A>0$ independent of $j$. If the measures $\mu_j=\mu_f(\theta_j,u_j)$ converge in total variation to a finite signed measure $\nu$, then $\nu=\mu_f(\theta,u)$.
\end{lemma}

\begin{proof}
Since $u_j\to u$ in capacity, the precompactness then implies that $\sup_Xu_j$ is uniformly bounded. Hence, up to extracting a subsequence, we may assume that $u_j\to u$ outside a $d$-polar set. Fix $k\ge1$ and set $h_k(t)=\min\{1,\max\{t+k,0\}\}$. By locality,
\[
 h_k(u_j)\mu_j
 =h_k(u_j)\mu_f(\theta_j,\max\{u_j,-k\}).
\]
For fixed $k$, the functions $\max\{u_j,-k\}$ are uniformly bounded and converge in capacity to $\max\{u,-k\}$. Applying Lemma~\ref{lem:weighted-coefficients} with the quasi-continuous weights $h_k(u_j)$, which converge in capacity to $h_k(u)$, and using locality again, we see that the right-hand side converges weakly to $h_k(u)\mu_f(\theta,u)$.

Since $\mu_j\to\nu$ in total variation and each $\mu_j$ puts no mass on $d$-polar sets, neither does $|\nu|$. Along the chosen subsequence, $u_j\to u$ outside a $d$-polar set. Hence $h_k(u_j)\to h_k(u)$ in $L^1(|\nu|)$ by dominated convergence, and
\[
 \|h_k(u_j)\mu_j-h_k(u)\nu\|_{\rm TV}
 \le \|\mu_j-\nu\|_{\rm TV}
 +\int_X |h_k(u_j)-h_k(u)|\,d|\nu|
 \longrightarrow0.
\]
It follows that $h_k(u)\nu=h_k(u)\mu_f(\theta,u)$. The desired equality follows as $k\to\infty$.
\end{proof}

\section{Relative uniform estimates}\label{sec:relative}

In this section, we prove a relative $L^\infty$ estimate in $f$-semistable classes with $Q>0$.

For $d<n$, the smooth approximating solutions are provided by
Fang-Ma~\cite[Theorem~1.1]{fangmaspectral2026}.
For smooth solvability in degree $d=n$, see Lin~\cite[Theorem~1.4]{Lin2023Solvability}.

\subsection{Smooth solvability}

{%
\begin{proposition}[{\cite{fangmaspectral2026,Lin2023Solvability}}]\label{prop:strict-equation}
Suppose that $[\theta]$ contains a smooth representative in $C_f(\chi)$.
Let $h\in C^\infty(X)$ be positive and satisfy
\[
 \int_X h\,\chi^n=\int_X\Phi_n^f(\theta,\chi).
\]
Then there is a unique $u\in C^\infty(X)$ such that
\[
 F(\theta+\ii\partial\bar\partial u)=h,\qquad
 \theta+\ii\partial\bar\partial u\in C_f(\chi),\qquad \sup_X u=0.
\]
\end{proposition}

}

For later application, we shall consider the following equation whose right-hand side depends on the potential.

\begin{lemma}\label{lem:proper-equation}
Suppose that $[\theta]$ contains a smooth representative in $C_f(\chi)$.
Let $H\in C^\infty(\mathbb R)$ satisfy $H'>0$, $H''\ge0$, and
$\lim_{x\to\pm\infty}H(x)=\pm\infty$. Then there is a unique
$q\in C^\infty(X)$ such that
\[
 \log F(\theta+\ii\partial\bar\partial q)=H(q),\qquad
 \theta+\ii\partial\bar\partial q\in C_f(\chi).
\]
\end{lemma}

\begin{proof}
By Proposition~\ref{prop:strict-equation}, there is a smooth representative
$\omega_*=\theta+\ii\partial\bar\partial q_*\in C_f(\chi)$ such that
\[
 F(\omega_*)=r_0,\qquad
 r_0=\frac{\int_X\Phi_n^f(\theta,\chi)}{\int_X\chi^n}>0.
\]
Set $h=e^{H/d}$ and $\Omega_s=(1-s)\omega_*+s\theta$. 
consider the set of $s\in[0,1]$ such that the equation
\[
 F^{1/d}(\eta_s)=h(q_s),\qquad
 \eta_s=\Omega_s+\ii\partial\bar\partial q_s\in C_f(\chi)
\]
is solvable.
At $s=0$ the constant $q_0=c:=H^{-1}(\log r_0)$ is a solution.
The linearized operator is elliptic with zeroth-order coefficient
$-h'(q_s)<0$, so the set of solvable parameters is open.

Since $\eta_s=\omega_*+\ii\partial\bar\partial(q_s-sq_*)$, comparison
with $\omega_*$ at the maximum and minimum of the function $q_s-sq_*$ yields
\[
 c+s(q_*-\sup_X q_*)\le q_s\le c+s(q_*-\inf_X q_*).
\]
In particular, $\|q_s\|_{C^0}$ is bounded independently of $s$.

For the second-order estimate, we use Sz\'ekelyhidi's test
function~\cite[Proposition~13]{szekelyhidi2018}, with gradient term
$\phi(|\partial q_s|_\chi^2)$ and zeroth-order term $\psi(q_s-sq_*)$.
The subsolution form is $\Omega_s+\ii\partial\bar\partial(sq_*)=\omega_*$
for every $s$. At the point under consideration, choose a $\chi$-unitary frame in which
$\eta_s$ is diagonal and $(\eta_s)_{1\bar1}$ is its largest eigenvalue,
denoted by $\lambda_1$. Since $h'>0$ and
$h''\ge0$, the terms obtained by differentiating the right-hand side
satisfy
\begin{align*}
 \operatorname{Re}\langle\partial(h(q_s)),\partial q_s\rangle_\chi
   &=h'(q_s)|\partial q_s|_\chi^2\ge0,\\
 (h(q_s))_{1\bar1}
   &=h'(q_s)\bigl(\lambda_1-(\Omega_s)_{1\bar1}\bigr)
     +h''(q_s)|(q_s)_1|^2.
\end{align*}
Since $\Omega_s$ is uniformly bounded, both terms are nonnegative when
$\lambda_1$ is sufficiently large. The subsolution form $\omega_*$ is fixed,
and the $C^0$ estimate for $q_s$ keeps $h(q_s)$ in a fixed compact
subinterval of $(0,\infty)$. Sz\'ekelyhidi's maximum principle argument~\cite[proof of Proposition~13]{szekelyhidi2018} therefore yields
\[
 \sup_X\lambda_1\le C(1+\|\partial q_s\|_\infty^2)
\]
with $C$ independent of $s$. By \eqref{eq:root-trace},
\[
 |\eta_s-c_{d-1}\chi|_\chi
 \le \operatorname{tr}_\chi(\eta_s-c_{d-1}\chi)
 \le n(\lambda_1-c_{d-1}).
\]
Subtracting $\Omega_s$ from $\eta_s$, we obtain
\begin{equation}\label{eq:fixed-hessian}
 \|\ii\partial\bar\partial q_s\|_\infty
 \le C(1+\|\partial q_s\|_\infty^2).
\end{equation}

For the gradient estimate, we use the blow-up argument. If
$R_j=\|\partial q_{s_j}\|_\infty\to\infty$, the right-hand side of the
rescaled polynomial equation is $R_j^{-2d}e^{H(q_{s_j})}$, which tends
uniformly to zero by the $C^0$ estimate. The same Liouville theorem
then yields a contradiction.

It follows from \eqref{eq:fixed-hessian} that the eigenvalues of $\eta_s$
are uniformly bounded. Since $F(\eta_s)=e^{H(q_s)}$ has a uniform positive
lower bound and $F$ vanishes on $\partial C_f$, the eigenvalues remain
in a compact subset of $C_f$. The $C^{2,\gamma}$ estimate of
\cite{TWWY15} now implies closedness. The solution at $s=1$ is smooth by
elliptic regularity, and uniqueness follows from the maximum principle
since $H'>0$.
\end{proof}

{%
As applications of Lemma~\ref{lem:proper-equation}, we will take $H(x)=\beta x$,
$\beta>0$, for the envelope approximation in Lemma~\ref{lem:continuous-envelope},
and $H(x)=x+\log\Lambda$ in \eqref{eq:exponential-family}.}

By Definition~\ref{def:f-positive-nef}, for each $t>0$ we can choose a smooth $b_t$ with $\theta+t\chi+\ii\partial\bar\partial b_t\in C_f(\chi)$.
 Proposition~\ref{prop:strict-equation} therefore applies in the class $\alpha+t\beta$.

\begingroup\emergencystretch=4em
To prove the continuity of $V_t$, we employ Berman's PDE method~\cite{Berman2019},
as in Lu-Nguyen~\cite[Theorem~3.2]{LN14Pure}.
\par\endgroup

\begin{lemma}\label{lem:continuous-envelope}
For fixed $t>0$, the envelope $V_t$ is continuous.
\end{lemma}

\begin{proof}
{%
Normalize the potential $b_t$ chosen above by $\sup_X b_t=0$, and set
\[
 m_t=\min_X F(\theta_t+\ii\partial\bar\partial b_t)>0.
\]
For $\beta>1$, Lemma~\ref{lem:proper-equation} provides a smooth admissible
solution $v_\beta$ of}
$F(\theta_t+\ii\partial\bar\partial v_\beta)=\ee^{\beta v_\beta}$.
{%
At a maximum point $x_\beta$ of $v_\beta$, write
$\eta=(\theta_t+\ii\partial\bar\partial v_\beta)(x_\beta)\in C_f(\chi_{x_\beta})$
and $\xi=-\ii\partial\bar\partial v_\beta(x_\beta)\ge0$.
By Lemma~\ref{lem:root-trace}(1)-(2), we obtain
\[
 \ee^{\beta v_\beta(x_\beta)}=F(\eta)
 \le F(\eta+\xi)=F(\theta_t)(x_\beta)\le K_t,\qquad
 K_t:=\max\{1,\sup_X|F(\theta_t)|\}.
\]
Thus $v_\beta-\beta^{-1}\log K_t\leq0$ and is admissible,
so by the definition of $V_t$ in \eqref{eq:admissible-envelope},
$v_\beta\le V_t+\beta^{-1}\log K_t$.

To bound $v_\beta$ from below, set
\[
 c_\beta=\beta^{-1}\log^+\frac{\beta^d}{m_t},\qquad
 w=(1-\beta^{-1})V_t+\beta^{-1} b_t-c_\beta.
\]
{Suppose $M=\max_X(w-v_\beta)>0$ and choose a point $x$ where the maximum is attained. Then}
\[
 \phi=\frac{v_\beta+M-\beta^{-1} b_t+c_\beta}{1-\beta^{-1}}
\]
touches $V_t$ from above at $x$. Proposition~\ref{prop:admissibility-equivalence}
shows that $(\theta_t+\ii\partial\bar\partial\phi)(x)\in\overline C_f(\chi_x)$.
Since
\[
 \theta_t+\ii\partial\bar\partial v_\beta
 =(1-\beta^{-1})(\theta_t+\ii\partial\bar\partial\phi)
   +\beta^{-1}(\theta_t+\ii\partial\bar\partial b_t),
\]
the concavity of $F^{1/d}$ implies
$\ee^{\beta v_\beta(x)}\ge\beta^{-d} F(\theta_t+\ii\partial\bar\partial b_t)(x)
\ge\beta^{-d} m_t$. Since $V_t,b_t\le0$, we have $w\le-c_\beta$ by definition.
At $x$, the preceding lower bound and $w(x)>v_\beta(x)$ therefore give
\[
 \beta^{-d}m_t
 \le \ee^{\beta v_\beta(x)}
 < \ee^{\beta w(x)}
 \le \ee^{-\beta c_\beta}
 = \min\{1,\beta^{-d}m_t\}
 \le \beta^{-d}m_t.
\]
This contradiction proves $w\le v_\beta$ on $X$. Together with the upper bound, we have}
\[
 V_t-\frac{\|b_t\|_\infty+d\log\beta+\log^+(1/m_t)}{\beta}
 \le v_\beta\le V_t+\frac{\log K_t}{\beta}.
\]
{%
For fixed $t>0$, these bounds give
$\sup_X|v_\beta-V_t|\to0$ as $\beta\to\infty$.
Since $v_{\beta}$ are continuous functions, $V_t$ is continuous.}
\end{proof}

\subsection{Auxiliary Monge-Ampère equations}

Write $W=\int_X\chi^n$ and $\rd\nu=\chi^n/W$. {Let $\psi$ be the quasi-psh
potential of Proposition~\ref{prop:common-locus}, normalized by
$\sup_X\psi=-1$. Since $\ii\partial\bar\partial\psi\ge-B\chi$ for some
$B>0$, the Skoda-Zeriahi integrability theorem
 provides constants $b,J>0$ such that
\begin{equation}\label{eq:anchor-integrability}
 \int_X\ee^{-b\psi}\,\rd\nu\le J.
\end{equation}}

{%
We use the following uniform estimate for complex Monge-Amp\`ere equations with $L\log^q L$ densities \cite{Kol1998,GuedjLuEnvelopes1,guedj2025,GPT2023}.

\begin{lemma}\label{lem:ordinary-entropy}
Let $q>n$ and let $H>0$ be smooth with $\int_X H\,\rd\nu=1$.
Let $v$ be the $\chi$-psh solution of
\[
 (\chi+\ii\partial\bar\partial v)^n=H\chi^n,\qquad \sup_X v=0.
\]
Then
\[
 \|v\|_\infty\le C(1+A)^{1/q},
 \qquad A=\int_X H\bigl(\log(1+H)\bigr)^q\,\rd\nu,
\]
where $C$ depends only on $(X,\chi)$ and $q$.
\end{lemma}
}

For $q>n$, let the function $\Psi_q(s)=s^{n/d}\log^q(1+s)$ for $s\ge0$.
The corresponding Luxemburg norm is
\begin{equation}\label{eq:luxemburg-norm}
 \|h\|_{\Psi_q,\nu}
 =\inf\left\{\lambda>0:\int_X\Psi_q(|h|/\lambda)\,\rd\nu\le1\right\}.
\end{equation}
\begingroup
\begin{proposition}\label{prop:relative-bound}
Fix $t_0>0$ and $q>n$. Let $u_t$ be smooth solutions, $0<t\le t_0$, such that
\[
 \sup_Xu_t=0,\qquad
 \theta_t+\ii\partial\bar\partial u_t\in C_f(\chi),\qquad
 F(\theta_t+\ii\partial\bar\partial u_t)=h_t>0.
\]
If
\[
 {\sup_{0<t\le t_0}\left\|h_t\right\|_{\Psi_q,\nu}\le M<\infty},
\]
then
\[
 V_t-{C}\le u_t\le V_t,
 \qquad
 \sup_{0<t\le t_0}\|u_t-V_t\|_{L^\infty(X)}\le{C}.
\]
{%
Here $C>0$ is a uniform constant depending only on $(X,\chi)$, $\theta$,
$f$, the choice of the barrier $\psi$ from Proposition~\ref{prop:common-locus},
$q$, $t_0$, and $M$.}
\end{proposition}
\endgroup

\begin{proof}
{With $r=n/d$ as in the definition of $\Psi_q$, the bound
$\|h_t\|_{\Psi_q,\nu}\le M$ implies
\[
 \int_X\Psi_q(h_t/M)\,\rd\nu\le1.
\]
Since $x\mapsto\log(1+x)$ is increasing, concave and vanishes at zero,
\[
 \log(1+Mz)\le\max\{1,M\}\log(1+z).
\]
Hence
\begin{align*}
 \int_X h_t^r\log^q(1+h_t^r)\,\rd\nu
 &\le r^q\int_X\Psi_q(h_t)\,\rd\nu\le r^q M^r\max\{1,M\}^q
       \int_X\Psi_q(h_t/M)\,\rd\nu\\
 &
\le r^q M^r\max\{1,M\}^q=:K.
\end{align*}
For $0<t\le t_0$, the trace inequality and Green's formula give
$\int_X(-u_t)\,\rd\nu\le L_0$.} {%
By \eqref{eq:anchor-integrability}, $\nu(\{\psi<-R\})\le J\ee^{-bR}$.
On $\{h_t^r>\ee^{bR/2}\}$ we have $\log(1+h_t^r)\ge bR/2$. Hence, for $R\ge1$,
\begin{equation}\label{eq:anchor-tail}
\begin{aligned}
 \int_{\{\psi<-R\}}h_t^r\,\rd\nu
 &\le \ee^{bR/2}\nu(\{\psi<-R\})
       +\int_{\{h_t^r>\ee^{bR/2}\}}h_t^r\,\rd\nu\\
 &\le J\ee^{-bR/2}
       +\left(\frac{2}{bR}\right)^q
        \int_X h_t^r\log^q(1+h_t^r)\,\rd\nu
 \le B_0R^{-q}.
\end{aligned}
\end{equation}}

{%
For fixed $t\in(0,t_0]$, let
\[
 E_l=\{V_t-u_t>l\},\qquad m(l)=\int_{E_l}h_t^r\,\rd\nu.
\]
Since $V_t$ is continuous by Lemma~\ref{lem:continuous-envelope} and $u_t$
is smooth, $\{V_t-u_t\ge l+\varepsilon\}$ is a compact subset of $E_l$
for every $\varepsilon>0$. Choose $\eta\in C^\infty_c(E_l)$ with
$0\le\eta\le1$ and $\eta=1$ on $\{V_t-u_t\ge l+\varepsilon\}$.
For $m=m(l)>0$, define
\[
 N=\int_X(\eta h_t^r+m)\,\rd\nu,\qquad
 H=\frac{\eta h_t^r+m}{N},\qquad m\le N\le2m.
\]}
{For $0<m\le1/2$ we have $m\le N\le1$ and
$\log(1+H)\le\log(2+h_t^r)+\log(1/N)$.
Since $\int_XH\,\rd\nu=1$ and $H\le(1+h_t^r)/N$, it follows that
\begin{align*}
 \int_X H\log^q(1+H)\,\rd\nu
 &\le \frac{2^{q-1}}{N}\int_X(1+h_t^r)\log^q(2+h_t^r)\,\rd\nu
       +2^{q-1}\log^q(1/N)\\
 &\le \frac{C_q(1+K)}{m}.
\end{align*}
Here $(1+h_t^r)\log^q(2+h_t^r)\le C_q(1+h_t^r\log^q(1+h_t^r))$ and
$N\log^q(1/N)\le C_q$ for $0<N\le1$.}

{Let $v$ be the $\chi$-psh solution of the auxiliary Monge-Ampère equation
\[
 (\chi+\ii\partial\bar\partial v)^n
 =\frac{\eta h_t^r+m}{N}\,\chi^n,\qquad \sup_X v=0.
\]}
Choose $\delta_*>0$ with $\theta+\ii\partial\bar\partial\psi-\delta_*\chi\in\overline C_f$ off $E_{\mathrm{nK}}^{f}$, set $a=\min(1,\delta_*/2)$, and fix $\Lambda$ with $(\Lambda a)^d\ge2$. {Let $s=\Lambda N^{1/n}$.} For small enough $m$, $s\le1/2$ and Lemma~\ref{lem:ordinary-entropy} yields
\[
 s\|v\|_\infty\le Cm^{1/n-1/q},\qquad s\le Cm^{1/n}.
\]
{%
For a real $(1,1)$-form $\gamma$, recall
\(E_d(\gamma)=\sigma_d(\lambda_\chi(\gamma))/\binom nd
        =\gamma^d\wedge\chi^{n-d}/\chi^n.
\)
For $\eta\in C_f(\chi)$ and $\gamma>0$, Lemma~\ref{lem:root-trace}(1)-(2) yields
\[
 G(\eta+\gamma)-G(\eta)
 =\int_0^1 D_\gamma G(\eta+\tau\gamma)\,\rd\tau
 \ge E_d(\gamma)^{1/d}.
\]
Using continuity, we obtain for
$\eta\in\overline C_f(\chi)$ and $\gamma\ge0$,
\begin{equation}\label{eq:polynomial-insertion}
 F(\eta+\gamma)\ge\bigl(F(\eta)^{1/d}+E_d(\gamma)^{1/d}\bigr)^d
 \ge F(\eta)+E_d(\gamma).
\end{equation}}
Consider ${w=(1-s)V_t+s\psi+sa v}$.
Since ${\psi\le V_t}$ and $v\le0$, one has ${w\le V_t}$.
If $q_0$ touches $w$ from above at a finite point outside $E_{\mathrm{nK}}^{f}$, then
\[
 {q_{V_t}}=\frac{q_0-s\psi-sa v}{1-s}
\]
touches ${V_t}$ from above. By the upper-test characterization in Proposition~\ref{prop:admissibility-equivalence}, we have ${\eta_{V_t}}:=\theta_t+\ii\partial\bar\partial {q_{V_t}}\in\overline C_f$. Moreover,
\begin{align*}
 \theta_t+\ii\partial\bar\partial q_0
   &=\eta_0+sa(\chi+\ii\partial\bar\partial v),\\*
 \eta_0
   &=(1-s)\eta_{V_t}+s(\theta_t+\ii\partial\bar\partial\psi-a\chi)\in C_f.
\end{align*}
{%
Indeed, by convexity of $\overline C_f$,
\[
 \eta_0-s(\delta_*-a+t)\chi
 =(1-s)\eta_{V_t}+s(\theta+\ii\partial\bar\partial\psi-\delta_*\chi)
 \in\overline C_f.
\]
Since $s(\delta_*-a+t)>0$, $\eta_0\in C_f$ by \eqref{eq:centered-branch}.
Write
$\chi_v=\chi+\ii\partial\bar\partial v$, we obtain from
\eqref{eq:polynomial-insertion} and Maclaurin's inequality
\[
 F(\theta_t+\ii\partial\bar\partial q_0)
 \ge (sa)^d E_d(\chi_v)
 \ge (sa)^d E_n(\chi_v)^{d/n}
 = (sa)^d H^{d/n}.
\]

Suppose $D=\max_X(w-u_t-l-\varepsilon)>0$, attained at $x$.
Since $w=-\infty$ on the pole set of $\psi$, $x$ lies outside it. Also,
\[
 V_t(x)-u_t(x)\ge w(x)-u_t(x)=l+\varepsilon+D,
 \qquad\text{so }\eta(x)=1.
\]
The function $q_0=u_t+l+\varepsilon+D$ touches $w$ from above at $x$.
At $x$, substituting $s=\Lambda N^{1/n}$ into the preceding inequality yields
\[
 h_t(x)=F(\theta_t+\ii\partial\bar\partial q_0)(x)
 \ge (sa)^d H(x)^{d/n}
 =(\Lambda a)^d(h_t(x)^r+m)^{d/n}\ge2h_t(x).
\]
This contradicts $h_t(x)>0$. Hence $w\le u_t+l+\varepsilon$.} On ${\{V_t-\psi\le R\}}$ it follows that
\[
 {V_t-u_t}\le l+\varepsilon+C\bigl(m^{1/n}R+m^{1/n-1/q}\bigr).
\]
The complement lies in $\{\psi<-R\}$, since ${V_t\le0}$. By \eqref{eq:anchor-tail}, choosing
\[
 R=m^{-(1/q+1/n)/2},\qquad
 \varepsilon=m^\tau,\qquad\tau=\tfrac12(1/n-1/q)>0
\]
produces
\begin{equation}\label{eq:relative-iteration}
 m(l+C_0m(l)^\tau)\le B_0m(l)^\lambda,
 \qquad\lambda=\tfrac12(1+q/n)>1.
\end{equation}

To start the iteration, we observe that for every $T>1$,
\[
 m(l)\le T\nu(E_l)+{\frac K{\log^q(1+T)}}
 \le\frac{TL_0}{l}+{\frac K{\log^q(1+T)}}.
\]
Choose $m_0>0$ so small that
\[
 2m_0\le1,\qquad \Lambda(2m_0)^{1/n}\le\tfrac12,
 \qquad B_0m_0^{\lambda-1}\le\tfrac12,
\]
and so that all preceding small-mass requirements hold. Then choose
$T,l_0$ so that $m(l_0)\le m_0$.

{%
To iterate \eqref{eq:relative-iteration}, define
\[
 l_{j+1}=l_j+C_0m(l_j)^\tau,\qquad j\ge0.
\]
Since $m(l_j)\le m_0$ and $B_0m_0^{\lambda-1}\le\tfrac12$, we obtain
\[
 m(l_{j+1})\le B_0m(l_j)^{\lambda-1}m(l_j)
 \le \tfrac12 m(l_j),\qquad m(l_j)\le 2^{-j}m_0.
\]
If $m(l_j)=0$, then $l_{j+1}=l_j$ and these bounds still hold.
Consequently, for every $j$,
\[
 l_j\le l_0+C_0m_0^\tau\sum_{k=0}^{\infty}2^{-k\tau}
 =l_0+\frac{C_0m_0^\tau}{1-2^{-\tau}}=:L.
\]
Thus $m(L)\le m(l_j)\to0$. Since $E_L=\{V_t-u_t>L\}$ is open and $h_t^r>0$,
$\int_{E_L}h_t^r\,\rd\nu=0$ implies $E_L=\varnothing$. Hence $u_t\ge V_t-L$.
Finally, $u_t$ is admissible and $\sup_Xu_t=0$, so $u_t\le V_t$ by
\eqref{eq:admissible-envelope}.}
\end{proof}

\section{Existence of weak solutions}\label{sec:limits}

{%
Throughout this section, $f$ is a monic right-Noetherian polynomial, the class $\alpha=[\theta]$ is $f$-semistable, and
$Q=\int_X\Phi_n^f(\alpha,\beta)>0$.
For $\rho$ as in Proposition~\ref{prop:weak-orlicz}, we will first construct
$u\in\mathcal A_f(\theta)$ satisfying $\mu_f(\theta,u)=Q\rho\,\chi^n$ on $X$.
In Section~6.2, we will show that
$F\bigl((\theta+\ii\partial\bar\partial u)_{\rm ac}\bigr)=Q\rho$ a.e.
}

\subsection{Existence and uniqueness}\label{subsec:nonnegative-uniqueness}

\begingroup
\begin{proposition}\label{prop:weak-orlicz}
{Fix $q>n$ and let $\rho\ge0$ satisfy $\int_X\rho\,\chi^n=1$ and
$\|\rho\|_{\Psi_q,\nu}<\infty$, where $\Psi_q(s)=s^{n/d}\log^q(1+s)$ and
$\rd\nu=\chi^n/\int_X\chi^n$.
Then the equation
\[
 \mu_f(\theta,u)=Q\rho\,\chi^n
\]
admits a solution $u\in\mathcal A_f(\theta)$ with $\sup_Xu=0$ and
$V_\theta-C\le u\le V_\theta$ for some $C>0$.}
\end{proposition}
\begin{proof}
{%
Choose smooth positive densities $\rho_j$ with $\int_X\rho_j\,\chi^n=1$ such that
$\rho_j\to\rho$ in $L^{\Psi_q}(X,\nu)$.
Up to extracting a subsequence, we may assume that $\sum_j\|\rho_j-\rho\|_{L^1}<\infty$.

Fix $0<t_j\le t_0$ with $t_j\searrow0$, and write
$g_j=Q_{t_j}\rho_j$ and $g=Q\rho$. By
Proposition~\ref{prop:strict-equation} there is a normalized smooth
solution $u_j$ with $F(\theta_{t_j}+\ii\partial\bar\partial u_j)=g_j$.
Since $Q_{t_j}\to Q$ and $\sup_j\|\rho_j\|_{\Psi_q,\nu}<\infty$,
the norms $\|g_j\|_{\Psi_q,\nu}$ are uniformly bounded. By Proposition~\ref{prop:relative-bound}, the solutions satisfy
\[
 V_{t_j}-C\le u_j\le V_{t_j}\le0
\]
with $C$ independent of $j$. Let $h_j=\inf_{k\ge j}g_k$. Then
$0\le h_j\nearrow g$ a.e. and in $L^1(X,\chi^n)$.

For $N\ge j$, set $w_{j,N}=\max_{j\le k\le N}u_k$ and
$w_j=(\sup_{k\ge j}u_k)^*$. Since $\theta_{t_j}-\theta_{t_k}=(t_j-t_k)\chi\ge0$
for $k\ge j$, each $u_k$ is admissible in the background $\theta_{t_j}$ and
$\mu_f(\theta_{t_j},u_k)\ge g_k\chi^n$ by~\eqref{eq:polynomial-insertion}.
Lemma~\ref{lem:maximum-principle} implies that
\[
 \mu_f(\theta_{t_j},w_{j,N})\ge\min_{j\le k\le N}g_k\,\chi^n\ge h_j\chi^n.
\]
For fixed $j$, we have $w_{j,N}\le0$, and $w_{j,N}$ increases to $w_j$
outside a $d$-polar set. Thus $w_{j,N}\to w_j$ in capacity. Lemma~\ref{lem:weighted-coefficients} yields
\[
 \mu_j:=\mu_f(\theta_{t_j},w_j)\ge h_j\chi^n.
\]
Since $w_j$ is bounded, integration by parts yields $\mu_j(X)=\int_X\Phi_n^f(\theta_{t_j},\chi)=Q_{t_j}$.
Thus
\begin{equation}\label{eq:tail-tv}
 \begin{aligned}
 \|\mu_j-g\chi^n\|_{\rm TV}
 &\le (\mu_j-h_j\chi^n)(X)+\int_X(g-h_j)\chi^n\\
 &= Q_{t_j}-Q+2\int_X(g-h_j)\chi^n\longrightarrow0.
 \end{aligned}
\end{equation}

Then $w_j\searrow u$ with $\sup_X u=0$. Since $V_\theta-C\le w_j\le V_{t_j}$ and
$V_{t_j}\searrow V_\theta$, we have $V_\theta-C\le u\le V_\theta$.
By Corollary~\ref{cor:admissible-closure}(i), $u\in\mathcal A_f(\theta)$.
Since $w_j\to u$ in $d$-capacity, it follows from Lemma~\ref{prop:tv-closure}
and~\eqref{eq:tail-tv} that $\mu_f(\theta,u)=g\chi^n$.}
\end{proof}
\endgroup

Following \cite{BEGZ10}, we establish the comparison and domination principles for $\mu_f$.

\begin{lemma}\label{lem:weak-comparison}
Let $u,v\in\mathcal A_f(\theta)$ satisfy $u=v+O(1)$ on $X$. Then
\[
 \int_{\{u<v\}}\mu_f(\theta,v)\le\int_{\{u<v\}}\mu_f(\theta,u).
\]
\end{lemma}

\begin{proof}
Set $w=\max\{u,v\}$ and $E=\{u<v\}$. Since $w=u+O(1)$, the total masses agree
by \cite[Theorem~3.7]{LN19Pure}. Applying Lemma~\ref{lem:maximum-principle}, we obtain
\[
 \mu_f(\theta,v)(E)
 \le \mu_f(\theta,w)(X)-\mu_f(\theta,u)(X\setminus E)
 =\mu_f(\theta,u)(E).
\]
\end{proof}

\begin{lemma}\label{lem:weak-domination}
Let $u,v\in\mathcal A_f(\theta)$ and assume that $v=V_\theta+O(1)$.
If $u\le v$ $\mu_f(\theta,v)$-a.e.,
then $u\le v$ on $X$.
\end{lemma}

\begin{proof}
Choose $\delta>0$ and $\psi\in\mathcal A_f(\theta-\delta\chi)$ as in
Proposition~\ref{prop:common-locus}. After subtracting a constant, we may assume
that $\psi\le v$. For $0<t<1$, set $u_t=(1-t)u+t\psi$ and $E_t=\{v<u_t\}$.
Since $u_t\in\mathcal A_f(\theta-t\delta\chi)$, expanding in $t\delta\chi$
and applying Proposition~\ref{prop:nonpolar-positivity}, we obtain
$\mu_f(\theta,u_t)\ge(t\delta)^d\chi^n$.

Since $u\le v+C$ for some $C>0$, we have $\max\{v,u_t\}=v+O(1)$.
Applying Lemma~\ref{lem:weak-comparison} to $v$ and $\max\{v,u_t\}$
and using locality on $E_t$, we obtain
\[
 (t\delta)^d\int_{E_t}\chi^n\le\int_{E_t}\mu_f(\theta,u_t)
 \le\int_{E_t}\mu_f(\theta,v)=0,
\]
where the last equality follows from $E_t\subset\{v<u\}$.
Letting $t\searrow0$, we obtain $u\le v$ on $X$.
\end{proof}

\begin{lemma}\label{lem:constructed-mixed-root}
Let $u,\widetilde u$ be solutions constructed as in Proposition~\ref{prop:weak-orlicz}
with the same sequence $t_j\searrow0$ and normalized densities
$\rho,\widetilde\rho$, respectively. Write $g=Q\rho$ and
$\widetilde g=Q\widetilde\rho$. For $0<t<1$,
\[
 \mu_f(\theta,(1-t)u+t\widetilde u)
 \ge \bigl((1-t)g^{1/d}+t\widetilde g^{1/d}\bigr)^d\,\chi^n .
\]
\end{lemma}

\begin{proof}
Use the notation of Proposition~\ref{prop:weak-orlicz} for $u$, and write
$\widetilde u_j,\widetilde g_j,\widetilde h_j,\widetilde w_{j,N},\widetilde w_j$
for the corresponding objects for $\widetilde u$. For $N\ge j$,
\[
 (1-t)w_{j,N}+t\widetilde w_{j,N}
 =\max_{j\le k,l\le N}\bigl((1-t)u_k+t\widetilde u_l\bigr).
\]

\noindent\hspace*{2em}In the common background $\theta_{t_j}$, the concavity of $F^{1/d}$ implies
\[
 \mu_f\bigl(\theta_{t_j},(1-t)u_k+t\widetilde u_l\bigr)
 \ge \bigl((1-t)g_k^{1/d}+t\widetilde g_l^{1/d}\bigr)^d\,\chi^n .
\]

\noindent\hspace*{2em}Applying Lemma~\ref{lem:maximum-principle} to the maximum over $k,l$, we obtain
\[
 \mu_f\bigl(\theta_{t_j},(1-t)w_{j,N}+t\widetilde w_{j,N}\bigr)
 \ge \bigl((1-t)h_j^{1/d}+t\widetilde h_j^{1/d}\bigr)^d\,\chi^n .
\]
For fixed $j$, the potentials $(1-t)w_{j,N}+t\widetilde w_{j,N}$ are uniformly bounded
and converge in capacity to $(1-t)w_j+t\widetilde w_j$.
Letting $N\to\infty$ and applying Lemma~\ref{lem:weighted-coefficients},
we obtain the same lower bound for $\mu_f(\theta_{t_j},(1-t)w_j+t\widetilde w_j)$.

By Proposition~\ref{prop:weak-orlicz}, $w_j\searrow u$ and $\widetilde w_j\searrow\widetilde u$
with uniform bounds on compact subsets of $X\setminus E_{\mathrm{nK}}^f$.
Since $\theta_{t_j}\to\theta$ smoothly and $h_j\nearrow g$, $\widetilde h_j\nearrow\widetilde g$ a.e., Bedford-Taylor continuity implies the desired inequality on $X\setminus E_{\mathrm{nK}}^f$. Since both sides put no mass on $E_{\mathrm{nK}}^f$, the inequality then holds on $X$.
\end{proof}

\begin{theorem}\label{thm:nonnegative-minimal-uniqueness}
Under the hypotheses of Proposition~\ref{prop:weak-orlicz}, the equation
$\mu_f(\theta,u)=Q\rho\,\chi^n$ has a unique solution $u$ such that $u=V_\theta+O(1)$ and $\sup_Xu=0$.
\end{theorem}

\begin{proof}
Let $v$ be the solution obtained in Proposition~\ref{prop:weak-orlicz}, and let $u$ be another normalized solution with minimal singularities. Set $g=Q\rho$. Assume by contradiction that $v-u$ is not constant $(g\chi^n)$-a.e. We can then choose $s\in\mathbb R$ such that $0<\int_E g\chi^n<Q$, where $E=\{v>u+s\}$.
Following Dinew-Lu~\cite{DinewLu2015}, choose $\eta>0$ so that $\int_E(g+\eta)\,\chi^n<Q$, and set
\[
 \widetilde g=\frac{Q\,\mathbf1_E(g+\eta)}{\int_E(g+\eta)\,\chi^n}.
\]
Thus $\widetilde g>g$ on $E$, $\int_X\widetilde g\,\chi^n=Q$, and
$\widetilde g\le C(g+1)$. Apply Proposition~\ref{prop:weak-orlicz}
with $\widetilde\rho=\widetilde g/Q$ and the same $t_j$ as for $v$.
Denote the resulting solution by $\widetilde v$. After subtracting a constant,
we may assume that $\widetilde v\le v$, since both have minimal singularities.

For $0<t<1$, let $v_t=(1-t)v+t\widetilde v$ and $E_t=\{v_t>u+s\}\subset E$.
Applying Lemma~\ref{lem:constructed-mixed-root} to $v,\widetilde v$ and
Lemma~\ref{lem:weak-comparison} to $u+s,v_t$, we obtain
\[
 \int_{E_t}\bigl((1-t)g^{1/d}+t\widetilde g^{1/d}\bigr)^d\,\chi^n
 \le \int_{E_t}g\,\chi^n .
\]
Since $\widetilde g>g$ on $E_t$, we have $\int_{E_t}g\,\chi^n=0$.
As $t\searrow0$, we have $v_t\nearrow v$ a.e., so
$\int_E g\,\chi^n=\lim_{t\searrow0}\int_{E_t}g\,\chi^n=0$, contrary to the choice of $s$.

Thus $u=v+c$ $(g\,\chi^n)$-a.e. for some $c\in\mathbb R$.
By Lemma~\ref{lem:weak-domination}, $u\le v+c$ and $v+c\le u$ on $X$.
Hence $c=\sup_Xu-\sup_Xv=0$.
\end{proof}

\subsection{The absolutely continuous part}

{%
Let $u$ be the solution of Proposition~\ref{prop:weak-orlicz} and
$T=\theta+\ii\partial\bar\partial u$. We prove that
$F(T_{\rm ac})=Q\rho$ $\chi^n$-a.e.}

\begin{lemma}\label{lem:ac-limit}
{Let $\theta_j\to\theta$ smoothly. For each $j$, suppose that
$v_{j,N}\to w_j$ in $L^1$ as $N\to\infty$ and that
$w_j\searrow u\not\equiv-\infty$. Assume that there is an $A>0$ such that
$v_{j,N}/A,w_j/A,u/A\in\mathrm{SH}_d(X,\chi)$ for all $j,N$. Write $T_{j,N}=\theta_j+\ii\partial\bar\partial v_{j,N}$
and $T=\theta+\ii\partial\bar\partial u$. If
\[
 (T_{j,N})_{\rm ac}\in\overline C_f(\chi),\qquad
 F((T_{j,N})_{\rm ac})\ge h_j,\qquad 0\le h_j\nearrow h
 \quad\text{a.e.},
\]
then $T_{\rm ac}\in\overline C_f(\chi)$ and $F(T_{\rm ac})\ge h$
a.e.}
\end{lemma}

\begin{proof}
{%
On a coordinate ball, fix a smooth $\chi$-unitary coframe and let $\rho_\delta$ be a radial mollifier. By convexity of $\overline C_f$, we then have $(T_{j,N})_{\rm ac}*\rho_\delta\in\overline C_f$ on each smaller ball for $\delta$ sufficiently small. Also, as $(T_{j,N})_{\rm s}\ge0$ by Lemma~\ref{lem:measure-hessian}, we deduce that $T_{j,N}*\rho_\delta\in\overline C_f$. Moreover, it follows from the monotonicity and concavity of $G=F^{1/d}$ that
\[
 G(T_{j,N}*\rho_\delta)
 \ge G\bigl((T_{j,N})_{\rm ac}*\rho_\delta\bigr)
 \ge G\bigl((T_{j,N})_{\rm ac}\bigr)*\rho_\delta
 \ge h_j^{1/d}*\rho_\delta .
\]
Letting first $N\to\infty$ and then $j\to\infty$, and using $h_j\nearrow h$, we obtain
\[
 T*\rho_\delta\in\overline C_f,\qquad
 G(T*\rho_\delta)\ge h^{1/d}*\rho_\delta .
\]
As $\delta\searrow0$, the Lebesgue differentiation theorem implies that
$T*\rho_\delta\to T_{\rm ac}$ and $h^{1/d}*\rho_\delta\to h^{1/d}$
a.e. This implies that $T_{\rm ac}\in\overline C_f$ and $F(T_{\rm ac})\ge h$ a.e.}
\end{proof}

\begin{proposition}\label{lem:constructed-ac-lower}
Let $u$ be the solution constructed in Proposition~\ref{prop:weak-orlicz}. Then
\begin{equation}\label{eq:weak-ac-lower}
 F\bigl((\theta+\ii\partial\bar\partial u)_{\rm ac}\bigr)
 \ge Q\rho\quad\text{a.e. on }X.
\end{equation}
\end{proposition}
\begin{proof}
To replace the role of the maximum principle in Lemma~\ref{lem:maximum-principle} as used in the proof of
Proposition~\ref{prop:weak-orlicz}, we argue as follows. Fix $j\le N$ and write
$m=\min_{j\le k\le N}g_k$ and $p=F((\theta_{t_j}+\ii\partial\bar\partial w_{j,N})_{\rm ac})$.
For each rational $c>0$, the set $\{m>c\}$ is open. There each $u_k$, $j\le k\le N$,
is $(f-c)$-admissible in the background $\theta_{t_j}$, since
$C_{f-c}=C_f\cap\{F>c\}$ by~\eqref{eq:partial-characterization}.
Their maximum is $(f-c)$-admissible there by Corollary~\ref{cor:admissible-closure}(ii).
Thus $\{0\le p<c<m\}$ has zero volume for each rational $c>0$. Since $p\ge0$ a.e. and
\[
 \{0\le p<m\}=\bigcup_{c\in\mathbb Q_{>0}}\{0\le p<c<m\}
\]
has Lebesgue measure zero, we thus have $p\ge m\ge h_j$ a.e. Applying Lemma~\ref{lem:ac-limit} with $h=Q\rho$ then proves the desired inequality.
\end{proof}

\begin{proposition}\label{lem:ac-product-bound}
Let $u\in\mathcal A_f(\theta)$ and set $T=\theta+\ii\partial\bar\partial u$. Then
\[
 \mu_f(\theta,u)\ge F(T_{\rm ac})\,\chi^n.
\]
\end{proposition}
\begin{proof}
Write $G=F^{1/d}$. For $\eta\in\overline C_f(\chi)$ and $\gamma\ge0$,
Lemma~\ref{lem:root-trace}(2) and the AM-GM inequality imply
\[
 G(\eta+\gamma)-G(\eta)\ge E_d(\gamma)^{1/d}\ge(\det\nolimits_\chi\gamma)^{1/n}.
\]

Work on a relatively compact coordinate ball $B$ and fix $c>0$.
We choose a smooth function $h$ on $B$ with $\theta+\ii\partial\bar\partial h=b\chi\in C_{f-c}(\chi)$,
and set $u_k=\max\{u,h-k\}$.

Let $K\subset B\cap\{F(T_{\rm ac})>c\}$ be compact. Choose open sets $K\subset U_j\Subset B$ with
$|U_j\setminus K|\to0$. For a constant $M>0$ to be chosen below, let $b_j$ be the
convex Alexandrov solution of
\[
 \operatorname{MA}_{\mathbb R}(b_j)=M\mathbf1_{U_j\setminus K}\,\rd x
 \quad\text{in }B,\qquad b_j=0\quad\text{on }\partial B,
\]
where $\rd x$ denotes Euclidean volume. By Alexandrov's maximum principle,
\[
 \|b_j\|_{L^\infty(B)}
 \le C_B\bigl(M|U_j\setminus K|\bigr)^{1/(2n)}\longrightarrow0.
\]
Moreover, $\det_{\mathbb R}D^2b_j=M$ a.e. on $U_j\setminus K$
(see \cite[Section~2]{DePhilippisFigalli2014}). Set $\gamma_j=(\ii\partial\bar\partial b_j)_{\rm ac}$. Since $b_j$ is convex, at these points we have
\[
 (\det\nolimits_\chi\gamma_j)^{1/n}
 \ge c_B(\det\nolimits_{\mathbb R}D^2b_j)^{1/(2n)}
 =c_BM^{1/(2n)}.
\]
Choose $M$ so that $c_BM^{1/(2n)}>c^{1/d}$. The estimate for $G$ then implies
$F(T_{\rm ac}+\gamma_j)>c$ there. On $K$ the same
inequality follows from $F(T_{\rm ac})>c$ and $\gamma_j\ge0$. Thus $u+b_j$ is
$f$-admissible on $B$ and $(f-c)$-admissible on $U_j$, since $C_{f-c}=C_f\cap\{F>c\}$.

Set $v_{j,k}=\max\{u+b_j,h-k\}$. By Corollary~\ref{cor:admissible-closure}(ii),
$v_{j,k}$ is $f$-admissible on $B$ and therefore $(f-c)$-admissible on $U_j$. It follows from Proposition~\ref{prop:nonpolar-positivity} that, on $U_j$, we have
\[
 \mu_f(\theta,v_{j,k})-c\chi^n=\mu_{f-c}(\theta,v_{j,k})\ge0.
\]
For fixed $k$, the potentials $v_{j,k}$ are locally uniformly bounded and
converge to $u_k$ uniformly. By Bedford-Taylor continuity, we have $\mu_f(\theta,v_{j,k})\rightharpoonup\mu_f(\theta,u_k)$ on $B$. Moreover, since $K$ is compact, we have
\[
 \mu_f(\theta,u_k)(K)\ge\limsup_{j\to\infty}\mu_f(\theta,v_{j,k})(K)
 \ge c\int_K\chi^n.
\]
The inner regularity then yields $\mu_f(\theta,u_k)\ge c\chi^n$ on $\{F(T_{\rm ac})>c\}$. Finally, let $E_k=\{u>h-k\}$. Then, by pluri-fine locality, on $\{F(T_{\rm ac})>c\}$ we have
\[
 \mathbf1_{E_k}\mu_f(\theta,u)
 =\mathbf1_{E_k}\mu_f(\theta,u_k)
 \ge c\,\mathbf1_{E_k}\chi^n.
\]
Letting $k\to\infty$, we obtain
$\mu_f(\theta,u)\ge c\chi^n$ on $\{F(T_{\rm ac})>c\}$.
Now denote by $p$ the density of the absolutely continuous part of $\mu_f(\theta,u)$ with respect to $\chi^n$. Then $p\ge c$ a.e. on $\{F(T_{\rm ac})>c\}$ for every rational $c>0$. Similar to the proof of Proposition~\ref{lem:constructed-ac-lower}, we also have that $p\ge F(T_{\rm ac})$ a.e. on $B$. Thus $\mu_f(\theta,u)\ge p\,\chi^n\ge F(T_{\rm ac})\,\chi^n$ on $B$.
Since $B$ is arbitrary, we complete the proof.
\end{proof}

For the solution $u$ constructed in Proposition~\ref{prop:weak-orlicz}, Proposition~\ref{lem:ac-product-bound}
and $\mu_f(\theta,u)=Q\rho\chi^n$ imply $F(T_{\rm ac})\le Q\rho$ a.e. Together with
\eqref{eq:weak-ac-lower}, this proves $F(T_{\rm ac})=Q\rho$ a.e. on $X$.

\begin{proposition}\label{prop:actual-family}
Let $t_j\searrow0$ and let $u_j$ be normalized smooth admissible solutions of
$F(\theta_{t_j}+\ii\partial\bar\partial u_j)=g_j$. Assume that
\[
 \sup_j\|g_j\|_\infty<\infty,\qquad g_j\to\rho\ \text{in }L^1(X,\chi^n),
 \qquad V_{t_j}-C\le u_j\le V_{t_j}.
\]
Then $u_j\to u$ in $L^1(X)$, where $u$ is the unique normalized weak solution such that $\mu_f(\theta,u)=\rho\chi^n$. Moreover,
$V_\theta-C\le u\le V_\theta$, $u$ is locally bounded on $X\setminus E_{\mathrm{nK}}^f$, and
\[
 F\bigl((\theta+\ii\partial\bar\partial u)_{\rm ac}\bigr)=\rho\quad\text{a.e.}
\]
\end{proposition}
\begin{proof}
Since $\int_X\rho\,\chi^n=\lim_j Q_{t_j}=Q$, we may apply Proposition~\ref{prop:weak-orlicz}
to $\rho/Q$. Let $u$ be the resulting solution.
Choose $A>0$ so that $u_j/A\in\mathrm{SH}_d(X,\chi)$ for all $j$. By $L^1$ compactness,
every subsequence has a further subsequence along which $u_j\to v$ in $L^1(X)$
and $g_j\to\rho$ a.e.

Along this subsequence, Hartogs' lemma implies that $w_j=(\sup_{k\ge j}u_k)^*\searrow v$.
The proof of Proposition~\ref{prop:weak-orlicz} therefore shows that $v$ is a normalized
weak solution of $\mu_f(\theta,v)=\rho\chi^n$ with $V_\theta-C\le v\le V_\theta$.
By Theorem~\ref{thm:nonnegative-minimal-uniqueness}, $v=u$.
Thus $u$ is the only $L^1$ cluster point, and $u_j\to u$ in $L^1(X)$.

Since $u\ge V_\theta-C\ge\psi-C$ and $\psi$ is smooth on $X\setminus E_{\mathrm{nK}}^f$,
$u$ is locally bounded there. Moreover, combining \eqref{eq:weak-ac-lower} and Proposition~\ref{lem:ac-product-bound}, we obtain
$F\bigl((\theta+\ii\partial\bar\partial u)_{\rm ac}\bigr)=\rho$ a.e. on $X$, finishing the proof.
\end{proof}

\section{Regularity outside the null locus}\label{sec:regularity}

Fix the potential $\psi\le-1$ from Proposition~\ref{prop:common-locus}, with pole set $Z=E_{\mathrm{nK}}^{f}$ and
\[
 \theta+\ii\partial\bar\partial\psi-\delta\chi\in\overline C_f,
 \qquad |\partial\psi|^2\le C\ee^{-\kappa\psi}
 \quad\text{on }X\setminus Z
\]
for fixed $\delta,C,\kappa>0$.
All calculations with this potential take place outside $Z$. 

Fix $0<\tau<1$. To obtain a lower bound for $u_t$ in terms of
$\tau\psi$, choose $\varepsilon_\tau>0$ so small that
$(1-\tau)\varepsilon_\tau\le\tau\delta/4$.
Since $\alpha=[\theta]$ is $f$-semistable, the class
$\alpha+\varepsilon_\tau[\chi]$ admits a smooth representative
$\omega_\tau\in C_f(\chi)$. Write
\[
 \omega_\tau=\theta+\varepsilon_\tau\chi
             +\ii\partial\bar\partial b_\tau,
 \qquad b_\tau\in C^\infty(X,\mathbb R).
\]
Then
\[
 \theta+\ii\partial\bar\partial(\tau\psi+(1-\tau)b_\tau)
       -\tau\delta\chi/2
 =\tau(\theta+\ii\partial\bar\partial\psi-\delta\chi)
 +(1-\tau)\omega_\tau+
       (\tau\delta/2-(1-\tau)\varepsilon_\tau)\chi
 \in\overline C_f.
\]

{Now set} $w_\tau=\tau\psi+(1-\tau)(b_\tau-\sup_X b_\tau)$. Then $w_\tau\le0$ is a candidate for the envelope $V_\theta$.
It follows from
Proposition~\ref{prop:relative-bound} that for any family of normalized solutions with uniformly bounded densities, $u_t\ge V_t-C$ with $C$
independent of $t$. Since  $V_t\ge V_\theta$, we obtain
\begin{equation}\label{eq:small-pole-bound}
 u_t\ge V_t-C\ge V_\theta-C\ge w_\tau-C\ge\tau\psi-C_\tau,
\end{equation}
where $C_\tau=C+(1-\tau)\operatorname{osc}_X b_\tau$ is uniform in $t$
for each fixed $\tau\in(0,1)$.

{Taking $\tau=1/2$ in \eqref{eq:small-pole-bound}, we have
$\psi-u_t\le\psi/2+C_{1/2}$. Thus for each $a>0$ the weight
$\ee^{a(\psi-u_t)}$ converges to 0 along $Z$, uniformly in $t$.}

\subsection{Weighted gradient estimates}

For a smooth $u\le0$, write
\[
 \theta_u=\theta+t\chi+\ii\partial\bar\partial u,\quad
 \eta_u=\theta_u-c_{d-1}\chi,\quad
 \Omega=\theta+(t-c_{d-1})\chi,\quad
 P=|\partial u|_\chi^2.
\]
Let $\lambda$ be the largest eigenvalue of $\eta_u$. For $a>0$, set
\[
 \zeta=\ee^{a(\psi-u)},\qquad
 K_a=\sup_{X\setminus Z}\zeta P,\qquad
 H_a=\sup_{X\setminus Z}\zeta(\lambda+aP/2).
\]

We first obtain an estimate of the complex Hessian in terms of the gradient.

\begin{proposition}\label{prop:weighted-hessian}
Suppose $\theta_u\in C_f$, $u\ge\psi/2-C_0$, and either
\[
 F(\theta_u)^{1/d}=h(u),\qquad
 0<h\le H,\quad h',h''\ge0,
\]
or
\[
 F(\theta_u)^{1/d}=h(x),\qquad
 0<h\le H,\quad
 |\ii\partial\bar\partial h|+\frac{|\partial h|^2}{h}
       \le C_h\ee^{-r\psi}.
\]
There exists $a_0>0$ such that, for every compact interval $I\subset[a_0,\infty)$, we have
\begin{equation}\label{eq:weighted-hessian}
 H_a\le C_I(1+K_a),\qquad a\in I.
\end{equation}
\end{proposition}

\begin{proof}
We write $G=F^{1/d}$ and $g=G(\theta_u)$. In a $\chi$-unitary frame diagonalizing $\eta_u$, we denote its eigenvalues by $\lambda=\lambda_1\ge\cdots\ge\lambda_n$ and the eigenvalue derivatives of $G$ at $\theta_u$ by $r_i$. We use the notation
\[
 Lw=\sum_i r_iw_{i\bar i},\qquad W=\sum_i r_i,\qquad \Gamma(w)=\sum_i r_i|w_i|^2.
\]
Since $\eta_u\in\Gamma_d$ by \eqref{eq:centered-branch}, Lemma~\ref{lem:root-trace}(2) implies that
\[
 W=D_\chi G(\theta_u)\ge1,\qquad Lu=D_{\eta_u}G(\theta_u)-D_\Omega G(\theta_u)\ge-C_\Omega W.
\]
Moreover, by the concavity of $G$ and the choice of $\Xi=\theta+\ii\partial\bar\partial\psi-\delta\chi\in\overline C_f$, we obtain
\[
 L(\psi-u)=D_{\Xi-\theta_u}G(\theta_u)+(\delta+t)W\ge G(\Xi)-g+(\delta+t)W.
\]
Together with \eqref{eq:root-trace} and $g\le H$, these inequalities imply
\begin{equation}\label{eq:weighted-structure}
 W\ge1,\qquad \lambda\le nHW^{d-1},\qquad Lu\ge-C_\Omega W,\qquad L(\psi-u)\ge\delta W-H.
\end{equation}

For $a\in I$, we set $b=a/2$ and consider
\[
 Q=\log S+a(\psi-u)+J(u),\qquad S=\lambda+bP,\qquad J(u)=\frac{j}{1-u},
\]
where $0<j\le\min\{1,a_0/2\}$. Since $S$ is bounded for each smooth solution and
\[
 \zeta=e^{a(\psi-u)}\le e^{aC_0}e^{a\psi/2},
\]
the function $e^Q=e^J\zeta S$ attains its maximum at a point $x_0\in X\setminus Z$.

We choose K\"ahler normal coordinates at $x_0$ in which $\eta_u$ is diagonal. If $\lambda$ has multiplicity $k$, we perturb $\eta_u$ by
\[
 B=\varepsilon\,\ii\sum_{j=2}^{k}dz^j\wedge d\bar z^j,\qquad \widetilde\eta=\eta_u-B,\qquad 0<\varepsilon\le1.
\]
The largest eigenvalue $\widetilde\lambda$ of $\widetilde\eta$ with respect to $\chi$ is simple at $x_0$, and $\widetilde\lambda\le\lambda$, with equality there. We may therefore replace $\lambda$ by $\widetilde\lambda$ when applying the maximum principle at $x_0$.

Following the calculation in \cite[(6.10)-(6.11)]{GuoSong2026Gradient} (see also \cite[(3.9)-(3.11)]{HuangZhangZhang2026} for the $k$-Hessian equation), we write $\mu_j=\widetilde\eta_{j\bar j}(x_0)$. The eigenvalue differentiation formula \cite[(64), (70)]{szekelyhidi2018} then gives
\[
 \partial_{\bar p}\partial_p\widetilde\lambda
 =(\nabla_{\bar p}\nabla_p\widetilde\eta)_{1\bar1}
 +\sum_{j>1}\frac{|\nabla_p\widetilde\eta_{j\bar1}|^2+|\nabla_p\widetilde\eta_{1\bar j}|^2}{\lambda-\mu_j}.
\]
We also set $T=\nabla_1\eta_u=\nabla_1\theta_u$. Since $\nabla B(x_0)=0$ and $\eta_u$ is closed, we have
\[
 \partial_i\widetilde\lambda=\nabla_i(\eta_u)_{1\bar1}=\nabla_1(\eta_u)_{i\bar1}=T_{i\bar1}.
\]
Multiplying the preceding second-derivative identity by $r_p$ and summing over $p$, we obtain
\begin{align*}
 L\widetilde\lambda-\sum_p r_p(\nabla_{\bar p}\nabla_p\widetilde\eta)_{1\bar1}
 &=\sum_{p=1}^n\sum_{j>1}\frac{r_p\bigl(|\nabla_p\widetilde\eta_{j\bar1}|^2+|\nabla_p\widetilde\eta_{1\bar j}|^2\bigr)}{\lambda-\mu_j}\\
 &\ge\sum_{i>k}\frac{r_1|T_{i\bar1}|^2}{\lambda-\lambda_i},
\end{align*}
where we used $\mu_i=\lambda_i$ for $i>k$.
In these coordinates, $R_{p\bar p1\bar1}=-\partial_{\bar p}\partial_p\chi_{1\bar1}$.
Commuting covariant derivatives of $\eta_u$, we have
\begin{align*}
 (\nabla_{\bar p}\nabla_p\widetilde\eta)_{1\bar1}
 &= (\nabla_{\bar1}\nabla_1\eta_u)_{p\bar p}
  +(\lambda-\lambda_p)R_{p\bar p1\bar1}
  -(\nabla_{\bar p}\nabla_pB)_{1\bar1},\\
 \left|\sum_p r_p(\lambda-\lambda_p)R_{p\bar p1\bar1}\right|
 &\le C\lambda W,\qquad
 \left|\sum_p r_p(\nabla_{\bar p}\nabla_p B)_{1\bar1}\right|
 \le C\varepsilon W,
\end{align*}
where $|\eta_u|_\chi\le\operatorname{tr}_\chi\eta_u\le n\lambda$
and $|\nabla^2B(x_0)|_\chi\le C\varepsilon$.
Thus
\begin{subequations}
\begin{equation}\label{eq:eigen-variation-bound}
 L\widetilde\lambda
 \ge \sum_p r_p(\nabla_{\bar1}\nabla_1\eta_u)_{p\bar p}
    +\sum_{i>k}\frac{r_1|T_{i\bar1}|^2}{\lambda-\lambda_i}
    -C(\lambda+\varepsilon)W .
\end{equation}
With $T^*$ denoting the Hermitian adjoint of $T$, differentiation of
$G(\theta_u)=g$ and the concavity formula \cite[(67)]{szekelyhidi2018}
give
\begin{equation}\label{eq:eigen-differentiated-bound}
\begin{aligned}
 \sum_p r_p(\nabla_{\bar1}\nabla_1\eta_u)_{p\bar p}
 &=g_{1\bar1}-D^2G_{\theta_u}[T,T^*]\ge g_{1\bar1}
      +\sum_{i>k}\frac{(r_i-r_1)|T_{i\bar1}|^2}{\lambda-\lambda_i}.
\end{aligned}
\end{equation}
By \eqref{eq:eigen-variation-bound}-\eqref{eq:eigen-differentiated-bound} and using
$T_{i\bar1}=\partial_i\widetilde\lambda$, we obtain
\begin{equation*}\tag{\theparentequation}\label{eq:weighted-eigen}
 L\widetilde\lambda\ge g_{1\bar1}
   +\sum_{i>k}\frac{r_i|\partial_i\widetilde\lambda|^2}{\lambda-\lambda_i}
   -C_R(1+\lambda)W.
\end{equation*}
\end{subequations}

For the gradient term, write
$U_i=\sum_q|u_{qi}|^2$ and $V_i=\sum_q|u_{i\bar q}|^2$.
Differentiating the equation once, we find
$\sum_i r_i\nabla_q u_{i\bar i}=g_q-\sum_i r_i\nabla_q\Omega_{i\bar i}$.
The K\"ahler case of the Bochner formula
\cite[(3.9)]{GuoSong2026Gradient} then reads
\begin{align*}
 LP&=\sum_i r_i(U_i+V_i)+2\operatorname{Re}\langle\partial g,\partial u\rangle_\chi
     -2\operatorname{Re}\sum_{i,q}r_i\overline{u_q}\nabla_q\Omega_{i\bar i}
     +E_R,\\
 |E_R|&\le CPW,\qquad
 \left|\sum_{i,q}r_i\overline{u_q}\nabla_q\Omega_{i\bar i}\right|
 \le C|\nabla\theta|_\chi\sqrt P\,W.
\end{align*}
In particular,
\begin{equation}\label{eq:weighted-bochner}
 LP\ge\sum_i r_i(U_i+V_i)+2\operatorname{Re}\langle\partial g,\partial u\rangle_\chi
       -C_B(P+\sqrt P)W.
\end{equation}

Let $M=e^{Q(x_0)}$. By $0<J\le j$ and $M=e^J\zeta S$ at $x_0$,
\begin{equation}\label{eq:weighted-max-ratios}
 H_a\le M\le e^jH_a,\qquad
 \frac{P}{S}=\frac{e^J\zeta P}{M}\le\frac{e^jK_a}{M}.
\end{equation}
It suffices to consider $M\ge A_I(1+K_a)$, where $A_I$ can be
increased below. Until $j$ is chosen, all constants $C_I$ are uniform
for $0<j\le\min\{1,a_0/2\}$. At $x_0$,
\[
 \frac{bP}{S}\le\frac{be^j}{A_I}\le\frac14,\qquad
 \lambda\ge\frac{3S}{4},\qquad
 \frac1S\le \frac{C_I}{M}e^{a\psi/2}.
\]
Choose $a_0/2>\max\{2\kappa,r\}$, with $r=0$ for $g=h(u)$.
By \eqref{eq:barrier-gradient}, $|\partial\psi|^2\le C e^{-\kappa\psi}$.
Since $(1-\psi)^3e^{(a_0/2-\kappa)\psi}$ is bounded for $\psi\le-1$,
\begin{equation}\label{eq:weighted-barrier-ratios}
 \frac{|\partial\psi|^2}{S}\le\frac{C_I}{M(1-\psi)^3},\qquad
 \frac{1+|\partial\psi|^4+e^{-r\psi}}{S}\le\frac{C_I}{M}.
\end{equation}

At $x_0$,
\begin{equation}\label{eq:weighted-log-identity}
 L\log S=\frac{L\widetilde\lambda+bLP}{S}-\frac{\Gamma(S)}{S^2}.
\end{equation}
Write $p_i=u_i$, $\xi_i=\psi_i$, $j_1=J'(u)$ and $\alpha=a-j_1$.
Since $\partial_iQ=0$, we have
\begin{equation}\label{eq:weighted-first-derivative}
 z_i:=\frac{\partial_i S}{S}=\alpha p_i-a\xi_i,\qquad
 \frac{\Gamma(S)}{S^2}=\sum_i r_i|z_i|^2 .
\end{equation}
With $t_i=\sum_q u_{qi}\overline{p_q}$ and
$e_i=\sum_q\Omega_{i\bar q}p_q$, we have
\begin{align*}
 \partial_iP&=t_i+\lambda_i p_i-e_i,\qquad |e_i|\le C_\Omega\sqrt P,\\
 \partial_i\widetilde\lambda+b t_i
 &=\partial_i S-b\lambda_i p_i+b e_i=S z_i-b\lambda_i p_i+b e_i .
\end{align*}

For $i>k$, let $D_i=\lambda-\lambda_i>0$.
Expanding $|\partial_i\widetilde\lambda-D_i z_i|^2$ and
$\sum_q|u_{qi}-p_qz_i|^2$, we find
\[
 \frac{|\partial_i\widetilde\lambda|^2}{D_i}
 \ge 2\operatorname{Re}(\partial_i\widetilde\lambda\,\overline{z_i})
       -D_i|z_i|^2,\qquad
 U_i\ge 2\operatorname{Re}(t_i\overline{z_i})-P|z_i|^2 .
\]
Since $D_i+bP=S-\lambda_i$ and $z_i-2bp_i=-j_1p_i-a\xi_i$,
\begin{equation}\label{eq:low-index-square}
\begin{aligned}
 \frac{|\partial_i\widetilde\lambda|^2}{SD_i}
      +\frac bS U_i-|z_i|^2
 &\ge \frac{2}{S}\operatorname{Re}
       \bigl((\partial_i\widetilde\lambda+bt_i)\overline{z_i}\bigr)
       -\left(1+\frac{D_i+bP}{S}\right)|z_i|^2\\
 &=-\frac{\lambda_i}{S}\operatorname{Re}
       \bigl((j_1p_i+a\xi_i)\overline{z_i}\bigr)
       +\frac{2b}{S}\operatorname{Re}(e_i\overline{z_i}).
\end{aligned}
\end{equation}
By $|z_i|\le a(|p_i|+|\xi_i|)$ and Young's inequality,
\begin{align*}
 \frac{|\lambda_i|}{S}|j_1p_i+a\xi_i|\,|z_i|
 &\le \frac{b\lambda_i^2}{4S}
     +\frac{C_I}{S}\bigl(j_1^2|p_i|^4+|p_i|^2|\xi_i|^2+|\xi_i|^4\bigr),\\
 \frac{2b}{S}|e_i|\,|z_i|&\le\frac{C_I}{S}(P+|\xi_i|^2),\\
 \frac bS V_i &\ge\frac bS|\lambda_i-\Omega_{i\bar i}|^2
          \ge\frac{b\lambda_i^2}{2S}-\frac{C_I}{S}.
\end{align*}
As $|p_i|^2/S\le P/S\le e^j/A_I$, it follows that for $i>k$,
\[
 \frac{|\partial_i\widetilde\lambda|^2}{SD_i}
 +\frac bS(U_i+V_i)-|z_i|^2
 \ge -C_I j_1^2|p_i|^2-\frac{C_I}{S}
      \bigl(|p_i|^2|\xi_i|^2+|\xi_i|^4+P+|\xi_i|^2+1\bigr).
\]
For $i\le k$, we have $\lambda_i=\lambda\ge 3S/4$.
Increase $A_I$ so that $S\ge4\sup_X|\Omega|_\chi$. Then
\[
 \frac bS V_i\ge\frac bS|\lambda-\Omega_{i\bar i}|^2\ge\frac{bS}{4},
 \qquad
 \frac{|z_i|^2}{S}
 \le 2a^2\frac{P+|\partial\psi|^2}{S}
 \le\frac{C_I(1+K_a)}{M}\le\frac b4 .
\]
We therefore have $bV_i/S-|z_i|^2\ge0$ for $i\le k$. Setting $\mathcal D=g_{1\bar1}+2b\operatorname{Re}\langle\partial g,\partial u\rangle_\chi$ and substituting \eqref{eq:weighted-eigen} and \eqref{eq:weighted-bochner} into \eqref{eq:weighted-log-identity}, we obtain
\begin{align*}
 L\log S\ge{}&\frac{\mathcal D}{S}
 +\sum_{i>k}r_i\frac{|\partial_i\widetilde\lambda|^2}{SD_i}
 +\sum_i r_i\left(\frac bS(U_i+V_i)-|z_i|^2\right)\\
 &-\frac{C_R(1+\lambda)+bC_B(P+\sqrt P)}{S}W.
\end{align*}
Combining the estimates for $i\le k$ and $i>k$, and using $\lambda/S\le1$, we infer that
\begin{equation}\label{eq:weighted-log-lower}
\begin{aligned}
 L\log S\ge{}&
 \frac{\mathcal D}{S}
 -C_I\left(j_1^2+\frac{|\partial\psi|^2}{S}\right)\Gamma(u)
 -C_R W\\
 &-\frac{C_I}{S}
     \bigl(1+|\partial\psi|^4+|\partial\psi|^2+P+\sqrt P\bigr)W.
\end{aligned}
\end{equation}
Here $C_R$ is the curvature constant from \eqref{eq:weighted-eigen} and is independent of $a$ and $j$. When $g=h(u)$, we have $u_{1\bar1}=\lambda-\Omega_{1\bar1}\ge0$ at $x_0$, and hence
\[
 \mathcal D=h'(u)(\lambda-\Omega_{1\bar1})+h''(u)|u_1|^2+2b h'(u)P\ge0.
\]
If $g=h(x)$, then $g_{1\bar1}\ge-C_h e^{-r\psi}$ and
$|\partial g|^2\le C_h H e^{-r\psi}$. Hence
\[
 \frac{\mathcal D}{S}
 \ge-\frac{C_h e^{-r\psi}+2b\sqrt{C_h H P}\,e^{-r\psi/2}}{S}
 \ge-\frac{C_I}{S}(e^{-r\psi}+P).
\]
Now $\sqrt P\le1+P$. By \eqref{eq:weighted-max-ratios} and
\eqref{eq:weighted-barrier-ratios},
\[
 \frac{1+|\partial\psi|^4+|\partial\psi|^2+P+\sqrt P+e^{-r\psi}}{S}
 \le\frac{C_I(1+K_a)}{M}.
\]
Since $W\ge1$, it follows from \eqref{eq:weighted-log-lower} in both cases that
\begin{equation}\label{eq:weighted-log-absorbed}
 L\log S\ge
 -\left(C_Ij_1^2+\frac{C_I}{M(1-\psi)^3}\right)\Gamma(u)
 -\left(C_R+\frac{C_I(1+K_a)}{M}\right)W .
\end{equation}

Let $j_2=J''(u)=2j/(1-u)^3$. By \eqref{eq:weighted-structure},
\begin{equation}\label{eq:weight-linearization}
\begin{aligned}
 L\bigl(a(\psi-u)+J(u)\bigr)
 &=aL(\psi-u)+j_1Lu+j_2\Gamma(u)\\
 &\ge(a\delta-C_\Omega j_1)W-aH+j_2\Gamma(u).
\end{aligned}
\end{equation}
Choose $a_0$ also so that $C_R\le a_0\delta/4$. For the given
interval $I$, first choose $j$ small, then $A_I$ large.
Indeed, $u\le0$ and $u\ge\psi/2-C_0$ imply
\[
 j_1\le j,\qquad
 \frac{j_1^2}{j_2}=\frac{j}{2(1-u)}\le\frac j2,\qquad
 j_2\ge\frac{c j}{(1-\psi)^3}.
\]
We may therefore choose $j$ and $A_I$ so that
\[
 C_\Omega j_1\le\frac{a\delta}{8},\quad C_I j_1^2\le\frac{j_2}{4},
 \qquad
 \frac{C_I(1+K_a)}{M}\le\frac{a\delta}{8},\quad
 \frac{C_I}{M(1-\psi)^3}\le\frac{j_2}{4}.
\]
Adding \eqref{eq:weighted-log-absorbed} and \eqref{eq:weight-linearization}
at $x_0$, we have
\[
 0\ge LQ\ge\frac{a\delta}{2}W-aH+\frac{j_2}{2}\Gamma(u),
 \qquad
 W\le\frac{2H}{\delta},\qquad
 \lambda\le nH\left(\frac{2H}{\delta}\right)^{d-1}.
\]
Consequently, when $M\ge A_I(1+K_a)$,
\[
 H_a\le M=e^J\zeta S\le\frac43 e^{j+aC_0}\lambda\le C_I .
\]
When $M<A_I(1+K_a)$, we have $H_a\le M<A_I(1+K_a)$.
This proves \eqref{eq:weighted-hessian}.
\end{proof}

{It remains to bound $K_a$. Guo-Song
\cite[Proposition~3.1]{GuoSong2026Gradient} estimate the gradient by
comparison with a smooth admissible function close to $u$ in $L^\infty$.
(For the family solving \eqref{eq:exponential-family}, we use a $C^1$
comparison function and need a bound linear in its weighted gradient.)}

\begin{lemma}\label{lem:gradient-transfer}
Let $u$ satisfy the hypotheses of Proposition~\ref{prop:weighted-hessian}.
For each compact interval $I\subset[a_0,\infty)$ there are constants
$\eta_I>0$, $C_{1}(I)$ and $C_2$ such that, for every $a\in I$ and every
real $v\in C^1_{\rm loc}(X\setminus Z)$ with
$\|u-v\|_{L^\infty(X\setminus Z)}\le\eta_I$,
\begin{equation}\label{eq:gradient-transfer}
 \sqrt{K_a(u)}\le C_{1}(I)+C_{2}
     \sup_{X\setminus Z}e^{a(\psi-v)/2}|\partial v|_\chi .
\end{equation}
{The constants $\eta_I$ and $C_1(I)$ may depend on $I$.} The coefficient $C_2$ depends only on the dimension and the bounds comparing $\chi$ with the Euclidean metric in the fixed finite atlas. In particular, the same $C_2$ works for all compact intervals $I\subset[a_0,\infty)$.

\end{lemma}

\begin{proof}
Fix $a\in I$ and set
\[
 D_v=\sup_{X\setminus Z}e^{a(\psi-v)/2}|\partial v|_\chi,
 \qquad \zeta=e^{a(\psi-u)}.
\]
We may assume $D_v<\infty$. For now suppose $\|u-v\|_\infty\le\eta$,
where $0<\eta\le1/\max I$. We choose $\eta=\eta_I$ at the end.
In the estimates below, $C$ depends only on the dimension and the fixed atlas.

Fix $x\in X\setminus Z$ and work in a coordinate chart from a fixed
finite atlas, using smaller charts that still cover $X$. By
\eqref{eq:barrier-gradient}, $\ell=e^{\kappa\psi/2}$ has bounded
gradient off $Z$ and extends by zero to a Lipschitz function.
On $B_{r_I\ell(x)}(x)$ we have
$|\ell(z)/\ell(x)-1|\le \operatorname{Lip}(\ell)r_I$.
Taking $r_I$ small makes the ball disjoint from $Z$ and implies that
$a|\psi(z)-\psi(x)|=(2a/\kappa)|\log(\ell(z)/\ell(x))|\le1$.
Since $u\ge\psi/2-C_0$, $a\ge a_0>2\kappa$, and $\psi\le0$, we have
\[
 \sqrt{\zeta(x)}\le e^{aC_0/2}e^{a\psi(x)/4}\le C_I\ell(x).
\]
We can therefore choose $c_I>0$ sufficiently small so that $B_\rho=B_\rho(x)$ lies in the preceding coordinate ball whenever
\begin{equation}\label{eq:transfer-radius}
 0<\rho\le\frac{c_I\sqrt{\zeta(x)}}{1+D_v}.
\end{equation}

To estimate the oscillation of $v$ on $B_\rho$, we set $V=e^{-av/2}$. By the definition of $D_v$, we have $|\partial V|_\chi\le(a/2)D_ve^{-a\psi/2}$. Integrating along the segment from $x$ to $z\in B_\rho$, and using $a|\psi-\psi(x)|\le1$ on $B_\rho$, we obtain
\[
 \left|\frac{V(z)}{V(x)}-1\right|
 \le CaD_v\rho\,e^{-a(\psi(x)-v(x))/2}
 \le CaD_v\frac{\rho}{\sqrt{\zeta(x)}},
\]
where the last inequality follows from $a|v(x)-u(x)|\le1$. After decreasing $c_I$, we may assume that the right-hand side is at most $1/4$. Taking logarithms, we then obtain
\[
 |v(z)-v(x)|\le\frac{CD_v\rho}{\sqrt{\zeta(x)}},\qquad a|v(z)-v(x)|\le2\log(4/3).
\]
Since $\operatorname{osc}_{B_\rho}u\le\operatorname{osc}_{B_\rho}v+2\eta$ and $\log\zeta=a(\psi-u)$, these estimates, together with $a\eta\le1$ and $a|\psi-\psi(x)|\le1$, imply that
\begin{equation}\label{eq:transfer-local}
 \operatorname{osc}_{B_\rho}u\le2\eta+\frac{CD_v\rho}{\sqrt{\zeta(x)}},\qquad C^{-1}\le\frac{\zeta(z)}{\zeta(x)}\le C.
\end{equation}
Both constants are independent of $I$.

We next estimate the Euclidean Laplacian of $u$ on $B_\rho$. By \eqref{eq:root-trace}, we have
\[
 |\eta_u|_\chi\le\operatorname{tr}_\chi\eta_u\le n\lambda.
\]
Since $\eta_u=\Omega+\ii\partial\bar\partial u$ and $\Delta_{\rm Euc}u=4\sum_i u_{i\bar i}$, it follows that $|\Delta_{\rm Euc}u|\le C\lambda+C_{\rm bg}$. Using $\zeta\lambda\le H_a$, $\zeta\le e^{aC_0}$, and the preceding comparison between $\zeta(z)$ and $\zeta(x)$, we conclude that
\[
 \zeta(x)\sup_{B_\rho}|\Delta_{\rm Euc}u|\le C(H_a+B_I),\qquad B_I=1+C_{\rm bg}e^{(\max I)C_0}.
\]
Applying the interior gradient estimate for Poisson's equation
\cite[Theorem~3.9]{GilbargTrudinger2001} to $u-\inf_{B_\rho}u$, we find
\[
 |\partial u(x)|_\chi\le C\left(
   \frac{\operatorname{osc}_{B_\rho}u}{\rho}
   +\rho\sup_{B_\rho}|\Delta_{\rm Euc}u|\right).
\]
Multiplying by
$\sqrt{\zeta(x)}$ and substituting the two bounds above, we obtain
\begin{equation}\label{eq:transfer-scaled}
 \sqrt{\zeta(x)}|\partial u(x)|_\chi
 \le C\left(D_v+\frac{\eta}{s}+(H_a+B_I)s\right),
 \qquad s=\frac{\rho}{\sqrt{\zeta(x)}}\le\frac{c_I}{1+D_v}.
\end{equation}

Take $\eta\le c_I$, {set} $A=H_a+B_I$, and choose
$s=\min\{\sqrt{\eta/A},c_I/(1+D_v)\}$. Then $As\le\sqrt{\eta A}$ and
\[
 \frac{\eta}{s}
 =\max\left\{\sqrt{\eta A},\frac{\eta}{c_I}(1+D_v)\right\}
 \le \sqrt{\eta A}+1+D_v .
\]
We square \eqref{eq:transfer-scaled} and take the supremum over $x$ to obtain
\begin{equation}\label{eq:transfer-before-absorption}
 K_a\le C_P(1+D_v^2)+C_P\eta(H_a+B_I),
\end{equation}
where $C_P$ depends only on the coordinate geometry.
Let $C_I$ be the constant in \eqref{eq:weighted-hessian} and take
$0<\eta_I\le\min\{1/\max I,c_I,(2C_PC_I)^{-1}\}$.
Substituting $H_a\le C_I(1+K_a)$ into
\eqref{eq:transfer-before-absorption} with $\eta=\eta_I$, we conclude that
\[
 K_a\le 2C_PD_v^2+2C_P\bigl[1+\eta_I(C_I+B_I)\bigr].
\]
Taking square roots proves \eqref{eq:gradient-transfer} with
$C_{2}=\sqrt{2C_P}$.
\end{proof}

\subsection{Stability}

Along the density path in Section~\ref{subsec:density-path} we will know that two nearby potentials are close in $L^1$. To apply Lemma~\ref{lem:gradient-transfer}, we need them to be close in $L^\infty$.

We estimate the volumes of the sublevel sets of their difference, following Dinew-Ko{\l}odziej \cite[Theorem~3.1]{DK2014}.
The lower barrier $\psi_0$ makes the estimate uniform in $t$.
\par\nopagebreak[4]

\begin{proposition}\label{prop:smooth-stability}
Fix a quasi-psh $\psi_0\le0$, smooth off a proper analytic pole set $Z$ and tending to $-\infty$ there, and fix $a,H>0$ and $p>n/d$. Suppose a smooth closed background $\theta$ satisfies
$\theta+\ii\partial\bar\partial\psi_0-a\chi\in\overline C_f$ off $Z$.
For smooth admissible $u,v$ in this background, assume
\[
 \psi_0\le v\le0,\qquad
 \|F(\theta+\ii\partial\bar\partial u)\|_{L^p}\le H.
\]
There are fixed $C,\gamma,L_*>0$ such that if
\[
 L:=\int_X(v-u)_+\chi^n\le L_*
 ,\] then \[ \sup_X(v-u)_+\le CL^\gamma.
\]
\end{proposition}

\begin{proof}
Let $p'=p/(p-1)$, with $p'=1$ for $p=\infty$, and choose
$p'<\beta<n/(n-d)$. Write $|E|:=\int_E\chi^n$ for the volume of a Borel set $E\subset X$,
 set $E_t=\{u<v-t\}$ and $V(t)=|E_t|$.
For $0<\delta\le1$, $K\ge1$, and a smooth $\chi$-psh $q\in[-1,0]$, set
\[
 s=\frac{\delta}{4(K+2a+1)},\qquad
 w=(1-s)v+s(\psi_0+a(q-1)).
 \]
 {
Write $\chi_q=\chi+\ii\partial\bar\partial q\ge0$.
On $X\setminus Z$, we have by definition of $w$
$\theta+\ii\partial\bar\partial w=\eta+sa\chi_q$, where
\[
 \eta=(1-s)(\theta+\ii\partial\bar\partial v)
      +s(\theta+\ii\partial\bar\partial\psi_0-a\chi).
\]
By the admissibility of $v$ and the barrier assumption on $\psi_0$,
both forms in this convex combination lie in $\overline C_f(\chi)$.
Thus $\eta\in\overline C_f(\chi)$ and $F(\eta)\ge0$.
Taking $\gamma=sa\chi_q$ in \eqref{eq:polynomial-insertion}, we have
\[
\begin{aligned}
 F(\theta+\ii\partial\bar\partial w)\chi^n
 &= F(\eta+sa\chi_q)\chi^n\ge \bigl(F(\eta)+(sa)^d E_d(\chi_q)\bigr)\chi^n\\
 &\ge (sa)^d\chi_q^d\wedge\chi^{n-d}
  = (sa)^d(\chi+\ii\partial\bar\partial q)^d\wedge\chi^{n-d}.
\end{aligned}
\]
}

{
Fix $0<\varepsilon<\delta/8$ and consider the regularized maximum
$U_\varepsilon:=M_\varepsilon(u,w-t-\delta/4)$.
Since $U_\varepsilon=u$ near $Z$, it extends to a smooth admissible function on $X$.
On $X\setminus E_t$ we have $w-t-\delta/4\le u-\delta/4$, so $U_\varepsilon=u$
in a neighborhood of $X\setminus E_t$. On
$B=E_{t+\delta}\cap\{\psi_0\ge-K\}$ we have
\[
 v-w=s\bigl(v-\psi_0+a(1-q)\bigr)\le s(K+2a)<\delta/4,
 \qquad w-t-\delta/4>u+\delta/2.
\]
Thus $U_\varepsilon=w-t-\delta/4$ near $B$.

Write $\mu_\phi=F(\theta+\ii\partial\bar\partial\phi)\chi^n$.
By Stokes' theorem on $X$, $\mu_{U_\varepsilon}$ and $\mu_u$ have the same total mass.
They agree off $E_t$, so the pointwise lower bound and H\"older's inequality give
\[
\begin{aligned}
 (sa)^d\int_B\chi_q^d\wedge\chi^{n-d}
 &\le \int_B\mu_{U_\varepsilon}
 \le \int_{E_t}\mu_{U_\varepsilon}
 = \int_{E_t}\mu_u\\
 &\le \|F(\theta+\ii\partial\bar\partial u)\|_{L^p(X,\chi^n)}
        V(t)^{1/p'}
 \le H V(t)^{1/p'}.
\end{aligned}
\]
}
{
The global volume-capacity
estimate of Dinew-Ko{\l}odziej \cite[Section~3, p.~238]{DK2014}
and the preceding upper bound give, for a suitable smooth $\chi$-psh $q$,
\[
 c_\beta |B|^{1/\beta}
 \le \int_B(\chi+\ii\partial\bar\partial q)^d\wedge\chi^{n-d}
 \le (sa)^{-d}H V(t)^{1/p'}.
\]
Since $(sa)^{-1}=4(K+2a+1)/(a\delta)\le C_a K/\delta$ for $K\ge1$,
we raise the estimate to the power $\beta$ to obtain}
\begin{equation}\label{eq:truncated-sublevel-volume}
 |E_{t+\delta}\cap\{\psi_0\ge-K\}|
 \le C\delta^{-d\beta}K^{d\beta}V(t)^{\beta/p'}.
\end{equation}
{The auxiliary Monge-Amp\`ere construction in the proof of
\cite[Proposition~2.1]{DK2014}, applied on $X$ with a smooth positive
density, provides a smooth $\chi$-psh test function $q$.}

{
Set $k=d\beta$, $\tau=(1+\beta/p')/2>1$, and $\eta=(\tau-1)/k$.
For $r=\tau/\eta$, the quasi-psh potential $\psi_0$ satisfies
$B_r:=\int_X(-\psi_0)^r\chi^n<\infty$.
By \eqref{eq:truncated-sublevel-volume} and Chebyshev's inequality,
\[
\begin{aligned}
 V(t+\delta)
 &\le |E_{t+\delta}\cap\{\psi_0\ge-K\}|+|\{\psi_0<-K\}|\\
 &\le C\delta^{-k}K^kV(t)^{\beta/p'}+B_rK^{-r}.
\end{aligned}
\]
If $0<V(t)\le1$, take $K=V(t)^{-\eta}$. Since
$\beta/p'-k\eta=\eta r=\tau$ and $0<\delta\le1$, we obtain
\begin{equation}\label{eq:volume-iteration}
 V(t+\delta)\le A\delta^{-k}V(t)^\tau.
\end{equation}

For $0<L\le L_*$, set $\gamma=\eta/(1+\eta)$ and $t_0=L^\gamma$.
By Markov's inequality, $V(t_0)\le L/t_0=L^{1/(1+\eta)}$.
Let $D_0=(2A)^{1/k}$ and define $t_{j+1}=t_j+D_0V(t_j)^\eta$.
Choose $L_*$ small enough that $V(t_0)\le1$ and $D_0V(t_0)^\eta\le1$.
Whenever $V(t_j)>0$, it follows from \eqref{eq:volume-iteration} that
\[
 V(t_{j+1})\le A D_0^{-k} V(t_j)^{\tau-k\eta}
              =\tfrac12 V(t_j).
\]
If $V(t_j)=0$, then $t_{j+1}=t_j$ and the same bound holds.
Hence $V(t_j)\le 2^{-j}V(t_0)$ and
\[
 t_\infty:=\lim_{j\to\infty}t_j
 \le t_0+\frac{D_0 V(t_0)^\eta}{1-2^{-\eta}}
 \le \left(1+\frac{D_0}{1-2^{-\eta}}\right)L^\gamma.
\]
Since $E_{t_\infty}\subset E_{t_j}$ for every $j$, the open set
$E_{t_\infty}$ has zero volume and is empty. Thus
$\sup_X(v-u)_+\le t_\infty\le C L^\gamma$.
The case $L=0$ follows from continuity.
}
\end{proof}

{Suppose $\psi_0\le u,v\le0$ and both densities have $L^p$ norm at most $H$.
If $\|u-v\|_{L^1(X,\chi^n)}\le L_*$, applying Proposition~\ref{prop:smooth-stability}
to $(u,v)$ and $(v,u)$ yields}
\[
 \|u-v\|_{L^\infty(X)}
 \le C\|u-v\|_{L^1(X,\chi^n)}^\gamma.
\]

\subsection{An auxiliary family}

The preceding stability estimate will allow us to apply
Lemma~\ref{lem:gradient-transfer} to neighboring solutions along a path of
densities. To start, we need a uniform weighted gradient bound at one
endpoint. We obtain this bound for an auxiliary family $q_t$ by comparing
its solutions at nearby values of $t$.

For $0<t\le\bar t$, we write $\theta_t=\theta+t\chi$ and fix
\[
 C_\theta=\max\left\{1,\sup_{0<t\le\bar t}\sup_X|F(\theta_t)|\right\},\qquad \Lambda>C_\theta.
\]
We consider the equation
\begin{equation}\label{eq:exponential-family}
 F(\theta_t+\ii\partial\bar\partial q_t)=\Lambda e^{q_t}.
\end{equation}
Since $[\theta_t]$ contains a smooth strict representative, Lemma~\ref{lem:proper-equation}, applied with $H(s)=\log\Lambda+s$, implies that this equation admits a unique smooth admissible solution $q_t$.

We denote $m_t=\sup_Xq_t$ and choose $x_t\in X$ such that $q_t(x_t)=m_t$. Since $\ii\partial\bar\partial q_t(x_t)\le0$, ellipticity implies that
\[
 \Lambda e^{m_t}\le F(\theta_t)(x_t)\le C_\theta<\Lambda.
\]
Thus $q_t<0$. We set $r_t=\Lambda e^{q_t}$. Integrating the equation and using Stokes' theorem, we obtain
\begin{equation}\label{eq:exponential-mass}
 0<r_t\le C_\theta,\qquad \int_Xr_t\,\chi^n=Q_t.
\end{equation}
It follows that
\[
 \log\frac{Q_t}{\Lambda\int_X\chi^n}\le m_t\le\log\frac{C_\theta}{\Lambda}<0.
\]
As $Q_t\to Q>0$, we have $m_t\ge-C_0$ uniformly in $t$.

We now normalize $q_t$ by setting $u_t=q_t-m_t$. Then $\sup_Xu_t=0$ and
\[
 \theta_t+\ii\partial\bar\partial u_t=\theta_t+\ii\partial\bar\partial q_t\in C_f(\chi),\qquad F(\theta_t+\ii\partial\bar\partial u_t)=r_t.
\]
For any fixed $l>n$, the bound $r_t\le C_\theta$ implies that $\|r_t\|_{\Psi_l,\nu}\le C_\theta\|1\|_{\Psi_l,\nu}$. Applying Proposition~\ref{prop:relative-bound} to $u_t$, we obtain $u_t\ge V_t-C_1$, with $C_1$ independent of $t$. By the definition of the envelope in \eqref{eq:admissible-envelope}, we have $V_t\ge V_\theta$. Moreover, since $w_{1/2}$ from \eqref{eq:small-pole-bound} is a competitor for $V_\theta$,
\[
 V_\theta\ge w_{1/2}=\frac{\psi}{2}+\frac{b_{1/2}-\sup_Xb_{1/2}}{2}\ge\frac{\psi}{2}-\frac12\operatorname{osc}_Xb_{1/2}.
\]
Combining these inequalities with $m_t\ge-C_0$, we conclude that
\[
 q_t=u_t+m_t\ge\frac{\psi}{2}-C,\qquad C=C_0+C_1+\frac12\operatorname{osc}_Xb_{1/2}.
\]

To study the dependence on $t$, we write $A_t=\theta_t+\ii\partial\bar\partial q_t$ and $G=F^{1/d}$. The equation becomes
\[
 G(A_t)=h(q_t),\qquad h(s)=\Lambda^{1/d}e^{s/d},\qquad h'(s)=\frac{h(s)}d,\quad h''(s)=\frac{h(s)}{d^2}.
\]
We denote the derivative of $G$ at $A_t$ in the direction $B$ by $DG_{A_t}[B]$, and set
\[
 L_tw=DG_{A_t}[\ii\partial\bar\partial w],\qquad W_t=DG_{A_t}[\chi].
\]
In a $\chi$-unitary frame diagonalizing $A_t$, we have
\[
 L_tw=\sum_i r_iw_{i\bar i},\qquad W_t=\sum_i r_i,\qquad r_i=\frac{\partial G}{\partial\lambda_i}\bigl(\lambda_\chi(A_t)\bigr)>0.
\]
For each fixed $t>0$, the linearized operator $L_t-h'(q_t)$ is elliptic and has a strictly negative zeroth-order coefficient, so it is invertible on $X$. By the implicit function theorem, $q_t$ depends smoothly on $t$.

Set $p_t=\partial_t q_t$. Since
$\partial_t A_t=\chi+\ii\partial\bar\partial p_t$, differentiating
$G(A_t)=h(q_t)$ with respect to $t$ yields
\[
 h'(q_t)p_t
 =DG_{A_t}[\chi+\ii\partial\bar\partial p_t]
 =W_t+L_t p_t .
\]
Thus $(L_t-h')p_t=-W_t$, where $h'=h'(q_t)=h(q_t)/d>0$.
At a minimum point of $p_t$, we have $L_t p_t\ge0$, so
$h'p_t\ge W_t>0$ there. Hence $p_t>0$ on $X$.

Using the bound $q_t\ge\psi/2-C$, subtract $2C+2d+1$ from $\psi$
and keep the same notation. With $T=1-\psi$ we then have
\[
 z_t:=q_t-\psi\ge d+\frac{T}{2}\qquad\text{on }X\setminus Z .
\]
{Set} $B=\theta+\ii\partial\bar\partial\psi-\delta\chi\in\overline C_f$.
By concavity of $G$, $DG_{A_t}[A_t-B]\le G(A_t)-G(B)\le h(q_t)$.
Since $h'(q_t)=h(q_t)/d$, it follows that
\[
 \begin{aligned}
 (L_t-h'(q_t))z_t
 &=DG_{A_t}[A_t-B]-(\delta+t)W_t-h'(q_t)z_t\\
 &\le h(q_t)\left(1-\frac{z_t}{d}\right)-(\delta+t)W_t
 \le-(\delta+t)W_t .
 \end{aligned}
\]

Since $(L_t-h'(q_t))p_t=-W_t$, we have
$(L_t-h'(q_t))(p_t-z_t/(\delta+t))\ge0$ on $X\setminus Z$.
For fixed $t>0$, $p_t-z_t/(\delta+t)$ tends to $-\infty$ along $Z$.
By the maximum principle, $p_t\le z_t/(\delta+t)$.
Integrating the bound for $p_s$ from $s=\varepsilon$ to $s=t$ leads to the second inequality below.
\begin{equation}\label{eq:parameter-velocity}
 0<p_t\le\frac{q_t-\psi}{\delta+t}\le\frac{1-\psi}{\delta},
 \qquad
 0\le q_t-q_\varepsilon\le\frac{t-\varepsilon}{\delta}(1-\psi)
 \quad(0<\varepsilon\le t\le\bar t).
\end{equation}

On the range of \(q_t\), we have \(h(q_t)\le C_\theta^{1/d}\) and
\(h',h''\ge0\), so Proposition~\ref{prop:weighted-hessian} and
Lemma~\ref{lem:gradient-transfer} apply. The constant \(C_2\) in
\eqref{eq:gradient-transfer} is independent of the exponent interval.
Choose \(A>0\) so that \(C_2e^{-A/8}<1/4\).
Fix
$I=[a_0,a_0+A]$, {let $\eta_I$ be the constant in}
Lemma~\ref{lem:gradient-transfer}, and choose
\[
 0<\eta\le\min\{\eta_I,\bar t/(2\delta),1\},\qquad
 c=\delta\eta,\qquad t_0=2c,\qquad
 a(t)=a_0+\frac{A}{t_0}(t_0-t)\quad (0<t\le t_0).
\]

For \(0<\varepsilon\le t_0/2\) and \(x\in X\setminus Z\), set
\[
 s_\varepsilon(x)=\varepsilon+\frac{c}{T(x)},
 \qquad v_\varepsilon(x)=q_{s_\varepsilon(x)}(x).
\]

Write \(s=s_\varepsilon(x)\). Since \(T\ge1\) and \(c=t_0/2\), we have
\(\varepsilon<s\le t_0\). Evaluating \eqref{eq:parameter-velocity} at \(x\), we find
\[
 0\le v_\varepsilon(x)-q_\varepsilon(x)
 \le\frac{T(x)}{\delta}(s-\varepsilon)=\frac c\delta=\eta .
\]
Thus $v_\varepsilon$ satisfies the assumption in Lemma~\ref{lem:gradient-transfer}. To estimate its gradient, we differentiate with respect to $x$ and obtain
\[
 \partial v_\varepsilon(x)=(\partial_x q_s)(x)+p_s(x)\frac{c}{T(x)^2}\,\partial\psi(x).
\]
We write
\[
 N_t=\sup_{X\setminus Z}e^{a(t)(\psi-q_t)/2}|\partial q_t|_\chi.
\]
Since $q_s-\psi\ge T/2$ and $a(\varepsilon)-a(s)=A(s-\varepsilon)/t_0=A/(2T)$, we have
\[
 \frac{e^{a(\varepsilon)(\psi-q_s)/2}}{e^{a(s)(\psi-q_s)/2}}
 =\exp\!\left[-\frac{A}{4T}(q_s-\psi)\right]\le e^{-A/8}.
\]
Hence the first term satisfies
\[
 e^{a(\varepsilon)(\psi-q_s)/2}|\partial_x q_s|_\chi\le e^{-A/8}N_s.
\]
For the second term, we use $p_s\le T/\delta$, $q_s-\psi\ge T/2$, and \eqref{eq:barrier-gradient} to obtain
\[
 e^{a(\varepsilon)(\psi-q_s)/2}\frac{c\,p_s}{T^2}|\partial\psi|_\chi
 \le\frac{C\eta}{T}\exp\!\left[-\left(\frac{a(\varepsilon)}4-\frac{\kappa}{2}\right)T\right]\le C,
\]
where the last inequality follows from $a(\varepsilon)\ge a_0>2\kappa$ and $T\ge1$.

We now fix $0<\varepsilon_*\le t_0/2$ and set $M_{\varepsilon_*}=\sup_{\varepsilon_*\le t\le t_0}N_t$. Since $q_t$ depends smoothly on $t$, we have $M_{\varepsilon_*}<\infty$. Moreover, for $\varepsilon_*\le\varepsilon\le t_0/2$, we have $\varepsilon_*\le s_\varepsilon(x)\le t_0$, so $N_{s_\varepsilon(x)}\le M_{\varepsilon_*}$. Applying
Lemma~\ref{lem:gradient-transfer} to \(q_\varepsilon\) and \(v_\varepsilon\), we obtain
\[
 N_\varepsilon\le C_*+kM_{\varepsilon_*},
 \qquad k=C_2e^{-A/8}<1/4.
\]
By smoothness, \(N_t\le C_*'\) on \([t_0/2,t_0]\).
Taking the supremum over the two intervals yields
\(M_{\varepsilon_*}\le\max\{C_*',C_* /(1-k)\}\), independently of \(\varepsilon_*\).

It follows that $K_{a(t)}(q_t)=N_t^2$ is uniformly bounded for
$0<t\le t_0$. Proposition~\ref{prop:weighted-hessian} then bounds
$H_{a(t)}(q_t)$ on the same interval. Set $A_*=a_0+A$. Since $q_t<0$,
$\psi\le0$, and $a(t)\le A_*$, we have
$e^{A_*\psi}\le e^{a(t)(\psi-q_t)}$. Combining these bounds with
smoothness on $[t_0,\bar t]$ yields
\begin{equation}\label{eq:seed-derivatives}
 \sup_{0<t\le\bar t}\sup_{X\setminus Z}
 \ee^{A_*\psi}\bigl(|\partial q_t|_\chi^2
                 +|\ii\partial\bar\partial q_t|_\chi\bigr)\le C .
\end{equation}

\subsection{A path of densities}\label{subsec:density-path}

Let $g>0$ be the smooth density in Theorem~\ref{thm:regularity}.
Starting with the bound \eqref{eq:seed-derivatives} for $q_t$, we now
keep $\theta_t$ fixed and vary the density from $r_t$ in \eqref{eq:exponential-mass} to
$Q_tg$. Proposition~\ref{prop:smooth-stability} will allow us to choose
the same finite subdivision for every $t$. Set
\[
 \rho_{t,s}:=(1-s)r_t+sQ_tg,\qquad 0\le s\le1 .
\]
For each $s\in[0,1]$, the density $\rho_{t,s}$ is smooth and positive.
By \eqref{eq:exponential-mass} and the normalization of $g$,
\[
 \int_X \rho_{t,s}\,\chi^n=(1-s)Q_t+sQ_t=Q_t .
\]
Since $[\theta_t]$ admits a smooth representative in $C_f(\chi)$ for $t>0$,
Proposition~\ref{prop:strict-equation} provides a unique $U_{t,s}\in C^\infty(X)$
satisfying
\[
 \begin{aligned}
 F(\theta_t+\ii\partial\bar\partial U_{t,s})&=\rho_{t,s},\quad
 \theta_t+\ii\partial\bar\partial U_{t,s}\in C_f(\chi),
 \qquad \sup_X U_{t,s}=0 .
 \end{aligned}
\]
At $s=0$ the normalized potential $q_t-\sup_Xq_t$ has density $r_t$,
and at $s=1$ the prescribed density is $Q_tg$. It follows from uniqueness that
\[
 U_{t,0}=q_t-\sup_Xq_t,\qquad U_{t,1}=u_{t,g}.
\]

By \eqref{eq:parameter-velocity}, $q_t$ decreases as $t\searrow0$.
Since $q_t\ge\psi/2-C$, its limit $q_0$ is finite on
$X\setminus Z$. Thus $r_t\to r_0:=\Lambda e^{q_0}>0$ a.e.
By \eqref{eq:exponential-mass} and dominated convergence,
\[
 r_t\longrightarrow r_0\ \text{in }L^p(X,\chi^n)\quad(1\le p<\infty),
 \qquad \int_Xr_0\,\chi^n=Q .
\]
All $\rho_{t,s}$ have a common upper bound, so
Proposition~\ref{prop:relative-bound} provides uniform relative bounds for $U_{t,s}$.

{
We claim that
\begin{equation}\label{eq:path-l1}
 \sup_{\substack{0<t\le\bar t\\ s,s'\in[0,1],\,|s-s'|\le h}}
 \|U_{t,s}-U_{t,s'}\|_{L^1}\longrightarrow0
 \qquad (h\searrow0).
\end{equation}
By compactness of $[0,\bar t]\times[0,1]$, it suffices to show that
$\|U_{t_j,s_j}-U_{t_j,s'_j}\|_{L^1}\to0$ whenever $s_j,s'_j\to s$
and either $t_j\searrow0$ or $t_j\to t_*>0$.

If $t_j\searrow0$, both $\rho_{t_j,s_j}$ and $\rho_{t_j,s'_j}$
converge in $L^1$ to $(1-s)r_0+sQg$, which is positive a.e.
 Proposition~\ref{prop:actual-family} therefore yields the
same $L^1$ limit for the two sequences of normalized solutions.
{
If $t_j\to t_*>0$, the normalized solutions depend continuously on
$(t,s)$ by the implicit function theorem applied modulo constants.
Hence $U_{t_j,s_j}$ and $U_{t_j,s'_j}$ both converge uniformly to
$U_{t_*,s}$ when $s_j,s'_j\to s$.
}
}

By \eqref{eq:small-pole-bound} with $\tau=1/2$, there is a constant
$B>0$, independent of $t,s$, such that
\[
 \psi_0:=\psi-B\le \tfrac12\psi-B\le U_{t,s}\le0 .
\]
On $X\setminus Z$ we have
\[
 \theta_t+\ii\partial\bar\partial\psi_0-\delta\chi
 =(\theta+\ii\partial\bar\partial\psi-\delta\chi)+t\chi
 \in\overline C_f(\chi).
\]
Since the densities are uniformly bounded, Proposition~\ref{prop:smooth-stability}
applied to $(U_{t,s},U_{t,s'})$ in both orders yields
\[
 \|U_{t,s}-U_{t,s'}\|_{L^\infty(X)}
 \le C\|U_{t,s}-U_{t,s'}\|_{L^1(X,\chi^n)}^\gamma
 \quad\text{if }\ \|U_{t,s}-U_{t,s'}\|_{L^1(X,\chi^n)}\le L_* .
\]
Together with \eqref{eq:path-l1}, the preceding estimate implies that, for any $\eta>0$, we can choose $N$ independently of $t$ such that
\begin{equation}\label{eq:path-subdivision}
 \sup_{0<t\le\bar t}\max_{1\le j\le N}\|U_{t,s_j}-U_{t,s_{j-1}}\|_{L^\infty(X)}\le\eta,\qquad s_j=j/N.
\end{equation}
Since $\log r_t=\log\Lambda+q_t$, we obtain from \eqref{eq:seed-derivatives} that
\[
 |\partial\log r_t|_\chi^2+|\ii\partial\bar\partial\log r_t|_\chi\le Ce^{-A_*\psi}.
\]

For $\alpha=\log r_t$, $\beta=\log(Q_tg)$ and
$\omega=(1-s)\ee^\alpha/\rho_{t,s}$,
\begin{align*}
 \partial\log\rho_{t,s}
 &=\omega\partial\alpha+(1-\omega)\partial\beta,\\
 \ii\partial\bar\partial\log\rho_{t,s}
 &=\omega\ii\partial\bar\partial\alpha
 +(1-\omega)\ii\partial\bar\partial\beta\\
 &\quad+\omega(1-\omega)
       \ii\partial(\alpha-\beta)\wedge\bar\partial(\alpha-\beta).
\end{align*}
Since $0\le\omega\le1$ and the spatial derivatives of $\beta$ are those
of $\log g$, these formulas bound
$|\partial\log\rho_{t,s}|^2+|\ii\partial\bar\partial\log\rho_{t,s}|$
by $C e^{-A_*\psi}$. Differentiating
$\rho_{t,s}^{1/d}=\exp((\log\rho_{t,s})/d)$ and using the common upper
bound for $\rho_{t,s}$, we obtain
\[
 |\ii\partial\bar\partial(\rho_{t,s}^{1/d})|
 +\frac{|\partial(\rho_{t,s}^{1/d})|^2}{\rho_{t,s}^{1/d}}
 \le C\ee^{-A_*\psi}.
\]

{Let $a_0$ be the constant in Proposition~\ref{prop:weighted-hessian} for $h(x)=\rho_{t,s}^{1/d}$. Fix $a\ge\max\{2A_*,a_0\}$ and set $I=[a,a+1]$.} At $s=0$ the quantity $K_a(U_{t,0})$ is uniformly bounded, since \eqref{eq:small-pole-bound} and \eqref{eq:seed-derivatives} give
\[
 \ee^{a(\psi-U_{t,0})}|\partial U_{t,0}|^2
 \le C\ee^{(a/2-A_*)\psi}
            \ee^{A_*\psi}|\partial q_t|^2.
\]

Choose the subdivision in \eqref{eq:path-subdivision} with $\eta=\eta_I$.
For each $1\le j\le N$, we have $\|U_{t,s_j}-U_{t,s_{j-1}}\|_{L^\infty(X)}\le\eta_I$. Applying Lemma~\ref{lem:gradient-transfer} to $U_{t,s_j}$ and $U_{t,s_{j-1}}$, we obtain
\[
 \sqrt{K_a(U_{t,s_j})}\le C_1(I)+C_2\sqrt{K_a(U_{t,s_{j-1}})}.
\]
We then have
\[
 \sqrt{K_a(U_{t,s_N})}\le C_2^N\sqrt{K_a(U_{t,s_0})}+C_1(I)\sum_{l=0}^{N-1}C_2^l.
\]
Since \eqref{eq:seed-derivatives} bounds the initial term and $N$ is independent of $t$, this proves a uniform bound for $K_a(u_{t,g})$. Here we used $s_N=1$ and $U_{t,1}=u_{t,g}$. By \eqref{eq:weighted-hessian}, we also obtain a uniform bound for $H_a(u_{t,g})$.

Writing $u=u_{t,g}$, we have $e^{a\psi}\le e^{a(\psi-u)}$, since $u\le0$. Together with $\ii\partial\bar\partial u=\eta_u-\Omega$ and $|\eta_u|_\chi\le n\lambda$, this implies that
\[
 e^{a\psi}|\partial u|_\chi^2\le K_a(u),\qquad e^{a\psi}|\ii\partial\bar\partial u|_\chi\le nH_a(u)+C.
\]
We therefore conclude that
\begin{equation}\label{eq:final-derivatives}
 \sup_{0<t\le\bar t}\sup_{X\setminus Z}\ee^{a\psi}\bigl(|\partial u_{t,g}|^2+|\ii\partial\bar\partial u_{t,g}|\bigr)\le C.
\end{equation}

\begin{proof}[Proof of Theorem~\ref{thm:regularity}]
We fix a compact subset $K\subset X\setminus Z$. By \eqref{eq:small-pole-bound} with $\tau=1/2$ and \eqref{eq:final-derivatives}, we have
\[
 \frac12\min_K\psi-C\le u_{t,g}\le0,\qquad |\partial u_{t,g}|_\chi^2+|\ii\partial\bar\partial u_{t,g}|_\chi\le Ce^{-a\min_K\psi}\quad\text{on }K.
\]
We denote the eigenvalues of $\theta_t+\ii\partial\bar\partial u_{t,g}$ with respect to $\chi$ by $\lambda_t$. Since $Q_t\to Q>0$ and $g>0$, we have
\[
 F(\lambda_t)=Q_tg\ge c,\qquad c=\frac12Q\min_Xg>0,
\]
for all sufficiently small $t$. The preceding Hessian bound and the identity $F=0$ on $\partial C_f$ then imply that
\[
 \lambda_t(x)\in\{\lambda\in\overline C_f:|\lambda|\le R_K,\ F(\lambda)\ge c\}\Subset C_f,\qquad x\in K,
\]
where $R_K$ is independent of $t$. Thus the equation is uniformly elliptic on $K$: for $G=F^{1/d}$, we have
\[
 0<c_K\le\partial_iG(\lambda_t(x))\le C_K,\qquad x\in K,\quad 1\le i\le n.
\]
These bounds also hold on a compact neighborhood of $K$ in $X\setminus Z$,
with constants depending on that neighborhood. The interior $C^{2,\alpha}$
estimate \cite{TWWY15} and the Schauder estimates for the differentiated
equation therefore give
\[
 \|u_{t,g}\|_{C^m(K)}\le C_{K,m},\qquad m\ge0 .
\]

Since $Q_tg\to Qg>0$ uniformly and $V_t-C\le u_{t,g}\le V_t$ by
Proposition~\ref{prop:relative-bound}, Proposition~\ref{prop:actual-family}
yields the same normalized $L^1(X)$ limit $u_g$ along every sequence
$t_j\searrow0$, with
\[
 \mu_f(\theta,u_g)=Qg\,\chi^n,\qquad V_\theta-C\le u_g\le V_\theta .
\]
Thus $u_{t,g}\to u_g$ in $L^1(X)$ as $t\searrow0$. The local $C^m$ bounds
give smooth precompactness on $X\setminus Z$, and every subsequential
limit there agrees with $u_g$.
 
Hence $u_{t,g}\to u_g$ in $C^\infty_{\rm loc}(X\setminus Z)$, where $Z=\operatorname{Null}_f(\alpha,\beta)$ by Theorem~\ref{thm:null}.
\end{proof}

\end{document}